\documentclass[12pt]{article}
\usepackage{amssymb,amsfonts,amsmath,amsthm,comment,dsfont,xcolor,tikz,tikz-cd}
\usepackage[colorlinks, citecolor=purple,linkcolor=blue,urlcolor=teal]{hyperref}
\usepackage[lite,alphabetic]{amsrefs}
\usepackage{cleveref}
\usepackage[a4paper, portrait, margin=25mm]{geometry}

\newcommand{\1}{\mathds{1}}
\newcommand{\0}{\mathds{O}}

\newcommand{\Z}{\mathbb{Z}}
\newcommand{\Q}{\mathbb{Q}}
\newcommand{\R}{\mathbb{R}}
\newcommand{\I}{\mathbb{I}}

\newcommand{\N}{\mathbb{N}}

\newcommand{\Ba}{\mathrm{B}}

\newcommand{\8}{\infty}

\newcommand{\spa}{\mathrm{span}}

\newcommand{\ext}{\mathrm{ext}}
\newcommand{\vlt}{\mathrm{vlt}}

\newcommand{\supp}{\mathrm{supp}}

\newcommand{\Co}{\mathcal{C}}

\newcommand{\Fo}{\mathcal{F}}

\newcommand{\Po}{\mathcal{P}}

\newcommand{\Ho}{\mathcal{H}}

\newcommand{\todo}{\xrightarrow[]{\mathrm{do}}}
\newcommand{\too}{\xrightarrow[]{\mathrm{o}}}
\newcommand{\touo}{\xrightarrow[]{\mathrm{uo}}}

\newcommand{\toru}{\xrightarrow[]{\mathrm{ru}}}

\newtheoremstyle{mytheorem}
  {12pt}   
  {4pt}   
  {\itshape} 
  {}       
  {\bfseries} 
  {.}      
  {.5em}   
  {}       

\newtheoremstyle{mydefinition}
  {12pt}        
  {8pt}       
  {\normalfont} 
  {}           
  {\bfseries}   
  {.}           
  {.5em}        
  {}

\newcounter{dummy} \numberwithin{dummy}{section}
\theoremstyle{mytheorem}
\newtheorem{theorem}[dummy]{Theorem}
\newtheorem{lemma}[dummy]{Lemma}

\newtheorem{proposition}[dummy]{Proposition}
\newtheorem{corollary}[dummy]{Corollary}
\newtheorem{question}[dummy]{Question}
\theoremstyle{mydefinition}
\newtheorem{remark}[dummy]{Remark}
\newtheorem{example}[dummy]{Example}
\usepackage{scalerel}
\DeclareMathOperator*{\bigplus}{\scalerel*{+}{\sum}}

\begin{document}

\title{Extensions of operators between ordered vector spaces}

\author{Eugene Bilokopytov\thanks{Department of Mathematical and Statistical Sciences,  University of Alberta, Edmonton, Canada and Instituto de Ciencias Matem\'aticas, Madrid, Spain (\texttt{bilokopy@ualberta.ca, ievgen.bilokopytov@icmat.es}).}}

\maketitle

\begin{abstract}
We prove analogues of several classical theorems about operators between vector lattices in the context of directed Archimedean ordered vector spaces. In particular, we show that every order continuous operator between pervasive Archimedean ordered vector spaces can be extended to an order continuous operator between their Dedekind completions. Additionally, we introduce the universal completion of a directed Archimedean ordered vector space, and show that every supremum-preserving operator between directed Archimedean ordered vector spaces can be extended to an order continuous homomorphism between their universal completions. We also establish some universal properties for the relative uniform completion.\medskip


MSC2020 46A19, 46A40, 46B42, 46E25, 54A20.
\end{abstract}

{\footnotesize\tableofcontents}\medskip

\section{Introduction}

This article is a follow-up to \cite{brtgww}, where various results on order convergence, which are well-known in the vector lattice theory, have been generalized to the level of directed Archimedean ordered vector spaces. In particular, it was established in that paper that the class of so-called pervasive ordered vector spaces is best suited for work with order convergence, since the latter behaves on such spaces in a way similar to its behavior on vector lattices.

The topic which was left unexplored in the aforementioned paper is operator theory, in particular various notions of continuity. We now fill this gap, and study order continuous operators and their extensions. We are primarily inspired by Veksler's theorem, which essentially allows extending order continuous operators ``by continuity'', despite the fact that the vector lattices in question might not be endowed with topologies.\medskip

The starting point of our investigation was a comparison of Veksler's and Riesz-Kantorovich theorems, as well as Fremlin's theorem on the extension of an order continuous homomorphism from an order dense sublattice. It turned out that the main steps in the proofs of these results are very similar. Namely, using a certain formula, an additive map is constructed, which then gives rise to the operator which we were looking for (an extension in the case of Veksler's and Fremlin's theorems, and the positive part in the case of Riesz-Kantorovich theorem). It was possible to replicate this method in some cases when the domain of the operator is not necessarily a lattice.

Another source of inspiration the result from \cite{aerz} which says that any supremum-preserving operator between directed Archimedean ordered vector spaces can be extended to an order continuous homomorphism between their Dedekind completions. Note that while Riesz-Kantorovich-Veksler-style results will be applicable to a relatively broad class of operators on a relatively narrow class of spaces, the aforementioned theorem guarantees extensions of a very narrow class of operators but on a broader class of spaces.

Using the opportunity we also present a number of extension results in other contexts, including extension of some special classes of operators on ordered vector spaces to their Riesz completion, and extension of positive operators to the relative uniform completion. We also consider order convergence and its variants on the space of operators, and study relative uniform and universal completions of directed Archimedean ordered vector spaces.\medskip

Let us describe the content of the paper. In Section \ref{secoc} we present basic results about order convergence and its relatives on ordered vector spaces, and on spaces of operators between them; we also present some elementary extension results. In particular, in Subsection \ref{sboc} we examine various notions of order convergence, continuity and completeness for ordered vector spaces. In Subsection \ref{sbop} we focus on one of these convergences on the space of operators and show that sets of operators with certain properties are closed in various senses (Propositions \ref{preserve} and \ref{bsc}). In Subsection \ref{supco} we explore ``finitary'' versions of the aforementioned result from \cite{aerz}. In Subsection \ref{s2} we continue studying the spaces of operators by going through the proof of the classical Riesz-Kantorovich theorem, keeping track of the precise assumptions needed for each step. We also enrich the configuration by considering operators defined on a subspace of an ordered vector space, but preserving the boundedness inherited from the entire space. Subsection \ref{sbob} is dedicated to relative uniform convergence and extension of some operators ``by continuity'' with respect to this convergence.

The core of the article is Section \ref{secperv} focused on extensions of operators on pervasive ordered vector spaces. We present some basic facts about these spaces in Subsection \ref{sbperv}. The main results of the paper are contained in Subsection \ref{s4}. Namely, in Theorem \ref{veksler0} we show that the collection of order continuous operators from a pervasive ordered vector space into a Dedekind complete vector lattice is itself a Dedekind complete vector lattice, and in Theorem \ref{veksler} we also establish isomorphism between that vector lattice and the vector lattice of order continuous operators from the Dedekind completion of the domain. We introduce unbounded order convergence in Subsection \ref{s6} and prove an extension result (Theorem \ref{uocont2}) for uo-continuous operators on pervasive spaces.

In Section \ref{secot} is dedicated to completions of directed Archimedean vector spaces. In particular, in Subsection \ref{sbuc} we introduce the concept of the universal completion of a directed Archimedean ordered vector space, and show in Theorem \ref{supcont4} that a supremum-preserving operator can be extended to an order continuous homomorphism between the universal completions. In Subsection \ref{sbru} we generalize the material about the relative uniform completion of a vector lattice studied in \cite{bt} to the level of directed Archimedean ordered vector spaces. Subsection \ref{cat} contains the categorical overview of the obtained results, while in Subsection \ref{sublim} we present a ``purely categorical'' construction of the relative uniform completion. Some counterexamples announced in the main text are not presented right away in order to avoid an interruption of the flow. They are instead ``parked'' in Section \ref{sbce}.

\section{Convergences and spaces of operators}\label{secoc}
\subsection{Order convergence and continuity}\label{sboc}

Let $P$ be a poset. If $q\in P$ and $Q\subset P$, we write $Q\uparrow q$ if $Q$ is directed upwards and $q=\bigvee Q$. The meaning of the notation $Q\downarrow q$ is analogous. We denote $Q^{\Uparrow}:=\left\{p\in P,~ p\ge Q\right\}$, i.e. the collection of all upper bounds of $Q$. If $R\subset P$, we further denote $Q^{\Uparrow_{R}}:=Q^{\Uparrow}\cap R$. The meaning of the notations $Q^{\Downarrow}$ and $Q^{\Downarrow_{R}}$ is similar. We say that $Q$ is \emph{supremum-dense} in $P$ if $p=\bigvee p^{\Downarrow_{Q}}$, for every $p\in P$, and define \emph{infimum-denseness} analogously. If $Q\subset R\subset P$, then the suprema of $Q$ within $P$ and $R$ may differ; we denote them by $\bigvee_{P}Q$ and $\bigvee_{R}Q$, in case of possible ambiguity. The meaning of $\bigwedge_{P}Q$ and $\bigwedge_{R}Q$ is analogous.

If $Q$ is such that $\left[p,q\right]\subset Q$ whenever $p,q\in Q$ and $p\le q$, then we call $Q$ \emph{full}; we call it a \emph{lower} set if $p\le q\in Q$ implies $p\in Q$. \emph{Upper} sets are defined analogously. Upper sets and lower sets are full.

If $\varphi:P\to R$ is a map between posets $P$ and $R$, we call it \emph{supremum-continuous}, if it preserves all existing suprema in $P$, i.e. if $Q\subset P$ and $q\in P$ are such that $q=\bigvee Q$, then $\varphi\left(q\right)=\bigvee\varphi\left(Q\right)$. Every order isomorphism is supremum-continuous. More generally, if $Q\subset P$ is supremum-dense and infimum-dense, then the embedding of $Q$ into $P$ is supremum-continuous (see \cite[Proposition 2.6]{brtgww}).

We also call $\varphi$ \emph{directed supremum-continuous} if it preserves all existing suprema of upward directed sets in $P$, and \emph{bounded supremum-continuous} if it preserves all existing suprema of order bounded sets in $P$. Every bounded supremum-continuous map is directed supremum-continuous. If $p\le q$, taking $Q=\left\{p,q\right\}$ yields $\varphi\left(q\right)=\varphi\left(p\right)\vee \varphi\left(q\right)$, in other words $\varphi\left(p\right)\le \varphi\left(q\right)$. It follows that every directed supremum-continuous (including bounded supremum-continuous and supremum-continuous) map is order preserving. Note that even though every set which has a supremum is bounded from above, it might not be bounded from below. Hence, the concept of bounded supremum-continuity is formally weaker than that of supremum-continuity. In fact, Example \ref{bscnr} exhibits a bounded supremum-continuous functional on a directed Archimedean ordered vector space which fails to be supremum-continuous.\medskip

Throughout the subsection $E,F$ are ordered vector spaces. A linear subspace $H\subset E$ is supremum-dense if and only if it is infimum-dense. Recall that $E$ is a \emph{Dedekind complete vector lattice} if it is directed and every order bounded set has a supremum in $E$ (in particular, every set of two elements). If $E$ is directed and Archimedean, then it has a \emph{Dedekind completion} $E^{\delta}$, which is the essentially unique Dedekind complete vector lattice which contains $E$ as a supremum-dense subspace (see \cite[Section 2.1]{kv}). It is easy to see that if $H$ is supremum-dense in $E$, then it is also supremum dense in $E^{\delta}$, and so $H^{\delta}=E^{\delta}$.

If $E,F$ are ordered vector spaces $T\in L\left(E,F\right)$ is supremum-continuous if and only if it preserves the sets whose infimum is zero.\footnote{In the literature these operators are often called complete Riesz homomorphisms, however such a name would be inconsistent with the nomenclature that we are developing in this paper.} In the case when $E,F$ are vector lattices this condition is also equivalent to the fact that $T$ is an order continuous vector lattice homomorphism (see \cite[Theorem 5.2]{erz}). The following result was presented in \cite{aerz}.

\begin{theorem}\label{supcont1}
Assume that $E$ and $F$ are directed and Archimedean. For an operator $T:E\to F$ the following conditions are equivalent:
\item[(i)] $T$ is supremum-continuous;
\item[(ii)] $T$ is positive and $T\left(A^{\Uparrow}\right)^{\Downarrow}\subset \left(TA\right)^{\Uparrow\Downarrow}$, for every $A\subset E$;
\item[(iii)] $T\left(A^{\Uparrow\Downarrow}\right)\subset \left(TA\right)^{\Uparrow\Downarrow}$, for every $A\subset E$;
\item[(iv)] There is a supremum-continuous operator $T^{\delta}:E^{\delta}\to F^{\delta}$ which extends $T$;
\item[(v)] There are Archimedean vector lattices $\widetilde{E},\widetilde{F}$ and a supremum-continuous operator $\widetilde{T}:\widetilde{E}\to \widetilde{F}$ such that $E$ embeds as a subspace of $\widetilde{E}$ in a supremum-continuous way, $F$ embeds as a subspace of $\widetilde{F}$, and $\widetilde{T}$ extends $T$.\medskip

Moreover, the extension in (iv) is unique even among positive operators. It is injective if and only if $T$ is an order embedding.
\end{theorem}

\begin{remark}\label{ries}
Note that $T$ is positive if and only if $T\left(A^{\Uparrow}\right)^{\Downarrow}\supset \left(TA\right)^{\Uparrow\Downarrow}$, for every $A\subset E$, and if and only if $T\left(A^{\Uparrow}\right)^{\Downarrow}\supset \left(TA\right)^{\Uparrow\Downarrow}$, for $A=\left\{0_{E}\right\}$. Hence, the condition (ii) in Theorem \ref{supcont1} can be restated as $T\left(A^{\Uparrow}\right)^{\Downarrow}= \left(TA\right)^{\Uparrow\Downarrow}$, for every $A\subset E$.
\qed\end{remark}

Several notions of convergence will appear throughout the text of the article. For the general theory we refer to \cite{bctv} and \cite{btw}, however we only need a few concepts and results. Namely, a map is called \emph{continuous} if it preserves limits of nets. A linear operator is \emph{bounded}, if it preserves bounded sets (see the definition in \cite[Section 6]{bctv}). Every continuous linear map is bounded. An \emph{adherence} of a set is the collection of all existing limits of nets in this set. A set $A$ is \emph{closed}, if it is equal to its adherence. The \emph{closure} of a set is the intersection of all closed sets which contain the given one. The closure can also be obtained by iterating the adherence. If $\eta$ is a convergence structure on a set $X$ and $A\subset X$, we denote the adherence and closure of $A$ with respect to $\eta$ by $\overline{A}_{\eta}^{1}$ and $\overline{A}_{\eta}$, respectively.\medskip

We say that a net $\left(e_{\alpha}\right)_{\alpha\in A}\subset E$ converges \emph{in order} to $e\in E$ (denoted $e_{\alpha}\too e$) if there is $P\subset E$ with $\bigwedge P=0_{E}$ which \emph{dominates the tails} of $\left(e_{\alpha}-e\right)_{\alpha\in A}$, i.e. for each $p\in P$ there is $\alpha_{p}$ such that $-p\le e_{\alpha}-e\le p$, whenever $\alpha\ge \alpha_{p}$. If $P$ was also downward directed, we say that $\left(e_{\alpha}\right)_{\alpha\in A}\subset E$ \emph{directed order} converges to $e\in E$, and denote this fact by $e_{\alpha}\todo e$. Clearly, directed order convergence is \emph{stronger} than order convergence, in the sense that every directed order convergent net order converges to the same limit. Every monotone net (directed) order converges to its extremum (if it exists, see \cite[Corollary 2.4]{brtgww}).

If $E$ is directed and Archimedean, then a set is order bounded if and only if it is bounded with respect to order convergence (see \cite[Proposition 2.23(iii)]{brtgww}). It follows that order continuous operators are order bounded. Also, order convergence on $E$ is the restriction of order convergence on $E^{\delta}$ (see \cite[Proposition 2.19]{brtgww}).

It is easy to see that $T\in L\left(E,F\right)$ is directed supremum-continuous if and only if it preserves the nets which decrease to zero. We will use the term \emph{monotone-continuous} for these operators.

\begin{proposition}\label{4.10}
Let $T\in L\left(E,F\right)$. Then:
\item[(i)] If $T$ is supremum-continuous, then it is order continuous.
\item[(ii)] $T$ is monotone continuous if and only if it is positive and directed order continuous if and only if it is positive and continuous from directed order convergence on $E$ into order convergence on $F$.
\item[(iii)] If $T$ is positive and order continuous, then it is directed order continuous.
\end{proposition}
\begin{proof}
(i) follows from \cite[Proposition 2.2]{brtgww}. (iii) follows from (ii).

(ii): If $T$ is monotone continuous, it is order preserving, hence positive. By \cite[Proposition 2.15]{brtgww} it is also directed order continuous. Since directed order convergence is stronger than order convergence, it follows that if $T$ is positive and directed order continuous then it is continuous from directed order convergence on $E$ into order convergence on $F$. Finally, assume that the latter condition is satisfied and $T$ is positive. If $P\downarrow 0_{E}$, then $P\todo 0_{E}$, hence $TP\too 0_{F}$. As $TP$ is directed downwards, we conclude that $TP\downarrow 0_{F}$.
\end{proof}

\begin{remark}
Positivity is an essential assumption in part (iii) of Proposition \ref{4.10}, while a positive directed order continuous operator needs not be order continuous. There is no analogue of the characterization (ii) in Proposition \ref{4.10} for non-positive operators. Namely, there is an ordered vector space $E$ and a regular linear functional $\nu:E\to\R$ which fails to be order continuous, but such that whenever $P\downarrow 0_{E}$ it then follows that $\nu\left(p\right)\to 0$. The corresponding counterexamples will be presented in \cite{bb}.
\qed\end{remark}

We will say that $E$ is \emph{monotone complete} if every directed bounded from above set in $E$ has a supremum. We will call a subspace $H\subset E$ \emph{monotone closed} if it is closed in $E$ with respect to taking suprema of increasing nets. It is easy to see that it is enough to check both of these conditions for directed subsets of $E_{+}$ and $H_{+}$, respectively. It is clear that a directed order closed subspace must be monotone closed. The converse is false even for full directed subspaces (a counterexample will be presented in \cite{bb}). Recall that $H$ is \emph{majorizing} in $E$ if for every $e\in E$ there is $h\in H$ such that $e\le h$.

\begin{proposition}\label{monotone}
Assume that $E$ is monotone complete and $H\subset E$ is a linear subspace. Then:
\item[(i)] If $H$ is monotone closed then it is monotone complete in its own right, and the embedding of $H$ into $E$ is monotone continuous.
\item[(ii)] If $H$ is full, the same conclusion as in (i).
\item[(iii)] If additionally $H$ is majorizing, then the converse to (i) holds.
\end{proposition}
\begin{proof}
(i): Let $P\subset H_{+}$ be directed and bounded from above. Since $E$ is monotone complete, there is $p\in E$ such that $p=\bigvee_{E} P$. As $H$ is monotone closed, it follows that $p\in H$ and consequently $p=\bigvee_{H} P$. We conclude that $H$ is monotone complete. If we started with the assumption that $P$ has a supremum in $H$, this supremum would have to be $p$. This shows that the embedding of $H$ into $E$ is monotone continuous.\medskip

(ii): Let $P$ and $p$ be as above. Let $h\in H$ be such that $h\ge P$ and let $q\in P$. It follows that $p\in\left[q,h\right]\subset H$ and consequently $p=\bigvee_{H} P$. The second part of the argument is the same as in (i).\medskip

(iii): Assume that $P\subset H_{+}$ is directed and let $p\in E$ be such that $p=\bigvee_{E}P$. As $H$ is majorizing it follows that $P$ is bounded from above in $H$, and since $H$ is monotone complete there is $q\in H$ such that $q=\bigvee_{H}P$. Because the embedding of $H$ into $E$ is monotone continuous, it follows that $p=q\in H$. Thus, $H$ is monotone closed.
\end{proof}

\begin{corollary}\label{monoton}
$E$ is monotone complete if and only if it is monotone closed in $E^{\delta}$.
\end{corollary}

\subsection{Directed order convergence on $L\left(E,F\right)$}\label{sbop}

Throughout the subsection $E$ and $F$ are ordered vector spaces, $E$ directed. We will study directed order convergence on $L\left(E,F\right)$. In order to do so we will first look at the last part of Riesz-Kantorovich theorem \cite[Theorem 1.59]{at} which deals with the supremum of a monotone net of operators. We will refer to semigroup homomorphisms of cones as \emph{additive maps}. Recall the following Kantorovich's Extension Lemma.

\begin{proposition}[Lemma 1.26, \cite{at}]\label{kext}
Assume that $E,F$ are directed and Archimedean. Any additive map $\tau:E_{+}\to F_{+}$ uniquely extends to a positive operator $T_{\tau}:E\to F$.
\end{proposition}

Let us denote the space of all maps from $E_{+}$ into $F$ by $F^{E_{+}}$. It is easy to see that any existing supremum of a set in $F^{E_{+}}$ is evaluated pointwise. Let $R:L\left(E,F\right)\to F^{E_{+}}$ be the restriction map. Clearly, $R$ is a linear order embedding (i.e. $T\in L\left(E,F\right)_{+}$ if and only if $R\left(T\right)\ge 0_{F^{E_{+}}}$).

\begin{proposition}\label{rkm0}
The image of $R$ is monotone closed in $F^{E_{+}}$.
\end{proposition}
\begin{proof}
Let $\left(T_{\alpha}\right)_{\alpha\in A}\subset L\left(E,F\right)_{+}$ be an increasing net, such that $\tau:=\bigvee\limits_{\alpha\in A} R\left(T_{\alpha}\right)$ exists in $F^{E_{+}}$. Then, $\tau\left(e\right)=\bigvee\limits_{\alpha}T_{\alpha}e$, for every $e\in E_{+}$. Since $T_{\alpha}$'s are positive, it follows that $\tau\left(e\right)\in F_{+}$, for every $e\in E_{+}$.

Fix $e,g\in E_{+}$. On the one hand, for every $\alpha$ we have $T_{\alpha}\left(e+g\right)=T_{\alpha}e+T_{\alpha}g\le \tau\left(e\right) + \tau\left(g\right)$, and so $\tau\left(e+g\right)\le \tau\left(e\right) + \tau\left(g\right)$. On the other hand, for every $\alpha,\beta$ there is $\gamma\ge\alpha,\beta$, so that $T_{\alpha}e+T_{\beta}g\le T_{\gamma}e+T_{\gamma}g=T_{\gamma}\left(e+g\right)\le \tau\left(e+g\right)$, hence $\tau\left(e\right) + \tau\left(g\right)\le \tau\left(e+g\right)$. Thus, $\tau$ is additive, and so by Proposition \ref{kext} it extends to a positive operator $T:E\to F$. It follows that $R\left(T_{\alpha}\right)\uparrow R\left(T\right)$, which shows that the image of $R$ is monotone closed.
\end{proof}

We can now recover the last part of Riesz-Kantorovich theorem (and even more).

\begin{corollary}\label{rkm}
Assume that $F$ is monotone complete. Then:
\item[(i)] Let $\left(T_{\alpha}\right)_{\alpha\in A}\subset L\left(E,F\right)$ be a bounded from above increasing net. Then, it has a supremum $T\in L\left(E,F\right)$ such that $Te=\bigvee\limits_{\alpha\in A}T_{\alpha}e$, for every $e\in E_{+}$.
\item[(ii)] If $\left(T_{\alpha}\right)_{\alpha\in A}\subset L\left(E,F\right)$ is directed order null, then $T_{\alpha}e\todo 0_{F}$, for every $e\in E$.
\end{corollary}
\begin{proof}
It is easy to see that if $F$ is monotone complete, then so is $F^{E_{+}}$. Now (i) follows from applying part (i) of Proposition \ref{monotone} to the image of $R$.

(ii): Let $\left(S_{\beta}\right)_{\beta\in B}\subset L\left(E,F\right)$ with $S_{\beta}\downarrow 0_{L\left(E,F\right)}$ dominate the tails of $\left(T_{\alpha}\right)_{\alpha\in A}$. For $e\in e_{+}$ it is not hard to check that $\left(S_{\beta}e\right)_{\beta\in B}$ dominates the tails of $\left(T_{\alpha}e\right)_{\alpha\in A}$, hence $T_{\alpha}e\todo 0_{F}$. Since $E$ is directed we in fact have $T_{\alpha}e\todo 0_{F}$, for every $e\in E$.
\end{proof}

\begin{remark}
The converse to part (ii) in Corollary \ref{rkm} is false. Assume that $E=L_{2}$ and $F=\R$. Then, $e_{n}\in E\simeq L\left(E,F\right)$ defined by $e\left(t\right):=\sin 2n\pi t$ is an orthogonal, hence weakly-null order bounded sequence. However, it does not converge to $\0$ almost everywhere, and so it is not (directed) order-null.
\qed\end{remark}

A net $\left(e_{\alpha}\right)_{\alpha\in A}\subset E$ converges \emph{relatively uniformly (ru)} to $e\in E$ (we denote this by $e_{\alpha}\toru e$) if there is $h\in E_{+}$ such that for every $\varepsilon>0$ there is $\alpha_{0}$ such that $-\varepsilon h\le e-e_{\alpha}\le\varepsilon h$ as soon as $\alpha\ge\alpha_{0}$. If $E$ is Archimedean, then $\mathrm{ru}$ is stronger than $\mathrm{do}$. Relative uniform convergence will be discussed in greater detail in Subsection \ref{sbob}.

\begin{proposition}\label{preserve}
Let $P\subset E_{+}$, and let $L_{P}$ be the set of operators $T\in L\left(E,F\right)_{+}$ such that $\bigwedge TP=0_{F}$. Then:
\item[(i)] $L_{P}$ is a lower set in $L\left(E,F\right)_{+}$.
\item[(ii)] If $P$ is bounded from above then $L_{P}$ is $\mathrm{ru}$-closed in $L\left(E,F\right)_{+}$.
\item[(iii)] If $P$ is bounded from above and $F$ is monotone complete then $L_{P}$ is $\mathrm{do}$-closed in $L\left(E,F\right)_{+}$.
\end{proposition}
\begin{proof}
(i): Assume that $S,T\in L\left(E,F\right)_{+}$ are such that $T\le S\in L_{P}$. Then, $\bigwedge SP=0_{F}$, and for every $p\in P$ we have that $0_{F}\le Tp\le Sp$. We conclude that $\bigwedge TP=0_{F}$, and so $T\in L_{P}$.\medskip

The proof of (ii) is analogous to that of (iii), but simpler, hence we only prove the latter.

(iii): Let $\left(T_{\alpha}\right)_{\alpha\in A}\subset L_{P}$ be directed order convergent to $T \in L\left(E,F\right)_{+}$. Let $f\in F$ be such that $f\le TP$. Our goal is to show that $f\le 0_{F}$. Let $e\in E$ be an upper bound for $P$.

Let $\left(S_{\beta}\right)_{\beta\in B}\subset L\left(E,F\right)$ with $S_{\beta}\downarrow 0_{L\left(E,F\right)}$ dominate the tails of $\left(T_{\alpha}-T\right)_{\alpha\in A}$. Fix $\beta\in B$ and find $\alpha\in A$ such that $-S_{\beta}\le T_{\alpha}-T\le S_{\beta}$. For every $p\in P\subset E_{+}$ we have $$f\le Tp\le \left(T_{\alpha}+S_{\beta}\right)p=T_{\alpha}p+S_{\beta}p\le T_{\alpha}p+S_{\beta}e,$$ hence $f-S_{\beta}e\le T_{\alpha}p$. As $T_{\alpha}\in L_{P}$, it follows that $f-S_{\beta}e\le 0_{F}$, so that $f\le S_{\beta}e$. As $\beta$ was arbitrary and according to part (i) of Corollary \ref{rkm} $S_{\beta}e\downarrow 0_{F}$, we conclude that $f\le 0_{F}$, as required.
\end{proof}

It is easy to see that $T\in L\left(E,F\right)$ is bounded supremum-continuous if and only if it preserves sets which are bounded from above and whose infimum is zero if and only if it preserves suprema of sets contained in $E_{+}$. Denote the collection of (directed) order continuous operators by $L_{oc}\left(E,F\right)$ ($L_{doc}\left(E,F\right)$), which is a subspace of $L\left(E,F\right)$. Recall that by part (ii) of Proposition \ref{4.10} we have that $L_{doc}\left(E,F\right)_{+}$ is precisely the set of monotone continuous operators.

A linear subspace $H\subset E$ is full if and only if $0_{E}\le e\le h\in H$ implies $e\in H$, in other words $H_{+}$ is a lower subset of $E_{+}$.

\begin{proposition}\label{bsc0}
$L_{doc}\left(E,F\right)$ is a full subspace of $L\left(E,F\right)$. The collection of (bounded) supremum-continuous operators is a lower subset of $L\left(E,F\right)_{+}$.
\end{proposition}
\begin{proof}
It follows immediately from part (i) of Proposition \ref{preserve} that the collections of (bounded / monotone) continuous operators are lower subsets of $L\left(E,F\right)_{+}$. Since the collection of monotone continuous operators equals $L_{doc}\left(E,F\right)_{+}$, we conclude that $L_{doc}\left(E,F\right)$ is a full subspace of $L\left(E,F\right)$.
\end{proof}

The following result may be viewed as a version of Ogasawara theorem \cite[Theorem 1.57]{ab} for ordered vector spaces.

\begin{proposition}\label{bsc}
Assume that $F$ is monotone complete. Then $L_{oc}\left(E,F\right)$, and the collections of monotone continuous, as well as bounded supremum-continuous operators are directed order closed in $L\left(E,F\right)$, while $L_{doc}\left(E,F\right)$ is monotone closed in $L\left(E,F\right)$.
\end{proposition}
\begin{proof}
Let $\left(T_{\alpha}\right)_{\alpha\in A}\subset L\left(E,F\right)$ be a net of bounded supremum-continuous operators which directed order converges to $T$. First, since $L\left(E,F\right)_{+}$ is directed order closed in $L\left(E,F\right)$ (see \cite[Proposition 2.3]{brtgww}), it follows that $T$ is positive. Moreover, it is bounded supremum-continuous, according to part (iii) of Proposition \ref{preserve}. Hence the collections of bounded supremum continuous operators is directed order closed in $L\left(E,F\right)$. Similarly, the collection of monotone continuous operators is also directed order closed in $L\left(E,F\right)$, hence monotone closed. As the latter set is equal to $L_{doc}\left(E,F\right)_{+}$, we conclude that $L_{doc}\left(E,F\right)$ is monotone closed in $L\left(E,F\right)$.\medskip

It is left to show that $L_{oc}\left(E,F\right)$ is directed order closed in $L\left(E,F\right)$. Suppose that $\left(T_{\alpha}\right)_{\alpha\in A}\subset L_{oc}\left(E,F\right)$ directed order converges to $T \in L\left(E,F\right)$. Let $\left(e_{\gamma}\right)_{\gamma\in\Gamma}\subset E$ be order-null, witnessed by $P\subset E_{+}$. In particular, without loss of generality we may assume that the net is order bounded, hence contained in $\left[-h,h\right]$, $h\in E_{+}$.

Let $\left(S_{\beta}\right)_{\beta\in B}\subset L\left(E,F\right)$ with $S_{\beta}\downarrow 0_{L\left(E,F\right)}$ dominate the tails of $\left(T_{\alpha}-T\right)_{\alpha\in A}$. Fix $\beta\in B$ and find $\alpha_{\beta}\in A$ such that $-S_{\beta}\le T_{\alpha_{\beta}}-T\le S_{\beta}$. For every $\gamma\in\Gamma$ we have $-S_{\beta}h\le \left(T_{\alpha_{\beta}}-T\right)e_{\gamma}\le S_{\beta}h$. Also, $T_{\alpha_{\beta}}e_{\gamma}\too 0_{F}$, witnessed by $Q_{\beta}\subset F_{+}$. For every $q\in Q_{\beta}$ we see that $-q\le T_{\alpha_{\beta}}e_{\gamma}\le q$, eventually, hence $-q-S_{\beta}h\le Te_{\gamma}\le q+S_{\beta}h$, eventually.

Let $Q:=\bigcup\limits_{\beta\in B}\left(Q_{\beta}+S_{\beta}h\right)\subset F_{+}$. First, let us show that $\bigwedge Q=0_{F}$. Assume that $f\le Q$, Then, for every $\beta\in B$ we have $f\le Q_{\beta}+S_{\beta}h$, and since $\bigwedge Q_{\beta}=0_{F}$, we get $f\le S_{\beta}h$. As $\beta$ was arbitrary, $S_{\beta}\downarrow 0_{L\left(E,F\right)}$, and $F$ is monotone complete, by part (i) of Corollary \ref{rkm} it follows that $S_{\beta}h\downarrow 0_{F}$, and so $f\le 0_{F}$. Next, the argument above shows that $Q$ dominates the tails of $\left(Te_{\gamma}\right)$. We conclude that $Te_{\gamma}\too 0_{F}$, and so $T$ is order continuous.
\end{proof}

\begin{remark}
Note that there is no analogue of Proposition \ref{bsc} for supremum-continuous operators. Indeed, in Example \ref{bscnoc} we will construct a sequence of supremum continuous functionals which relatively uniformly (hence directed order) converges to an operator which is necessarily bounded supremum-continuous, but not order continuous, and hence not supremum-continuous. This also shows that $L_{doc}\left(E,F\right)$ is not necessarily relatively uniformly closed in $L\left(E,F\right)$.

Monotone completeness is an essential assumption in part (ii) of Corollary \ref{rkm} and Proposition \ref{bsc}, even for vector lattices (see the discussion after the remark). On the other hand,  by part (ii) of Proposition \ref{preserve} the assumption of monotone completeness of $F$ is not needed to establish that the first three sets Proposition \ref{bsc} are closed with respect to relative uniform convergence.
\qed\end{remark}

It was proven in \cite{elliot} that the Riesz-Kantorovich formula for the positive part of an operator may fail in the case when the target space is not a Dedekind complete vector lattice. Namely, there is a compact space $K$ and an operator in $T\in L\left(L_{1},\Co\left(K\right)\right)$ which has a positive part which is not given by the Riesz-Kantorovich formula. Here we show that the Riesz-Kantorovich formula for an increasing net of operators may also fail. Namely, we construct a (non-closed) sublattice $F$ of $L_{1}$, and an increasing sequence of order continuous homomorphisms in $L\left(L_{1},F\right)$, whose supremum exists, but is not evaluated pointwise on $E_{+}$.

\begin{example}\label{mo}
Let $H:=\left\{\left(f_{n}\right)_{n\in\N}\subset L_{1},~ \exists \lim\limits_{n\to\8} f_{n}\right\}$. For $h\in H$ let $Rh$ be the limit in the definition of $H$. We view $H\oplus L_{1}\subset L_{1}^{\N_{\8}}$ and define $F:=\left\{\left(h,f\right)\in H\oplus L_{1},~ f-Rh\in L_{\8}\right\}$. For $n\in\N_{\8}$ let $P_{n}:F\to E:=L_{1}$ be the coordinate projection, and let $P:F\to H$ be the projection away from the last coordinate. Note that $L_{1}^{\N_{\8}}$ is lattice isometric to $L_{1}$.

For $0\le m<n\le \8$ let $T_{m}^{n}e:=\left(\0,...\0, e,...,e,\0,...\right)$, where the copies of $e$ are in the positions from $m+1$ to $n$, and $\0$ is on every other position. It is easy to see that all these operators are homomorphisms from $E$ into $F$. Since $E$ is order continuous, every positive operator on it order continuous, hence all the operators involved are supremum continuous.

A sequence $\left(T_{m}^{n}\right)_{n\in\N}$ is increasing with $T_{m}^{n}\le T_{m}^{\8}$, for every $m<n<\8$, but $T_{m}^{n}\1\not\uparrow T_{m}^{\8}\1$. Let us show that nevertheless we have $T_{m}^{n}\uparrow T_{m}^{\8}$, for every $m\in\N_{0}$. Since $T_{m}^{n}=T_{0}^{n}-T_{0}^{m}$, for every $0\le m<n\le \8$, it is enough to consider the case when $m=0$.

Assume that $S\in L\left(E,F\right)$ is such that $S\ge T_{0}^{n}$, for every $n\in\N$. In particular, $S\in L\left(E,F\right)_{+}$. For every $e\in E_{+}$ and $n\in\N$ we have $P_{n}Se\ge P_{n}T_{0}^{n}e=e$, and since $PSe\in H$, it follows that $e\le P_{n}Se\to RPSe$. Therefore, $RPS\ge\I_{E}$. Note that in order to prove $S\ge T_{0}^{\8}$, it is enough to show that $P_{\8}S\ge \I_{E}$.

We have that $U:=\left(RP-P_{\8}\right)S$ is a regular (hence continuous) operator on $E$ whose image is contained in $L_{\8}$. Its graph is closed in $E\times E$, hence in $E\times L_{\8}$, and so by the Closed Graph theorem, this operator is continuous from $E$ into $L_{\8}$. Therefore, for every $e\in E_{+}$ we have that $U\left[\0,e\right]$ is a norm bounded, hence an order bounded subset of $L_{\8}$, and so $U^{+}e=\bigvee U\left[\0,e\right]\in L_{\8}$ (here we are using Dedekind completeness of $L_{\8}$ and the Riesz-Kantorovich formula, see \cite[Theorem 1.59]{at}, or Subsection \ref{s2} below).

Let $Q:=U - \left(RPS-\I_{E}\right)=\I_{E}-P_{\8}S$. Then $Q\le \I_{E}$, $Q\le U$ yields $Q^{+}\le \I_{E}$ and $Q^{+}\le U^{+}$. The first inequality implies that there is $\0\le g\in L_{\8}$ such that $Q^{+}$ is the multiplication with $g$ (see \cite[Example 2.67]{ab}). Let us show that $g=\0$. Assume that is not the case. Then, there is a measurable $A\subset \left[0,1\right]$ of positive measure and $\varepsilon>0$ such that $g\ge \varepsilon\1_{A}$. Take an unbounded $f\in L_{1}\left(A\right)$, extend it by zero to $\left[0,1\right]$. Then, $L_{\8}\ni U^{+}f\ge Q^{+}f=gf\ge \varepsilon\1_{A}f$, contradicting the fact that $f$ is unbounded on $A$. It follows that $Q^{+}=0_{L\left(E,E\right)}$, hence $Q\le 0_{L\left(E,E\right)}$, and so $P_{\8}S\ge \I_{E}$, as required.
\qed\end{example}

Even if $H,F$ are sublattices of $L_{1}$ the supremum of an increasing sequence of order continuous homomorphisms in $L\left(H,F\right)$ may not be evaluated pointwise on $E_{+}$, and may fail to be order continuous.

\begin{example}\label{mo2}
Let $V_{n}:H\to F$ be defined by $V_{n}h:=\left(P_{1}h,...,P_{n}h,\0,...\right)$ and let $V:H\to F$ be defined by $Vh:=\left(P_{1}h,...,P_{n}h,...,Rh\right)$. It is clear that $\left(V_{n}\right)_{n\in\N}$ is an increasing sequence of supremum continuous operators dominated by $V$, which is a homomorphism. Let us show that $V_{n}\uparrow V$. Assume that $W\in L\left(H,F\right)$ is such that $W\ge V_{n}$, for every $n\in\N$. In particular, $W$ is positive, and hence continuous. Let $J_{n}:E\to H$ be defined by $J_{n}e:=\left(\0,...,\0,e,...\right)$, where the copies of $e$ are in the positions starting from $n\in\N$. Then, $WJ_{n}\ge V_{m}J_{n}=T_{n-1}^{m}$, for every $m\ge n$, hence $WJ_{n}\ge T_{n-1}^{\8}$, according to Example \ref{mo}.

For $h=\left(h_{n}\right)_{n\in\N}$ in $H_{+}$ and $m\in\N$ let $h^{m}=\left(h_{1},h_{2},...,h_{m-1},\0,...\right)$ and $g^{m}=\left(\0,...,\0,h_{m}-Rh,h_{m+1}-Rh,...\right)$. Then $g^{m}\to 0_{H}$, and also $V_{m}h^{m}+T_{m}^{\8}Rh\to Vh$. For every $m\in\N$ we have $Wh=Wh^{m}+WJ_{m}Rh+Wg^{m}\ge V_{m}h^{m}+T_{m}^{\8}Rh+Wg^{m}\to Vh$.\medskip

Let us now prove that $V$ is not order continuous. Indeed, for $J_{n}\1\downarrow \0$ in $H$, but $VJ_{n}\1= \left(\0,...,\0,...,\1\right)$ (meaning $\1$ is in the $\8$-th position), hence $VJ_{n}\1\not\downarrow 0_{F}$.
\qed\end{example}

\begin{remark}
It would be desirable to upgrade Examples \ref{mo} and \ref{mo2} to make all the spaces involved Banach lattices.
\qed\end{remark}

\subsection{Extensions of vector pre-lattice homomorphisms}\label{supco}

Throughout the section $E$ and $F$ are directed Archimedean ordered vector spaces. Let us obtain versions of Theorem \ref{supcont1} where we only consider finite sets in the conditions (ii) and (iii). If $E,F$ are ordered vector spaces, and $T:E\to F$ is a linear operator, we will call $T$ a \emph{vector pre-lattice homomorphism}, or just \emph{pre-homomorphism}\footnote{These are called Riesz* homomorphisms in the literature, but one of the aims of our project is to streamline the terminology for the topic.} if $T\left(A^{\Uparrow\Downarrow}\right)\subset \left(TA\right)^{\Uparrow\Downarrow}$, for every finite $A\subset E$, and \emph{strict vector pre-lattice homomorphism}, or just \emph{strict pre-homomorphism}\footnote{Known as Riesz homomorphisms in the literature.} if $T$ is positive and $T\left(A^{\Uparrow}\right)^{\Downarrow}\subset \left(TA\right)^{\Uparrow\Downarrow}$, for every finite $A\subset E$. As in the Remark \ref{ries} note that $T$ is a strict pre-homomorphism if and only if $T\left(A^{\Uparrow}\right)^{\Downarrow}= \left(TA\right)^{\Uparrow\Downarrow}$, for every finite $\varnothing\ne A\subset E$. It follows from Theorem \ref{supcont1} (see also \cite[Proposition 2.3.16]{kv}) that every supremum continuous operator between directed Archimedean ordered vector spaces is a strict pre-homomorphism, while every strict pre-homomorphism is a pre-homomorphism by \cite[Proposition 2.3.18]{kv}. Recall that the \emph{Riesz completion} $E^{\rho}$ of $E$ is the vector lattice generated by $E$ in $E^{\delta}$. The following characterizations complement results in \cite{kv} and \cite{kvgs}.

\begin{proposition}\label{riesz}For $T\in L\left(E,F\right)$ the following conditions are equivalent:
\item[(i)] $T$ is a pre-homomorphism;
\item[(ii)] There is a vector lattice homomorphism $T^{\rho}:E^{\rho}\to F^{\rho}$ which extends $T$;
\item[(iii)] There are Archimedean vector lattices $\widetilde{E},\widetilde{F}$ and a homomorphism $\widetilde{T}:\widetilde{E}\to \widetilde{F}$ such that $E$ embeds as a subspace of $\widetilde{E}$ in a supremum-continuous way, $F$ embeds as a subspace of $\widetilde{F}$, and $\widetilde{T}$ extends $T$.\medskip

Moreover, the extension in (ii) is unique. It is injective if and only if $T$ is an order embedding.
\end{proposition}
\begin{proof}
(i)$\Rightarrow$(ii) and the uniqueness clause follow from \cite[Theorem 2.3.12]{kv}. (ii)$\Rightarrow$(iii) is trivial, and the last claim are proven analogously to the corresponding part of the proof of Theorem \ref{supcont1} in \cite{aerz}.\medskip

(iii)$\Rightarrow$(i): Let $J:E\to \widetilde{E}$ be the embedding, which we assume to be supremum-continuous. Let $\varnothing\ne A\subset E$ be finite. It follows from Theorem \ref{supcont1} that $J\left(A^{\Uparrow\Downarrow}\right)\subset \left(JA\right)^{\Uparrow\Downarrow}=\left(\bigvee_{\widetilde{E}} JA\right)^{\Downarrow}$. Since $\widetilde{T}$ is a homomorphism and $T=\widetilde{T}J$, we get $$T\left(A^{\Uparrow\Downarrow}\right)=\widetilde{T}J\left(A^{\Uparrow\Downarrow}\right)\subset \widetilde{T}\left(\bigvee\nolimits_{\widetilde{E}} JA\right)^{\Downarrow}\subset \left(\widetilde{T}\bigvee\nolimits_{\widetilde{E}} JA\right)^{\Downarrow}=\left(\bigvee\nolimits_{\widetilde{F}} \widetilde{T}JA\right)^{\Downarrow}=\left(\bigvee\nolimits_{\widetilde{F}} TA\right)^{\Downarrow}=\left(TA\right)^{\Uparrow\Downarrow},$$
where we view $T$ as a map into $\widetilde{F}$. Since the desired property survives when the co-domain is restricted (see \cite[Lemma 3.8]{stennder}), the statement is proven.
\end{proof}

\begin{remark}
If $E$ was not directed, but $F$ was directed and Archimedean, the only pre-homomorphism from $E$ into $F$ is zero. Indeed, if $e\notin E_{+}-E_{+}$, then for $A=\left\{0_{E},e\right\}$ we have $A^{\Uparrow}=\varnothing$, and so $A^{\Uparrow\Downarrow}=E$. Hence, $TE=T\left(A^{\Uparrow\Downarrow}\right)\subset \left(TA\right)^{\Uparrow\Downarrow}=\left\{0_{F},Te\right\}^{\Uparrow\Downarrow}$. For $g\in E$ and $f\in \left\{0_{F},Te\right\}^{\Uparrow}$ we have $nTg\in TE \le f$, implying $Tg\le 0_{F}$, but then also $T\left(-g\right)\le 0_{F}$. We conclude that $TE=\left\{0_{F}\right\}$.
\qed\end{remark}

\begin{proposition}\label{riesz2}Let $E$ and $F$ be directed Archimedean ordered vector spaces and let $T\in L\left(E,F\right)_{+}$. Then, $T$ is a strict pre-homomorphism if and only if every positive extension $\widetilde{T}: E^{\rho}\to F^{\delta}$ of $T$ is a homomorphism.
\end{proposition}
\begin{proof}
Necessity follows from \cite[Proposition 6]{kvgs}. Let us prove sufficiency. By Kantorovich extension theorem \cite[Theorem 1.32]{ab} there is a positive $\widetilde{T}: E^{\rho}\to F^{\delta}$ which extends $T$. By our assumption $\widetilde{T}$ is a homomorphism.

Assume that there is a finite $\varnothing\ne A\subset E$ such that $T\left(A^{\Uparrow}\right)^{\Downarrow}\not\subset \left(TA\right)^{\Uparrow\Downarrow}$. Viewing $T$ as a map into $F^{\delta}$, this means that $\left\{\bigwedge T\left(A^{\Uparrow}\right)\right\}^{\Downarrow}\not\subset \left(\bigvee TA\right)^{\Downarrow}$, in other words $f:=\bigwedge T\left(A^{\Uparrow}\right)\not\le\bigvee TA$. As $T\left(A^{\Uparrow}\right)\ge TA$, we in fact have $f>\bigvee TA$.

Let $a:=\bigvee_{E^{\rho}}A$. If $a\in E$, then $a=\bigvee_{E}A$, therefore $A^{\Uparrow}=\left\{a\right\}^{\Uparrow}$, hence $f=Ta$, and since $\widetilde{T}$ is a homomorphism we get $f=Ta=\widetilde{T}a=\bigvee \widetilde{T}A=\bigvee TA<f$, contradiction. Hence, $a\notin E$; let $E_{1}:=E+\R a$. Define $T_{1}: E_{1}\to F^{\delta}$ by $T_{1}\left(e+ra\right):=Te+rf$. This is a well defined linear map which extends $T$. We will show momentarily that $T_{1}$ is positive. Once that is done we can apply Kantorovich extension theorem to $T_{1}$ to get positive $\widetilde{T}_{1}: E^{\rho}\to F^{\delta}$. Since then $\widetilde{T}_{1}$ also extends $T$ and $a\in E^{\rho}$, by our assumption it is a homomorphism, and so $f=T_{1}a=\widetilde{T}_{1}a=\bigvee \widetilde{T}_{1}A=\bigvee TA < f$, contradiction.\medskip

It is left to prove that if $e+ra\ge 0_{E_{1}}$, then $T_{1}\left(e+ra\right)=Te+rf\ge 0_{F^{\delta}}$. If $r=0$, it follows that $e\in E_{+}$, and so $T_{1}\left(e+ra\right)=Te\ge 0_{F^{\delta}}$. Hence, without loss of generality we may assume that $r=\pm 1$, and so it is enough to show that $e\ge a$ implies $Te\ge f$, and $-e\le a$ implies $-Te\le f$. In the first case we have $e \in a^{\Uparrow_{E}}=A^{\Uparrow}$, which yields $Te\ge \bigwedge T\left(A^{\Uparrow}\right) =f$. In the second case we have $-e\in a^{\Downarrow_{E}}=A^{\Uparrow\Downarrow}$, and so $-Te\in T\left(A^{\Uparrow\Downarrow}\right)\subset \left(TA\right)^{\Uparrow\Downarrow}=\left(\bigvee TA\right)^{\Downarrow}$. Hence, $-Te\le \bigvee TA<f$.
\end{proof}

\begin{remark}
One can show that if $E,F,G$ are Archimedean ordered vector spaces, and $T:E\to F$ and $J:F\to G$ are linear operators, such that $J$ is an order embedding, and $JT$ is positive, or (bounded, monotone) supremum-continuous, or is a (strict) pre-homomorphism, then the latter property also holds for $T$.
\qed\end{remark}

Note that $T$ is a strict pre-homomorphism iff $\bigwedge T\left(\left\{e,f\right\}^{\Uparrow}-\left\{e,f\right\}\right)=0_{F}$, for every $e,f\in E$ (see \cite[Lemma 2.3.23]{kv}). Hence, we get an immediate consequence of part (i) of Proposition \ref{preserve}.

\begin{corollary}\label{lower}
The collection of all strict pre-homomorphisms is a lower subset of $L\left(E,F\right)_{+}$.
\end{corollary}

There is no analogue of Corollary \ref{lower} for pre-homomorphisms.

\begin{example}\label{janko}
Let $E$ be the subspace of $\Co\left[-1,1\right]$ consisting of all functions $f$ such that $f\left(-1\right)+f\left(1\right)=2f\left(0\right)$. Then, $E$ is supremum-dense in $\Co\left[-1,1\right]$, and so the corresponding inclusion map is supremum-continuous. Moreover, $E$ generates $\Co\left[-1,1\right]$ as a vector lattice, and so $\Co\left[-1,1\right]=E^{\rho}$. It follows that $\delta_{0}$ is a pre-homomorphism, as a composition of a vector lattice homomorphism and a supremum continuous map. Note that $\delta_{0}=\frac{1}{2}\left(\delta_{-1}+\delta_{1}\right)$, and so $0_{E'}\le \nu\le 4\delta_{0}$, where $\nu:=\delta_{-1}+2\delta_{1}$. If $\nu$ was a pre-homomorphism, in the light of the fact that $\Co\left[-1,1\right]=E^{\rho}$, we would have a vector lattice homomorphism $\nu^{\rho}:\Co\left[-1,1\right]\to\R$ which extends $\nu$. Such a functional is of the form $r\delta_{t}$, for $r\ge 0$ and $t\in \left[-1,1\right]$. We leave it to the reader to verify that such a functional cannot restrict to $\nu$ on $E$.
\qed\end{example}

Nevertheless, the set of pre-homomorphisms enjoys a certain closedness property.

\begin{proposition}
Let $\left(T_{\alpha}\right)_{\alpha\in A}\subset L\left(E,F\right)$ and $T\in L\left(E,F\right)$ be such that $T_{\alpha}e\too Te$, for every $e\in E$. If $T_{\alpha}$ is a pre-homomorphism, for every $\alpha$, then the same is true for $T$.
\end{proposition}
\begin{proof}
Let $\varnothing\ne A\subset E$ be finite; it is enough to show that $T\left(A^{\Uparrow\Downarrow}\right)\subset \left(TA\right)^{\Uparrow\Downarrow}$. Assume that $e\in A^{\Uparrow\Downarrow}$, so that $e\le A^{\Uparrow}$, and let $f\in \left(TA\right)^{\Uparrow}$, so that $f\ge TA$. Our goal is to show that $f\ge Te$.

Let $A=\left\{a_{1},...,a_{n}\right\}$. Assume that $Q_{k}\subset F_{+}$ witnesses $T_{\alpha}a_{k}\too Ta_{k}$, for $k=1,...,n$, and that $Q_{0}$ witnesses $T_{\alpha}e\too Te$. It is easy to see that $Q:=Q_{0}+...+Q_{n}$ witnesses $T_{\alpha}a_{k}\too Ta_{k}$, for every $k=1,...,n$, as well as $T_{\alpha}e\too Te$.

Fix $q\in Q$ and find $\alpha$ such that $-q\le T_{\alpha}a_{k}-Ta_{k}\le q$, for every $k=1,...,n$, as well as $-q\le T_{\alpha}e-Te\le q$. We have $f+q\ge Ta_{k}+q\ge T_{\alpha}a_{k}$, for every $k=1,...,n$, and so $f+q\in \left(T_{\alpha}A\right)^{\Uparrow}$. Hence, $T_{\alpha}e\in T_{\alpha}\left(A^{\Uparrow\Downarrow}\right)\subset \left(T_{\alpha}A\right)^{\Uparrow\Downarrow}$ implies $Te-q\le T_{\alpha}e\le f+q$. It follows that $Te\le f+2q$. Since $q$ was arbitrary and $\bigwedge Q=0_{F}$, we conclude that $f\ge Te$.
\end{proof}

We say that $E$ has the \emph{Riesz decomposition property (RDP)} if $\left[0_{E},e\right]+\left[0_{E},g\right]\supset \left[0_{E},e+g\right]$, for every $e,g\in E_{+}$. Note that in fact we have equality, since the converse inclusion is always satisfied. Clearly, the RDP passes down to full subspaces. The RDP is also equivalent to the fact that whenever $A,B\subset E$ are finite and such that $A\le B$, there is $e\in E$ such that $A\le e\le B$ (see \cite[Definition 1.52 and Theorem 1.54]{at}).

\begin{lemma}\label{rdpa}
If $E$ has the RDP, $\varnothing\ne A\subset E$ is finite and $e\in A^{\Uparrow}$, then $\bigwedge\left(A^{\Uparrow}\cap e^{\Downarrow} -A\right)=0_{E}$. If $\varnothing\ne B\subset E_{+}$ is finite, then $\bigwedge\left(B-B^{\Downarrow}\cap E_{+}\right)=0_{E}$.
\end{lemma}
\begin{proof}
Let $P$ be the set in question and let $f\in P^{\Downarrow}$. We will find $\left(e_{n}\right)_{n\in\N}\subset A^{\Uparrow}\cap e^{\Downarrow}$ such that $e_{0}=e$ and $e_{n}\le e_{n-1}-f$, for every $n\in\N$. Assume that $e_{n-1}$ is selected. Then, $e_{n-1}-A\subset P\ge f$, hence $e_{n-1}-f\ge A$, and so we can apply the property mentioned above to $A$ and $\left\{e_{n-1}-f,e\right\}$ to get $e_{n}\in E$ such that $e_{n-1}-f,e\ge e_{n}\ge A$, as required. Fix $a\in A$. For every $n\in\N$ we have $a\le e_{n}\le e-nf$. Archimedean property now yields $f\le 0_{E}$. We can now conclude that $\bigwedge P=0_{E}$. The second claim follows from the first one.
\end{proof}

\begin{proposition}\label{rdpb}
Assume that $E$ has the RDP. Then:
\item[(i)] Every bounded supremum continuous operator from $E$ into $F$ is a strict pre-homomorphism.
\item[(ii)] Every bounded supremum continuous functional on $E$ is supremum continuous.
\end{proposition}
\begin{proof}
(i): Let $T:E\to F$ be bounded supremum continuous. Let $e,f\in E$, and let $g\in \left\{e,f\right\}^{\Uparrow}$ and let $P:=\left\{e,f\right\}^{\Uparrow}\cap g^{\Downarrow} -\left\{e,f\right\}$, which is order bounded. It follows from Lemma \ref{rdpa} that $\bigwedge P=0_{E}$, hence $\bigwedge TP=0_{F}$, and so $\bigwedge T\left(\left\{e,f\right\}^{\Uparrow}-\left\{e,f\right\}\right)=0_{F}$. Since $e,f$ were arbitrary, we conclude that $T$ is a strict pre-homomorphism.\medskip

(ii): Let $\nu:E\to\R$ be linear and bounded supremum continuous. Assume that $Q\subset E$ and $\varepsilon>0$ are such that $\bigwedge Q=0_{E}$, but $\nu\left(q\right)\ge\varepsilon$, for every $q\in Q$. Fix $q_{0}\in Q$. For every $r\in Q$ according to Lemma \ref{rdpa} we have that $\bigwedge\left(\left\{q_{0},r\right\}-\left\{q_{0},r\right\}^{\Downarrow}\cap E_{+}\right)=0_{E}$. Since the set in the infimum is order bounded it follows that $\bigwedge \nu\left(\left\{q_{0},r\right\}-\left\{q_{0},r\right\}^{\Downarrow}\cap E_{+}\right)=0$, and so there is $p_{r}\in \left\{q_{0},r\right\}^{\Downarrow}\cap E_{+}\subset\left[0_{E},q_{0}\right]$ such that either $\nu\left(q_{0}-p_{r}\right)\le \frac{\varepsilon}{2}$ or $\nu\left(r-p_{r}\right)\le \frac{\varepsilon}{2}$. Since $\nu\left(q_{0}\right),\nu\left(r\right)\ge\varepsilon$, we see that $\nu\left(p_{r}\right)\ge \frac{\varepsilon}{2}$. On the other hand, since $\bigwedge Q=0_{E}$ the set $R:=\left\{p_{r}\right\}_{r\in Q}\subset\left[0_{E},q_{0}\right]$ satisfies $\bigwedge R=0_{E}$, hence $\bigwedge \nu\left(R\right)=0$, contradiction.
\end{proof}

Note that the RDP assumption in Proposition \ref{rdpb} is not superfluous, see Example \ref{bscnr}. It is also not possible to strengthen the conclusion of the result and show that $T$ is supremum continuous, for any $F$, see Example \ref{tong}.

\subsection{Revisiting Riesz-Kantorovich theorem}\label{s2}

Throughout the subsection $E$ is a directed Archimedean ordered vector space, $H\subset E$ is a directed linear subspace and $F$ is a Dedekind complete vector lattice. If $e,g\in E$ with $e\le g$ we denote $\left[e,g\right]_{H}:=\left[e,g\right]\cap H$. We say that a linear operator $T:H\to F$ is \emph{$E$-bounded} if $T\left[0_{E},e\right]_{H}$ is order bounded in $F$, for every $e\in E_{+}$ (cf. \cite[Definition 2.1]{buskes}). Note that if $H$ is majorizing in $E$, then $T$ is $E$-bounded if and only if it is order bounded. In general, if $T\in L_{ob}\left(E,F\right)$, then $\left.T\right|_{H}$ is $E$-bounded. Let $L^{E}_{ob}\left(H,F\right)$ be the collection of all $E$-bounded linear operators; it is easy to see that $L^{E}_{ob}\left(H,F\right)$ is a full linear subspace of $L_{ob}\left(H,F\right)$.

\begin{proposition}\label{addit0}
If $T\in L^{E}_{ob}\left(H,F\right)$, then:
\item[(i)] A map $\tau_{T}:E_{+}\to F_{+}$ defined by $\tau_{T}\left(e\right):=\bigvee T\left[0_{E},e\right]_{H}$ is superadditive, i.e. $\tau_{T}\left(e+g\right)\ge \tau_{T}\left(e\right)+\tau_{T}\left(g\right)$, for every $e,g\in E_{+}$.
\item[(ii)] $\tau_{T}\left(h\right)\ge Th$, for every $h\in H_{+}$. If $T$ is positive, then it agrees with $\tau_{T}$ on $H_{+}$.
\item[(iii)] For every $e\in E_{+}$ we have $\tau_{T}\left(e\right)=\bigvee \tau_{T}\left(\left[0_{E},e\right]_{H}\right)$.
\end{proposition}
\begin{proof}
(i): First, $\tau_{T}$ is well defined since $F$ is Dedekind complete, and $T\left[0_{E},e\right]_{H}$ is order bounded in $F$. Superadditivity follows immediately from $\left[0_{E},e\right]_{H}+\left[0_{E},g\right]_{H}\subset \left[0_{E},e+g\right]_{H}$, and properties of the supremum in a vector lattice.

(ii) and (iii) are straightforward.
\end{proof}

Recall that every additive $\tau:E_{+}\to F_{+}$ can be uniquely extended to a positive operator $T_{\tau}:E\to F$. If $G\subset E$ is a directed linear subspace, then $T_{\left.\tau\right|_{G}}=\left.T_{\tau}\right|_{G}$.

\begin{proposition}\label{addit}
Assume that $T\in L^{E}_{ob}\left(H,F\right)$. Then:
\item[(i)] If $\tau_{T}$ is subadditive, then $T_{\tau_{T}}$ is the smallest among the positive operators $S:E\to F$ such that $Sh\ge Th$, for every $h\in H_{+}$. If additionally $T$ is positive, then $T_{\tau_{T}}$ is the smallest extension of $T$ to $E$.
\item[(ii)] If $\tau_{T}$ is subadditive on $H_{+}$, then $T_{\left.\tau_{T}\right|_{H}}=T^{+}$ in $L^{E}_{ob}\left(H,F\right)$. In particular, if $H=E$, then $T_{\tau_{T}}=T^{+}$ in $L\left(E,F\right)$.
\item[(iii)] If both $\tau_{T},\tau_{-T}$ are subadditive, then $T_{\tau_{T}}-T_{\tau_{-T}}$ is a regular extension of $T$ to $E$.
\end{proposition}
\begin{proof}
(i): First, if $\tau_{T}$ is subadditive, then by part (i) of Proposition \ref{addit0} it is in fact additive and so $T_{\tau_{T}}$ exists by Proposition \ref{kext}. By part (ii) of Proposition \ref{addit0} it satisfies $T_{\tau_{T}}h=\tau_{T}\left(h\right)\ge Th$, for every $h\in H_{+}$.

Assume that $S:E\to F$ is positive and also satisfies $Sh\ge Th$, for every $h\in H_{+}$. Then, for every $e\in E_{+}$ and $h\in\left[0_{E},e\right]_{H}$ we have $Se\ge Sh\ge Th$, hence $Se\ge \bigvee T\left[0_{E},e\right]_{H}=\tau_{T}\left(e\right)=T_{\tau_{T}}e$.

The last claim follows from part (ii) of Proposition \ref{addit0}.\medskip

(ii): Fix $e\in E_{+}$. There is $f\in F_{+}$ such that $T\left[0_{E},e\right]_{H}\subset\left[-f,f\right]$. For every $h\in \left[0_{E},e\right]_{H}$ we have $0_{F}\le \tau\left(h\right)=\bigvee T\left[0_{E},h\right]_{H}\le \bigvee T\left[0_{E},e\right]_{H}\le f$. Hence, $T_{\left.\tau_{T}\right|_{H}}\left[0_{E},e\right]_{H}=\tau_{T}\left(\left[0_{E},e\right]_{H}\right)\subset\left[0_{F},f\right]$. Since $e$ was arbitrary, we conclude that $T_{\left.\tau_{T}\right|_{H}}\in L^{E}_{ob}\left(H,F\right)_{+}$. The fact that it is $T^{+}$ in $L^{E}_{ob}\left(H,F\right)$ is proven analogously to (i).\medskip

(iii): According to (ii) we have that $\left.T_{\tau_{T}}\right|_{H}=T^{+}$ in $L^{E}_{ob}\left(H,F\right)$, and so $\left.T_{\tau_{-T}}\right|_{H}=\left(-T\right)^{+}=T^{-}$ in $L^{E}_{ob}\left(H,F\right)$. In particular, $\left.\left(T_{\tau_{T}}-T_{\tau_{-T}}\right)\right|_{H}=T^{+}-T^{-}=T$.
\end{proof}

\begin{corollary}\label{rk0}
Let $G\subset L^{E}_{ob}\left(H,F\right)$ be a linear subspace such that for every $T\in G$ we have that $\tau_{T}$ is subadditive on $H_{+}$ and $T_{\left.\tau_{T}\right|_{H}}\in G$. Then, $G$ is a vector lattice in which the lattice operations are given via Riesz-Kantorovich formulas (\cite[Theorem 1.59]{at}).
\end{corollary}

Note that ``other'' Riesz-Kantorovich formulas can be deduced from the formula for the positive part.

The following result may be interpreted in the context of the so-called monogenic operators (see \cite{lp}).

\begin{lemma}\label{rk}
Let $G\subset L_{ob}\left(E,F\right)$ and $G_{H}\subset L^{E}_{ob}\left(H,F\right)$ be linear subspaces which satisfy the following properties:
\item[(i)] For every $S\in G$ we have $\left.S\right|_{H}\in G_{H}$;
\item[(ii)] If $S\in G$ vanishes on $H$, then $S=0_{L\left(E,F\right)}$;
\item[(iii)] For every $T\in G_{H}$ we have that $\tau_{T}$ is subadditive and $T_{\tau_{T}}\in G$.\medskip

Then, both $G$ and $G_{H}$ are vector lattices in which the lattice operations are given via Riesz-Kantorovich formulas (\cite[Theorem 1.59]{at}) and the restriction to $H$ is an order isomorphism from $G$ onto $G_{H}$.\medskip

If moreover, either $G$ is full in $L\left(E,F\right)$, or $G_{H}$ is full in $L\left(H,F\right)$, then both of them are Dedekind complete vector lattices, where the suprema of an increasing nets of operators are evaluated pointwise.
\end{lemma}
\begin{proof}
The fact that $G_{H}$ is a vector lattice with the specified lattice operations follows from Corollary \ref{rk0}. Let $R:G\to G_{H}$ be the operator of restriction to $H$. Clearly, $R$ preserves positive operators, and so is itself a positive operator. It follows immediately from the condition (ii) that $R$ is an injection. Assume that $T\in G_{H}$ is positive. Then, by our assumption and by part (i) of Proposition \ref{addit} $T_{\tau_{T}}\in G$ extends $T$, in other words $T=RT_{\tau_{T}}$. Every injection which maps the positive cone onto the positive cone is an order isomorphism.\medskip

Now let $S\in G$. It follows from Corollary \ref{rk0} that $RT_{\tau_{RS}}$ is the positive part of $RS$ in $G_{H}$, and since $R$ is an order isomorphism, it follows that $S^{+}=T_{\tau_{RS}}$ in $G$. Note that for $e\in E_{+}$ we have $\tau_{RS}\left(e\right)=\bigvee S\left[0_{E},e\right]_{H}\le \bigvee S\left[0_{E},e\right] \le S^{+}e=\tau_{RS}\left(e\right)$, where the last inequality follows from part (i) of Proposition \ref{addit} applied to $E$. We conclude that $G$ is a vector lattice with the specified lattice operations.\medskip

We now prove the last claim. Recall that by part (i) of Corollary \ref{rkm} $L\left(E,F\right)$ is monotone complete. Also recall that for vector lattices monotone completeness is equivalent to Dedekind completeness. Hence, if $G$ is full in $L\left(E,F\right)$, then by part (ii) of Proposition \ref{monotone} it is a Dedekind complete vector lattice, and moreover the embedding into $L\left(E,F\right)$ is monotone continuous. The last property together with part (i) of Corollary \ref{rkm} implies that the suprema of increasing nets of operators in $G$ are evaluated pointwise. The case of $G_{H}$ is analogous.
\end{proof}

\begin{proposition}\label{rrdp}
Assume that $E$ has the RDP and $H\subset E$ is full. Then, $L^{E}_{ob}\left(H,F\right)$ is a Dedekind complete vector lattice in which the lattice operations as well as the suprema of increasing nets are evaluated via Riesz-Kantorovich formulas. Moreover, the map $U:L^{E}_{ob}\left(H,F\right)\to L_{ob}\left(E,F\right)$ defined by $U\left(T\right):= T_{\tau_{T}}-T_{\tau_{-T}}$ is a lattice homomorphism.
\end{proposition}
\begin{proof}
First, we prove that $\tau_{T}$ is subadditive, for every $T\in L^{E}_{ob}\left(H,F\right)$. To that end, it is enough to show that for every $e,g\in E_{+}$ we have that $\left[0_{E},e+g\right]_{H}\subset\left[0_{E},e\right]_{H}+\left[0_{E},g\right]_{H}$ in $E$. Take $h\in \left[0_{E},e+g\right]_{H}$. By the Riesz decomposition property there are $u\in \left[0_{E},e\right]$ and $v\in \left[0_{E},g\right]$ such that $u+v=h$. Then, $u,v\in\left[0_{E},h\right]\subset H$, where the latter inclusion follows from fullness of $H$. Thus, $u\in \left[0_{E},e\right]_{H}$ and $v\in \left[0_{E},g\right]_{H}$, therefore $h\in\left[0_{E},e\right]_{H}+\left[0_{E},g\right]_{H}$.

Now the main claim follows from Corollary \ref{rk0} and the argument from the proof of Lemma \ref{rk}.\medskip

We now prove the last claim. Let us show that if $T,S\in L^{E}_{ob}\left(H,F\right)_{+}$, then $\tau_{T+S}=\tau_{T}+\tau_{S}$. Let $e\in E_{+}$. For every $h\in \left[0_{E},e\right]_{H}$ we have $\left(T+S\right)h=Th+Sh\le\tau_{T}\left(e\right)+\tau_{S}\left(e\right)$. Taking the supremum on the left yields $\tau_{T+S}\left(e\right)\le\tau_{T}\left(e\right)+\tau_{S}\left(e\right)$.

On the other hand, let $g,h\in \left[0_{E},e\right]_{H}$. The RDP applied to $g\le h+\left(e-h\right)$ yield $u\in\left[0_{E},h\right]\subset H$ and $v\in\left[0_{E},e-h\right]$ such that $g=u+v$. Then, $v=g-u\in H$, $w:=h+v\in H$, $h\le w\le h+e-h=e$, and $g=u+v\le w$. It follows that $Tg+Sh\le \left(T+S\right)w\le \tau_{T+S}\left(e\right)$. Taking the supremum on the left yields $\tau_{T}\left(e\right)+\tau_{S}\left(e\right)\le \tau_{T+S}\left(e\right)$.\medskip

It follows that $T\mapsto\tau_{T}$ is an additive map from $L^{E}_{ob}\left(H,F\right)_{+}$ into $L_{ob}\left(E,F\right)_{+}$. This map extends uniquely to a positive linear map from $L^{E}_{ob}\left(H,F\right)$ into $L_{ob}\left(E,F\right)$, and the extension is equal to $U$ (see the proof of part (iii) of Proposition \ref{addit}). It is left to show that this map is a homomorphism. First, for $T\in L_{ob}^{E}\left(H,F\right)$ we have that $T^{+}$ in $L_{ob}^{E}\left(H,F\right)$ is given via $T^{+}h=\tau_{T}\left(h\right)$. Hence, for every $e\in E_{+}$ we have $\tau_{T^{+}}\left(e\right)=\bigvee \tau_{T}\left(\left[0_{E},e\right]_{H}\right)=\tau_{T}\left(e\right)$ (see part (iii) of Proposition \ref{addit0}), and so $U\left(T^{+}\right)= T_{\tau_{T}}$. It is clear that $T_{\tau_{T}}$ is positive and $T_{\tau_{T}}\ge T_{\tau_{T}}-T_{\tau_{-T}}=U\left(T\right)$. On the other hand, if $S\in L_{ob}\left(E,F\right)_{+}$ is such that $S\ge U\left(T\right)$, then, as the latter is an extension of $T$, it follows that $\left.S\right|_{H}\ge T$, and so $S\ge T_{\tau_{T}}$, by part (i) of Proposition \ref{addit}. We conclude that $U\left(T^{+}\right)= T_{\tau_{T}}=U\left(T\right)^{+}$, as required.
\end{proof}

Let us introduce a slightly weaker property. Namely, we will say that $E$ has the \emph{almost decomposition property (ADP)} if $\left[0_{E},e+g\right]$ is the relative uniform closure of $\left[0_{E},e\right]+\left[0_{E},g\right]$. A linear operator is order bounded if and only if it is relatively uniformly continuous (see \cite[Example 7.5 and Proposition 7.8]{bctv}; note that it is immaterial whether the underlying OVS's are lattices). It follows that a pre-image of every relatively uniformly closed set is relatively uniformly closed. This allows us to prove the following result.

\begin{proposition}\label{rdp}
Assume that $E$ has the ADP. Then:
\item[(i)] $L_{ob}\left(E,F\right)$ is a Dedekind complete vector lattice in which the lattice operations as well as the suprema of increasing nets are evaluated via Riesz-Kantorovich formulas.
\item[(ii)] If $H\subset E$ is full and such that every $T\in L^{E}_{ob}\left(H,F\right)$ admits an extension in $L_{ob}\left(E,F\right)$. Then, $L^{E}_{ob}\left(H,F\right)$ is a Dedekind complete vector lattice in which the lattice operations as well as the suprema of increasing nets are evaluated via Riesz-Kantorovich formulas.
\end{proposition}
\begin{proof}
Let $T\in L_{ob}\left(E,F\right)$. Let us show that $\tau_{T}$ is subadditive. Fix $e,g\in E_{+}$. Then $A:=\tau_{T}\left(e\right)+\tau_{T}\left(g\right)-F_{+}$ is a relatively uniform closed subset of $F$ which contains $T\left(\left[0_{E},e\right]+\left[0_{E},g\right]\right)$. Since $T$ is order bounded, it follows that $T^{-1}A$ is a relatively uniformly closed subset of $E$ which contains $\left[0_{E},e\right]+\left[0_{E},g\right]$. As $E$ has the ADP, $T^{-1}A$ contains $\left[0_{E},e+g\right]$, therefore $T\left[0_{E},e+g\right]\subset A\le \tau_{T}\left(e\right)+\tau_{T}\left(g\right)$. We conclude that $\tau_{T}\left(e+g\right)\le \tau_{T}\left(e\right)+\tau_{T}\left(g\right)$.

The result now follows from Corollary \ref{rk0} and Lemma \ref{rk} applied to the case $H=E$ and $G=L^{E}_{ob}\left(H,F\right)$.\medskip

Let $H$ satisfy the extension hypothesis. Let $T\in L^{E}_{ob}\left(H,F\right)$, let $\widehat{T}\in L_{ob}\left(E,F\right)$ be its extension, and let $\widehat{T}^{+}$ be the positive part of $\widehat{T}$ in $L_{ob}\left(E,F\right)$. By (i) we have that $\widehat{T}^{+}e=\bigvee \widehat{T}\left[0_{E},e\right]$, for every $e\in E_{+}$. In particular, since $H$ is full, for every $h\in H_{+}$ we have $\widehat{T}^{+}h=\bigvee \widehat{T}\left[0_{E},h\right]=\bigvee T\left[0_{E},h\right]=\tau_{T}\left(h\right)$. It follows that $\left.\tau_{T}\right|_{H_{+}}=\left.\widehat{T}^{+}\right|_{H_{+}}$ is subadditive on $H_{+}$. Applying Corollary \ref{rk0} completes the proof of the main claim.
\end{proof}

Order unit spaces with the ADP have some special features. We say that an order unit space $\left(E,e\right)$ has the \emph{approximate interpolation property (AIP)} if for any $\varepsilon>0$ and finite $P,Q\subset E$ such that $P\le Q$, there is $f\in E$ such that $P\le f+\varepsilon e$ and $f-\varepsilon e\le Q$.

\begin{proposition}\label{aip}
For an order unit space $\left(E,e\right)$ the following conditions are equivalent:
\item[(i)] The norm completion $\widehat{E}$ of $E$ has the RDP;
\item[(ii)] $E$ has the approximate interpolation property;
\item[(iii)] For any $P,Q\subset E$ consisting of two elements and such that $P\le Q$, there is $f\in E$ such that $P\le f+\varepsilon e$ and $f-\varepsilon e\le Q$.
\end{proposition}
\begin{proof}
(i)$\Rightarrow$(ii): There is $f_{0}\in \widehat{E}$ such that $P\le f_{0} \le Q$. There is $f\in E$ such that $\|f-f_{0}\|\le\varepsilon$. Then, $f_{0}\le f+\varepsilon e$ and $f_{0}\ge f-\varepsilon e$, which immediately yields the required inequalities.\medskip

(ii)$\Rightarrow$(i): First we show that $\widehat{E}$ has the AIP. Let $\varepsilon>0$ and let finite $P,Q\subset \widehat{E}$ be such that $P\le Q$. For every $p\in P$ and $q\in Q$ there are $p',q'\in E$ such that $\|p-p'\|_{e},\|q-q'\|_{e}\le\frac{\varepsilon}{3}$. Let $P':=\left\{p'-\frac{\varepsilon}{3}e,~ p\in P\right\}$ and $Q':=\left\{q'+\frac{\varepsilon}{3}e,~ q\in Q\right\}$, which are finite subsets of $E$ satisfying $p'-\frac{\varepsilon}{3}e\le p\le q\le q'+\frac{\varepsilon}{3}e$, for every $p\in P$ and $q\in Q$. Since $E$ has the AIP and $P'\le Q'$ there is $f\in E$ such that $P'\le f+\frac{\varepsilon}{3} e$ and $f-\frac{\varepsilon}{3} e\le Q'$. For every $p\in P$ and $q\in Q$ we then have $f\ge p'-\frac{\varepsilon}{3} e -\frac{\varepsilon}{3} e\ge p-\varepsilon e$ and $f\le q'+\frac{\varepsilon}{3} e+\frac{\varepsilon}{3} e\le q+\varepsilon e$, as required.\medskip

Next we prove that if $E$ is norm complete and has the AIP then it has the RDP. Let $P=\left\{p,r\right\}$ and $Q=\left\{q,w\right\}$ satisfy $P\le Q$. Choose $f_{0}\in E$ such that $P\le f_{0}+e$ and $f_{0}-e\le Q$. We now inductively construct a sequence $\left(f_{n}\right)_{n\in\N_{0}}$ such that $f_{n}+\frac{1}{2^{n}}e\ge P$, $f_{n}-\frac{1}{2^{n}}e\le Q$, and $\|f_{n}-f_{n-1}\|_{e}\le \frac{3}{2^{n}}$, for every $n\in\N$. If $f_{n-1}$ is selected, let $P_{n-1}:=P\cup\left\{f_{n-1}-2^{1-n}e\right\}$ and $Q_{n-1}:=Q\cup\left\{f_{n-1}+2^{1-n}e\right\}$. The assumption of induction yields $P_{n-1}\le Q_{n-1}$, and so there is $f_{n}\in E$ such that $f_{n}+\frac{1}{2^{n}}e\ge P_{n-1}$ and $f_{n}-\frac{1}{2^{n}}e\le Q_{n-1}$, which immediately validates the first two requirements. We also have $f_{n}+\frac{1}{2^{n}}e\ge f_{n-1}-2^{1-n}e$ and $f_{n}-\frac{1}{2^{n}}e\le f_{n-1}+2^{1-n}e$, hence $-\frac{3}{2^{n}}e\le f_{n}-f_{n-1}\le \frac{3}{2^{n}}e$, as required.  Since $E$ is complete, there is a norm limit $f$ of the sequence. It is easy to see that $P\le f\le Q$.\medskip

(ii)$\Rightarrow$(iii) is trivial. For the converse we assume that $Q$ has two elements, but $P$ has $n$ elements and run induction over $n$. The cases $n=0,1$ are trivial. Assume that the claim is proven for $n\in\N$ and let $P=\left\{p_{0},...,p_{n}\right\}$. Then, apply the assumption of induction to $Q$ and  $\left\{p_{1},...,p_{n}\right\}$ and find $h\in E$ such that $\left\{p_{1},...,p_{n}\right\}\le h+\frac{\varepsilon}{2} e$ and $h-\frac{\varepsilon}{2} e\le Q$. Then apply the previous step to $\left\{h-\frac{\varepsilon}{2}e,p_{0}\right\}$ and $Q$ to find $g\in E$ with $p_{0},h-\frac{\varepsilon}{2}e\le g+\frac{\varepsilon}{2} e$ and $g-\frac{\varepsilon}{2} e\le Q$. Finally, let $f:=g+\frac{\varepsilon}{2} e$. We have $f-\varepsilon e=g-\frac{\varepsilon}{2} e\le Q$, $f+\varepsilon e=g+\frac{3\varepsilon}{2} e \ge  g+\frac{\varepsilon}{2} e\ge p_{0}$ and $f+\varepsilon e=g+\frac{3\varepsilon}{2} e \ge  h+\frac{\varepsilon}{2}e\ge p_{1},...,p_{n}$.

Finally, use similar arguments to run induction over the the size of $Q$.
\end{proof}

\begin{proposition}\label{adpou}
Let $E$ be an order unit space with the unit $e$ and the ADP. Then:
\item[(i)] $E$ has the approximate interpolation property.
\item[(ii)] If $H$ is such that $T\left[-e,e\right]_{H}$ is order bounded in $F$, for every $e\in E_{+}$, then $T$ can be extended to an order bounded operator from $E$ into $F$.
\end{proposition}
\begin{proof}
(i): We may assume that both $P,Q$ have two elements each, namely $P=\left\{p,r\right\}$ and $Q=\left\{q,w\right\}$. We have that $0_{E}\le q-p\le \left(q-p\right)+\left(w-r\right)=\left(q-r\right)+\left(w-p\right)$. Since $q-r,w-p\in E_{+}$, the ADP yields $u\in \left[0_{E},q-r\right]$ and $v\in \left[0_{E},w-p\right]$ such that $-\varepsilon e\le q-p -u-v\le \varepsilon e$. Then, $f:=q-u$ satisfies $f\le q$, $f=q-u\ge q-\left(q-r\right)=r$ and  $-\varepsilon e\le f-p -v\le \varepsilon e$, so that $f\ge p+v-\varepsilon e\ge p-\varepsilon e$ and $f\le p+v+\varepsilon e\le p+\left(w-p\right)+\varepsilon e=w+\varepsilon e$.\medskip

(ii): The set $T\left[-e,e\right]_{H}$ is order bounded in $F$ by some $f\in F_{+}$. For every $h\in H$ we have that $h\in  \|h\|_{e}\left[-e,e\right]_{H}$, and so $Th\in \|h\|_{e}T\left[-e,e\right]_{H}\subset \|h\|_{e}\left[-f,f\right]$. It follows that $\|\cdot\|_{e}f$ is a sublinear function from $E$ into $F$, which satisfies $Th\le \|h\|_{e}f$, for every $h\in H$. By Hahn-Banach-Kantorovich theorem there is a linear $\widehat{T}:E\to F$ which extends $T$ and such that $\widehat{T}g\le \|g\|_{e}f$, for every $g\in E$. In particular, $\widehat{T}\left[-e,e\right]\subset\left[-f,f\right]$, and so $\widehat{T}\in L_{ob}\left(E,F\right)$, as required.
\end{proof}

The condition appearing in part (ii) of Proposition \ref{adpou} is not equivalent to $E$-boundedness: if it was, part (ii) of Proposition \ref{rdp} would be valid for any full directed subspace of an ordered unit space, which is false (see Example \ref{adp}). The converse to part (i) of Proposition \ref{adpou} is false:

\begin{example}\label{nadpou}
Let $K:=\left[0,1\right]\cup\left[2,3\right]$, let $P\left(K\right)$ be the space of restrictions of polynomials to $K$, let $f:=\1_{\left[0,1\right]}$, $g:=\1_{\left[2,3\right]}$ and let $E$ be the span inside of $\Co\left(K\right)$ of $P\left(K\right)\cup \left\{f\right\}$. Note that $g=\1-f\in E$, and in fact $E$ is also the span of $P\left(K\right)\cup \left\{g\right\}$. Since $P\left(K\right)$ is norm dense in $\Co\left(K\right)$, it follows that the latter is the norm completion of $E$, and so $E$ has the AIP, according to Proposition \ref{aip}.

To see that $E$ does not have the ADP, note that $\left[\0,f\right]=\left\{tf,~ 0\le t\le 1\right\}$ and $\left[\0,g\right]=\left\{tg,~ 0\le t\le 1\right\}$. Indeed, if $p\in P\left(K\right)$ and $s\in \R$ are such that $\0\le p+sf\le f$, it follows that $p$ vanishes on $\left[2,3\right]$, hence $p=\0$ (the case of $g$ is similar). It follows that $\left[\0,f\right]+\left[\0,g\right]=\left\{tf+sg,~ 0\le t,s\le 1\right\}$ is norm closed in $E$, hence relatively uniformly closed. On the other hand, it is clearly not equal to $\left[\0,f+g\right]=\left[\0,\1\right]$, which contains various non-constant polynomials.
\qed\end{example}

\subsection{Relative uniform convergence and extension of order bounded operators}\label{sbob}

We start with a more detailed study of relative uniform convergence. Throughout the section $H$ is a subspace of a directed Archimedean ordered vector space $E$. It is easy to see that order intervals are relatively uniformly closed, as well as $E_{+}$. If $H\subset E$ is majorizing, then relative uniform convergence of $H$ is the restriction of relative uniform convergence on $E$. The following results are well known for Banach lattices, but we state them for ordered vector spaces endowed with a complete metrizable locally full linear topologies.

\begin{proposition}\label{ru1}
Let $\tau$ be a complete metrizable linear topology on $E$ such that $E_{+}$ is closed. Then, every $\tau$-convergent net in $E$ contains a sequence, which ru-converges to the same limit.
\end{proposition}
\begin{proof}
Since $\tau$ is metrizable, it is given by a pseudonorm $\|\cdot\|$. Moreover, according to \cite{vt} we may assume that for every $e\in E$ there are $e',e''\in E_{+}$ such that $\|e'\|+\|e''\|\le 2\|e\|$ and $e=e'-e''$. Let $\left(e_{\alpha}\right)_{\alpha\in A}\subset E$ be $\tau$-null, and so there is an increasing sequence $\left(\alpha_{n}\right)_{n\in\N}\subset A$ such that $\|e_{\alpha_{n}}\|\le \frac{1}{n2^{n+1}}$. For every $n\in\N$ let $e'_{n},e''_{n}\in E_{+}$ be such that $e_{\alpha_{n}}=e'_{n}-e''_{n}$, and $\|e'_{n}\|+\|e''_{n}\|\le 2\|e_{\alpha_{n}}\|\le \frac{1}{n2^{n}}$. Then, $h:=\sum\limits_{n\in\N}n\left(e'_{n}+e''_{n}\right)$ exists in $E$, and moreover, $\pm e_{\alpha_{n}}\le e'_{n}+e''_{n}\le \frac{1}{n}h$ yields $e_{\alpha_{n}}\toru 0_{E}$.
\end{proof}

\begin{proposition}\label{locfull}
Let $\tau$ be a complete metrizable locally full linear topology on $E$ such that $E_{+}$ is closed. Then:
\item[(i)] $\tau$-closed sets are precisely ru-closed sets.
\item[(ii)] If $F$ is an ordered vector space endowed with a locally full linear topology, and $T:E\to F$ is order bounded, then it is topologically continuous.
\item[(iii)] If every net in $E$ which decreases to $0_{E}$ is $\tau$-null, then $\mathrm{ru}=\mathrm{do}$, and every order bounded operator from $E$ into any Archimedean ordered vector space is directed order continuous.
\end{proposition}
\begin{proof}
(i): Since $\mathrm{ru}$ is the strongest locally full linear convergence on $E$, it follows that every $\tau$-closed set is ru-closed. The converse follows from Proposition \ref{ru1}.\medskip

(ii): Let $C\subset F$ be closed. Then, it is ru-closed, and since $T$ is ru-continuous, it follows that $T^{-1}C$ is ru-closed in $E$. According to (i), $T^{-1}C$ is in fact $\tau$-closed. Hence, the preimage of every closed set is closed, and so $T$ is continuous.\medskip

(iii): For the first claim see the proof of \cite[Proposition 2.24]{brtgww}. If $F$ is an Archimedean ordered vector space and $T:E\to F$ is order bounded, then it is ru-continuous. Since $\mathrm{ru}=\mathrm{do}$, it follows that $T$ is do-to-ru-continuous, and thus do-continuous.
\end{proof}
We say that a net $\left(e_{\alpha}\right)_{\alpha\in A}\subset E$ converges to $e\in E$ \emph{with respect to} $P\subset E_{+}$ (denoted $e_{\alpha}\xrightarrow[]{P} e$) if $P$ \emph{dominates the tails} of $\left(e_{\alpha}-e\right)_{\alpha\in A}$. Of course a net is (directed) order convergent if and only if it converges with respect to (downward directed) set with infimum $0_{E}$, and it is relatively uniformly convergent if and only if it converges with respect to the set of the form $\left\{\frac{1}{n}h\right\}_{n\in\N}$, for some $h\in E_{+}$.

Note that $\left(e_{\alpha}\right)_{\alpha\in A}$ is Cauchy with respect to $P$ (i.e. if $e_{\alpha}-e_{\beta}\xrightarrow[]{P} 0_{E}$) if and only if for every $p\in P$ there is $\alpha_{p}$ such that whenever $\alpha,\beta\ge\alpha_{p}$ we have $e_{\alpha}-e_{\beta}\le p$.

A net $\left(f_{\beta}\right)_{\beta\in B}$ is a \emph{quasi-subnet} of $\left(e_{\alpha}\right)_{\alpha\in A}$ if for every $\alpha_{0}\in A$ there is $\beta_{0}\in B$ such that $\left\{f_{\beta}\right\}_{\beta\ge\beta_{0}}\subset\left\{e_{\alpha}\right\}_{\alpha\ge\alpha_{0}}$.

\begin{proposition}\label{ruu}
Let $\eta$ be a linear convergence structure on $E$ such that the order intervals are $\eta$-closed. Let $P\subset E_{+}$. If a net is Cauchy with respect to $P$ and has a $\eta$-convergent quasi-subnet, then it converges to the same limit with respect to $P$.
\end{proposition}
\begin{proof}
Let $\left(e_{\alpha}\right)_{\alpha\in A}\subset E$ be Cauchy with respect to $P\subset E_{+}$. Assume that $\left(f_{\beta}\right)_{\beta\in B}$ is its $\eta$-null quasi-subnet. Fix $p\in P$ and find $\alpha_{p}$ such that $-p\le e_{\alpha}-e_{\gamma}\le p$, whenever $\alpha,\gamma\ge\alpha_{p}$; in particular, $e_{\alpha}\in \left[e_{\gamma}-p,e_{\gamma}+p\right]$. Fix $\gamma\ge\alpha_{p}$ and let $\beta_{p}$ be such that $\left(f_{\beta}\right)_{\beta\ge\beta_{p}}\subset \left(e_{\alpha}\right)_{\alpha\ge\alpha_{p}}\subset \left[e_{\gamma}-p,e_{\gamma}+p\right]$. As the latter interval is $\eta$-closed, and $f_{\beta}\xrightarrow[]{\eta} 0_{E}$, it follows that $0_{E}\in\left[e_{\gamma}-p,e_{\gamma}+p\right]$. Hence, $-p\le e_{\gamma} \le p$, whenever $\gamma\ge\alpha_{p}$. We conclude that $e_{\alpha}\xrightarrow[]{P} 0_{E}$.
\end{proof}

\begin{corollary}\label{ruh}
If a $\|\cdot\|_{h}$-Cauchy net in $E$ relatively uniformly converges, for some $h\in E_{+}$, then $h$ can be chosen as the regulator of the convergence.
\end{corollary}

$E$ is called \emph{relatively uniformly complete} if every relatively uniformly Cauchy net has a relative uniform limit. Note that it is enough to only consider sequences in this definition. The following is easy to derive from Corollary \ref{ruh}.

\begin{corollary}\label{ruf}
Relative uniform completeness passes down to full subspaces.
\end{corollary}

\begin{corollary}\label{ruur}
Let $\eta$ be a linear convergence structure on $E$ such that the order intervals are $\eta$-closed. Let $\theta$ be one of order, directed order and relative uniform convergences. If $\theta$ is stronger than $\eta$, and $E$ is $\eta$-complete, then $E$ is $\theta$-complete.
\end{corollary}
\begin{proof}
Every $\theta$-Cauchy net is $\eta$-Cauchy, hence $\eta$-convergent, and therefore $\theta$-convergent to the same limit, according to Proposition \ref{ruu}.
\end{proof}

We now focus on ru-closures and ru-adherences. The rest of the results in this subsection are analogous to those in \cite{bt}, but generalized from the level of Archimedean vector lattices to that of directed Archimedean OVS.

\begin{proposition}[Cf. Proposition 3.3, \cite{bt}]\label{ruc}
Assume that $H$ is majorizing in $E$. Then, every element in $\overline{H}^{\mathrm{ru}}$ is a supremum (in $E$) of a countable set in $H$. Moreover, $\overline{H}^{\mathrm{ru}}_{+}=\overline{H_{+}}^{\mathrm{ru}}$.
\end{proposition}
\begin{proof}
First, assume that $e\in \overline{H}^{\mathrm{ru}}_{1}$. There is $g\in E_{+}$ such that for every $n\in\N$ there is $h_{n}\in H$ with $-\frac{1}{n}g\le e-h_{n}\le \frac{1}{n}g$. Since $H$ is majorizing in $E$, there is $f\in H$ such that $g\le f$. For $n\in\N$ define $e_{n}:=h_{n}+\frac{1}{n}f\in H$; we have $$e_{n}-\frac{2}{n}f=h_{n}-\frac{1}{n}f\le h_{n}-\frac{1}{n}g\le e \le h_{n}+\frac{1}{n}g\le h_{n}+\frac{1}{n}f=e_{n}.$$ Hence $e\le e_{n}\toru e$, which implies $e=\bigwedge\limits_{n\in\N}e_{n}$. In particular, $e\ge 0_{E}$ implies $e_{n}\in H_{+}$, for every $n\in\N$, and so $e\in\overline{H_{+}}^{\mathrm{ru}}_{1}$.

After this, to complete the proof run the transfinite induction as in the proof of \cite[Proposition 3.3(ii)]{bt}.
\end{proof}

\begin{proposition}[Cf. Proposition 4.2, \cite{bt}]
The intersection of any nonempty collection of relatively uniformly complete subspaces of an Archimedean ordered vector space is relatively uniformly complete.
\end{proposition}
\begin{proof}
Let $\Ho\ne\varnothing$ be a collection of relatively uniformly complete subspaces of $E$, and let $G:=\bigcap\Ho$. Assume that a net in $G$ is $\|\cdot\|_{g}$-Cauchy, for some $g\in G_{+}$. Then, this net is $\|\cdot\|_{g}$-Cauchy in each $H\in\Ho$, and so by Corollary \ref{ruh} it $\|\cdot\|_{g}$-converges there (and hence in $E$) to some $h\in H$. By uniqueness of the limit in $E$, $h$ does not depend on $H$, and so lies in $G$. It follows that the net $\|\cdot\|_{g}$-converges to $h$ in $G$.
\end{proof}

\begin{corollary}[Cf. Proposition 4.1, \cite{bt}]
Assume that $E$ is relatively uniformly complete. Then:
\item[(i)] If $H$ is relatively uniformly closed, then it is relatively uniformly complete.
\item[(ii)] If $H$ is majorizing in $E$, then its $\mathrm{ru}_{E}$-closure is the intersection of all relatively uniformly complete subspaces $H\subset G\subset E$. In particular, if $H$ is relatively uniformly complete, then it is relatively uniformly closed.
\end{corollary}

Recall that if $H$ is not majorizing in $E$ it can be simultaneously ``intrinsically'' relatively uniformly complete and relatively uniformly dense, e.g. $H=c_{00}\subset c_{0}=E$.

\begin{proposition}[Cf. Theorem 5.1, \cite{bt}]\label{ruex}
Assume that $H$ is majorizing in $E$. Let $F$ be a relatively uniformly complete directed ordered vector space, and let $T:H\to F$ be regular. Then, $T$ extends uniquely to a regular operator from $\overline{H}^{\mathrm{ru}}$ to $F$. Moreover, ``regular'' can be replaced with ``order bounded'', if $H$ has the ADP or an order unit.
\end{proposition}
\begin{proof}
The main claim is proven analogously to the cited result from \cite{bt}, so we only address the ``moreover'' clause. For the ADP case, according to part (i) of Proposition \ref{rdp}, if $J:F\to F^{\delta}$ is the inclusion, then there are positive $T^{\pm}:H\to F^{\delta}$ such that $T^{+}-T^{-}=JT$. After that, the proof proceeds as in \cite{bt}. For the order unit case, argue as in part (ii) of Proposition \ref{adpou}.
\end{proof}

In general there can be no extension of an order bounded operator, see Example \ref{noext}.

\begin{proposition}[Cf. Theorem 5.2, \cite{bt}]\label{ruexp}
Let $F$ be an Archimedean ordered vector space, let $T:E\to F$ be an order bounded operator, and assume that $H$ is majorizing in $E$ and such that $\overline{H}^{\mathrm{ru}}=E$. For each of the following properties, $T$ has it if and only if $\left.T\right|_{H}$ does: positivity, bi-positivity and supremum continuity.
\end{proposition}
\begin{proof}
It is clear that positivity and bi-positivity pass down from $T$ to $\left.T\right|_{H}$; the same is also true for supremum continuity, since according to Proposition \ref{ruc} $H$ is supremum-dense in $E$, and so the inclusion map from $H$ into $E$ is supremum continuous.

Conversely, bi-positivity and supremum continuity pass up from $T$ to $\left.T\right|_{H}$ according to \cite[Proposition 2.3]{aerz}. The case of positivity is handled analogously to \cite[Theorem 5.2]{bt}, in light of Proposition \ref{ruc}.
\end{proof}

\section{Extension of operators on pervasive spaces}\label{secperv}

This section is dedicated to the class of ordered vector spaces whose behavior in regards to order convergence is very close to that of vector lattices, namely pervasive spaces. We start with a review of order convergence and continuity on these spaces, and then proceed to the results about extensions of operators on them.

\subsection{Pervasive spaces}\label{sbperv}

Throughout the subsection $E$ is a directed Archimedean ordered vector space and $H\subset E$ is a linear subspace. We say that $H$ is \emph{order dense}\footnote{We use the this term in a way which is consistent with the vector lattice literature. In the literature on ordered vector spaces this term often means what we call ``supremum dense''.} if $H\cap\left(0_{E},e\right]\ne\varnothing$, for every $e>0_{E}$. We will say that $H$ is \emph{supremum dense in $E_{+}$} if $e=\bigvee e^{\Downarrow_{H}}$, for every $e\in E_{+}$. It is easy to see that $H$ is supremum dense in $E$ if and only if it is majorizing and supremum dense in $E_{+}$. If $H$ is order dense in $E$, then it is supremum-dense in $E_{+}$ (see \cite[Proposition 3.1]{brtgww}). We will need more results from \cite{brtgww}.

\begin{proposition}[Proposition 3.3, \cite{brtgww}, cf. Lemma 3.7.11, \cite{kv}]\label{4.4}
Assume that $E$ is embedded into an Archimedean vector lattice $F$ as an order dense subspace. Let $g\in F_{+}$ and let $P\subset F_{+}$ be directed and such that $g=\bigvee P$. Then, there is a directed $Q\subset \bigcup\limits_{p\in P}\left[0_{F},p\right]_{E}$ such that $g=\bigvee Q$.
\end{proposition}

In particular, for every $f\in F_{+}$ we can apply Proposition \ref{4.4} to $P=\left\{f\right\}$ and obtain a $Q\subset F_{+}$ such that $Q\uparrow f$.

\begin{corollary}[Corollary 3.4, \cite{brtgww}]\label{4.5}
Assume that $E$ is embedded into an Archimedean vector lattice $F$ as an order dense subspace. Let $P\subset F_{+}$ be directed downward with $\bigwedge P=0_{F}$ and such that there are $p_{0}\in P$ and $e_{0}\in E$ with $p_{0}\le e_{0}$. Then, there is $Q\subset E\cap \left[0_{F},e_{0}\right]\cap \left(P+F_{+}\right)$ which is directed downward and such that $\bigwedge Q=0_{F}$.
\end{corollary}

We will say that $E$ is \emph{pervasive} if it embeds as an order dense subspace of a vector lattice. According to \cite[Theorem 3.5]{brtgww}, this condition is equivalent to the fact that $E$ is order dense in $E^{\delta}$. We will shorten ``pervasive directed Archimedean ordered vector space'' to ``pervasive Archimedean OVS''. Clearly, an order dense subspace of a pervasive Archimedean OVS is pervasive. If $E$ is a pervasive Archimedean OVS, for every $e\in E^{\delta}_{+}$ there is $Q\subset E_{+}$ such that $Q\uparrow e$.

\begin{proposition}\label{sdsc}
Assume that $E$ is pervasive. Then:
\item[(i)] Order convergence agrees with directed order convergence on $E$. If $F$ is an ordered vector space, $T\in L\left(E,F\right)$ is monotone continuous if and only if $T\in L_{oc}\left(E,F\right)_{+}$ if and only if $T\in L_{doc}\left(E,F\right)_{+}$.
\item[(ii)] If $H$ is full in $E$, or is supremum-dense in $E_{+}$, then the embedding of $H$ into $E$ is supremum-continuous.
\item[(iii)] If $H$ is order dense in $E$ if and only if it is pervasive and supremum-dense in $E_{+}$.
\item[(iv)] If $H$ is order dense in $E$, then it is dense in $E$ with respect to order convergence. Moreover, $\left[-e,e\right]_{H}$ is dense in $\left[-e,e\right]$, for every $e\in E_{+}$.
\end{proposition}
\begin{proof}
(i): The first claim is \cite[Proposition 3.9]{brtgww}. For the second claim combine it with part (ii) of Proposition \ref{4.10}.\medskip

(ii): Let $P\subset H$ be such that $\bigwedge_{H}P=0_{E}$.

First, we assume that $H$ is full in $E$. Suppose that $e\in E$ is such that $e\le P$, but $e\not\le 0_{E}$. Since $E$ is order dense in $E^{\delta}$, there is $0_{E}<f\in E$ such that $f\le e^{+}$ (the positive part in $E^{\delta}$; by assumption it is strictly positive). Note that since $p\ge e$ and $p\ge 0_{E}$, it follows that $p\ge f$, for every $p\in P$. In particular, $0_{E}\le f\le p\in H$ yields $f\in H$, but then $\bigwedge_{H}P=0_{E}$ implies $f=0_{E}$, contradiction.\medskip

Next, assume that $H$ is supremum-dense in $E_{+}$. Suppose that $f:=\bigwedge_{E^{\delta}}P>0_{E^{\delta}}$. Since $E$ is pervasive, it is order dense in $E^{\delta}$ and so there is $e\in \left(0_{E^{\delta}},f\right]\cap E$. If $h\in e^{\Downarrow_{H}}$, then $h\le e\le f\le  P$ implies $h\le 0_{E}$. Hence, $e=\bigvee e^{\Downarrow_{H}}\le 0_{E}$. This contradiction shows that $\bigwedge_{E^{\delta}}P=f=0_{E^{\delta}}$, which implies that $\bigwedge_{E}P=0_{E}$.\medskip

(iii): Necessity: Since $E$ is pervasive, it is order dense in $E^{\delta}$. As $H$ is order dense in $E$, it is also order dense in $E^{\delta}$, hence it is pervasive. Moreover, by \cite[Proposition 3.1]{brtgww} $H$ is supremum dense in $E^{\delta}_{+}$, hence in $E_{+}$.

Sufficiency: Let $J:H\to E$ be the inclusion map, which is an order embedding, and is supremum-continuous according to (ii). By Theorem \ref{supcont1} there is an injective homomorphism $J^{\delta}:H^{\delta}\to E^{\delta}$. Hence, we may assume that $H\subset H^{\delta}\subset E^{\delta}$. Fix $e>0_{E}$. Since $H$ is supremum-dense in $E_{+}$, there is $h\in H$ such that $h\le e$ and $h\not\le 0_{H}$. Take $h^{+}\in H^{\delta}$, and using the fact that $H$ is order dense in $H^{\delta}$, find $g\in H$ such that $0_{H^{\delta}}<g\le h^{+}$. As $e\ge h$ and $e\ge 0_{E}$, we have $e\ge h^{+}\ge g$.\medskip

(iv): Let $f\in \left[-e,e\right]$, and let $f^{+}$ be its positive part in $E^{\delta}$. Since $H$ is order dense in $E^{\delta}$, by the comment after Proposition \ref{4.4} there is a directed $P\subset H$ such that $P\uparrow f^{+}$. Similarly, there is a directed $Q\subset H$ such that $Q\uparrow f^{-}$. Then $P-Q\subset H$ viewed as a net indexed by $P\times Q$ order converges to $f$ in $E^{\delta}$, hence in $E$. Finally, for every $p\in P$ and $q\in Q$ we have that $0_{E^{\delta}}\le p\le f^{+}\le e$ and $0_{E^{\delta}}\le q\le f^{-}\le e$ imply $p-q\in \left[-e,e\right]_{H}$.
\end{proof}

The following example exhibits a full and order dense subspace $F\subset E$, which is also supremum-dense in $E_{+}$, where $E$ is an order unit space, such that the embedding of $F$ into $E$ is not $\uparrow$-continuous. It also demonstrates a full and order dense subspace $G\subset E$, which is fails to be supremum-dense in $E_{+}$.

\begin{example}\label{odensen}
Let $\N_{\8}:=\N\cup \left\{\8\right\}$, -- the one point compactification of $\N$, and let $E:=\Co\left(\N_{\8}\right)+\R f\subset\R^{\N_{\8}}$, where $f:\N_{\8}\to\R$ is defined by $f\left(n\right)=n$, for every $n\in\N$, and $f\left(\8\right)=0$. As a linear subspace of an Archimedean vector lattice, $E$ is Archimedean. It is easy to see that $f+\1$ is an order unit for $E$.

Let $F:=\Co\left(\N_{\8}\right)$ and let $G:=\left\{f\in F,~ f\left(\8\right)=0\right\}$. Let us show that $F=E_{\1}$. The direct inclusion is clear, and so we only prove the converse one. Assume that $e\in\Co\left(\N_{\8}\right)$ and $r,s\in \R$ are such that $-r\1\le e+sf\le r\1$. It follows that $sf\le r\1-e$. Since $f$ is positive and unbounded, while $r\1-e$ is bounded, we conclude that $s\le 0$. Analogously, $s\ge 0$, and so $s=0$. It follows that $e+sf=e\in \Co\left(\N_{\8}\right)$. Thus, $F$ is full in $E$.

Next, we establish that $F$ is order dense in $E$ and supremum dense in $E_{+}$. Let $e\in\Co\left(\N_{\8}\right)$ and $s\in \R$ be such that $g:=e+sf\ge \0$. Arguing as before one shows that $s\ge 0$. If $s=0$, there is nothing to prove and so we assume that $s>0$. For $n> \|e\|$ we have $g\left(n\right)=n+e\left(n\right)>0$, and so $\0<g\left(n\right)\1_{\left\{n\right\}}\in \left[\0,g\right]_{F}$. Moreover, it is easy to see that $g=\bigvee\left(\left\{g\left(n\right)\1_{\left\{n\right\}}\right\}_{n\in\N}\cup\left\{e\right\}\right)$. These arguments prove the claim.

We now show that the embedding of $F$ into $E$ is not $\uparrow$-continuous. For $n\in\N$ let $e_{n}:=\1_{\left\{1,...,n\right\}}$. It is easy to see that $e_{n}\uparrow \1$ in $F$. However, $e_{n}\le f\not\ge\1$, for every $n\in\N$, and so $e_{n}\not\uparrow \1$ in $E$.

Clearly, $G$ is a full directed order dense subspace of $F$, hence of $E$. Let us show that $G$ fails to be supremum dense in $E_{+}$. Indeed, for every $g\in G$ such that $g\le \1$, it follows that $g\left(n\right)\le 1\le f\left(n\right)$, for every $n\in\N$. Hence, $\1^{\Downarrow_{G}}\le f\not\ge\1$, and so $\bigvee\1^{\Downarrow_{G}}\ne \1$.
\qed\end{example}

It turns out that order continuity of a positive operator can be ``lifted'' from a large enough subspace.

\begin{proposition}\label{ocodm}
Assume that $E$ is pervasive and $H$ is order dense and majorizing. Let $F$ be an ordered vector space and let $T\in L\left(E,F\right)_{+}$ be such that $\left.T\right|_{H}$ is order continuous. Then, $T$ is order continuous.
\end{proposition}
\begin{proof}
By part (i) of Proposition \ref{sdsc} for both $T$ and $\left.T\right|_{H}$ order continuity is equivalent to monotone continuity.

Let $P\downarrow 0_{E}$. By Corollary \ref{4.5} there is $Q\subset H\cap\left(P+E_{+}\right)$ such that $Q\downarrow 0_{E}$. As $\left.T\right|_{H}$ is monotone continuous it follows that $TQ\downarrow 0_{F}$. For every $q\in Q$ there is $p_{q}\in P$ such that $p_{q}\le q$. It follows that $0_{F}\le \bigwedge TP\le \bigwedge\limits_{q\in Q}Tp_{q}\le \bigwedge TQ=0_{F}$. Thus, $T$ is monotone continuous.
\end{proof}

Proposition \ref{ocodm} is not valid if $H$ is supremum-dense, but not order dense, even if $E$ is a vector lattice, $F=\R$ and $T$ is a homomorphism.

\begin{example}\label{onuor}
Let $F:=c$, let $\nu:=\delta_{\8}$ and $\mu:=\sum\limits_{n\in\N}\frac{1}{2^{n}}\delta_{n}$. Note that $\nu$ is a homomorphism, while $\mu$ is order continuous. Let $E:=\ker\left(\nu-\mu\right)$. We claim that $E^{\rho}=F$. First, let $f\in F$, and assume that $r:=\nu\left(f\right)-\mu\left(f\right)\ge0$. Then, let $e,g\in E$ be defined by $e\left(1\right):=f\left(1\right)+2r$, $g\left(2\right):=f\left(2\right)+4r$, $e\left(n\right):=f\left(n\right)$, $n>1$, and $g\left(n\right):=f\left(n\right)$, for $n\ne 2$. Then, $f=e\wedge g$. Assume that $r<0$. Then, for $1<m\in\N$ define $e_{m}\in E$ by $e_{m}\left(n\right):=f\left(n\right)$, when $n<m$ and $e_{m}\left(n\right):=f\left(n\right)+\frac{2^{m-1}\left|r\right|}{2^{m-1}-1}$, when $n\ge m$. Then,
\begin{align*} \nu\left(e_{m}\right)&=\nu\left(f\right)+\frac{2^{m-1}\left|r\right|}{2^{m-1}-1}=\mu\left(f\right)+\frac{2^{m-1}\left|r\right|}{2^{m-1}-1}+r\\
&=\mu\left(f\right)+\frac{\left|r\right|}{2^{m-1}-1}=\mu\left(f\right)+\frac{1}{2^{m-1}}\frac{2^{m-1}\left|r\right|}{2^{m-1}-1}=\mu\left(e_{m}\right),
\end{align*}
hence $e_{m}\in E$, and $f=\bigwedge\limits_{m>1}e_{m}$. We conclude that $E$ is supremum dense in $F$. Moreover, every $f\in F$ with $\nu\left(f\right)-\mu\left(f\right)\ge 0$ is an infimum of two elements of $E$, similarly if $\nu\left(f\right)-\mu\left(f\right)\le 0$, then $f$ is an supremum of two elements of $E$.

Since the embedding of $E$ into $F$ is order continuous, and $\mu$ is order continuous on $F$, its restriction is order continuous on $E$. Since $\nu$ is a homomorphism on $F$, its restriction is a pre-homomorphism on $E$. Since these two restrictions coincide on $E$, it follows that $\left.\nu\right|_{E}$ is an order continuous pre-homomorphism. However, $\left(\left.\nu\right|_{E}\right)^{\rho}=\nu$ is not order continuous, despite having an order continuous restriction to a supremum dense subspace.
\qed\end{example}

\begin{remark}
Proposition \ref{ocodm} is also not valid if $T$ is non-positive order bounded even if $E$ is a Dedekind complete Banach lattice, and $H$ is an order dense majorizing norm-closed sublattice, and $F=\R$ (see Example \ref{example}). However, there is an analogous results (Proposition \ref{rucompl}) where we drop positivity but add some additional assumptions about $H$.

We also note that Proposition \ref{ocodm} is not valid if $H$ is not majorizing, even if $E,F$ are Banach lattices, $E$ Dedekind complete, $H$ is an order dense ideal and $T$ is a surjective homomorphism. Indeed, take the quotient map $Q:\ell_{\8}\to \ell_{\8}\slash c_{0}$, which is not order continuous, but its restriction to $H=c_{0}$ is zero, hence order continuous.
\qed\end{remark}

The following result can be viewed as an approximate Riesz Decomposition property.

\begin{proposition}\label{ardp}
Let $F$ be an Archimedean vector lattice and let $E\subset F$ be order dense. For any $e,f\in F_{+}$ and $g\in \left[0_{F}, e+f\right]$ there are directed $P\subset\left[0_{F}, e\right]_{E}$ and $Q\subset\left[0_{F}, f\right]_{E}$ such that $P+Q\uparrow g$. In particular, $\overline{\left[0_{F}, e\right]_{E}+\left[0_{F}, f\right]_{E}}^{1}_{\mathrm{o}}=\left[0_{F}, e+f\right]$.
\end{proposition}
\begin{proof}
By the comment after Proposition \ref{4.4} there are $P\subset\left[0_{F},e\wedge g\right]_{E}$ and $Q\subset\left[0_{F},g-e\wedge g\right]_{E}\subset \left[0_{F}, f\right]_{E}$ such that $P\uparrow e\wedge g$ and $Q\uparrow g-e\wedge g$ in $F$. It follows that $P+Q\uparrow g$, hence $P+Q\too g$.
\end{proof}

Taking $F=E^{\delta}$ yields the following.

\begin{corollary}\label{ardp2}
Let $E$ be pervasive. For any $e,f\in E_{+}$ and $g\in \left[0_{E}, e+f\right]$ there are directed $P\subset\left[0_{E}, e\right]$ and $Q\subset\left[0_{E}, f\right]$ such that $P+Q\uparrow g$. In particular, $\overline{\left[0_{E}, e\right]+\left[0_{E}, f\right]}^{1}_{\mathrm{o}}=\left[0_{E}, e+f\right]$.
\end{corollary}

\subsection{Riesz-Kantorovich-Veksler theorem}\label{s4}

Throughout the subsection $E$ is a pervasive Archimedean OVS, $H\subset E$ is an order dense directed subspace and $F$ is a Dedekind complete vector lattice. We denote the collection of all order continuous operators from $E$ into $F$ by $L_{oc}\left(E,F\right)$. Recall that $L^{E}_{ob}\left(H,F\right)$ is the set of $E$-bounded operators $T:H\to F$, i.e. such that $T\left[0_{E},e\right]_{H}$ is order bounded in $F$, for every $e\in E_{+}$. For such operators we denote $\tau_{T}\left(e\right):=\bigvee T\left[0_{E},e\right]_{H}$, for $e\in E_{+}$, which is a superadditive function, see Section \ref{s2}. The following lemma will allow us to prove versions of both Riesz-Kantorovich and Veksler theorems simultaneously.

\begin{lemma}\label{rkv}
Let $T\in L^{E}_{ob}\left(H,F\right)$ be such that if $H\supset Q\downarrow 0_{E}$, then $\bigwedge\left|TQ\right|= 0_{F}$. Then, $\tau_{T}$ is subadditive and $T_{\tau_{T}}$ is order continuous.
\end{lemma}
\begin{proof}
Let $e,g\in E_{+}$ and let $h\in \left[0_{E},e+g\right]_{H}$. According to Proposition \ref{ardp} there are directed $P\subset\left[0_{E}, e\right]_{H}$ and $Q\subset\left[0_{E}, g\right]_{H}$ such that $P+Q\uparrow h$. We have that $h-P-Q\downarrow 0_{E}$, and so $\bigwedge\limits_{p\in P,~q\in Q}\left|T\left(h-p-q\right)\right|= 0_{F}$. For every $p\in P$ and $q\in Q$ we have $$Th=Tp+Tq+T\left(h-p-q\right)\le Tp+Tq+\left|T\left(h-p-q\right)\right|\le \tau\left(e\right)+\tau\left(g\right)+\left|T\left(h-p-q\right)\right|.$$ It follows that $Th\le \tau\left(e\right)+\tau\left(g\right)+\bigwedge\limits_{p\in P,~q\in Q}\left|T\left(h-p-q\right)\right|= \tau\left(e\right)+\tau\left(g\right)$. As $h$ was arbitrary, we conclude that $\tau\left(e+g\right)\le \tau\left(e\right)+\tau\left(g\right)$. Subadditivity is proven.\medskip

According to Proposition \ref{addit}, $\tau_{T}$ extends to a positive operator $T_{\tau_{T}}:E\to F$. By part (i) of Proposition \ref{sdsc}, in order to establish its order continuity it is enough to prove that if $0_{E}\le e_{\alpha}\uparrow e$ then $\bigvee\limits_{\alpha} \tau_{T}\left(e_{\alpha}\right)=\tau_{T}\left(e\right)=\bigvee T\left[0_{E},e\right]_{H}$. We always have $\tau_{T}\left(e\right)\ge \tau_{T}\left(e_{\alpha}\right)$, for every $\alpha$. Hence, it is sufficient to show that for every $h\in\left[0_{E},e\right]_{H}$ we have $\bigvee\limits_{\alpha} \tau_{T}\left(e_{\alpha}\right)\ge Th$, which is equivalent to $\bigwedge\limits_{\alpha}\left(Th-\tau_{T}\left(e_{\alpha}\right)\right)^{+}=0_{F}$.\medskip

For every $\alpha$ denote $g_{\alpha}:=h\wedge\left(e-e_{\alpha}\right)$ and fix $\alpha_{0}$. We have that $h\ge g_{\alpha}\downarrow 0_{E}$, and so according to Corollary \ref{4.5} there is $Q\subset \bigcup\limits_{\alpha\ge\alpha_{0}}\left[g_{\alpha},g_{\alpha_{0}}\right]_{H}$ such that $Q\downarrow 0_{E}$. By our assumption we have $\bigwedge\limits_{q\in Q}\left|Tq\right|= 0_{F}$. For every $q\in Q$ there is $\alpha_{q}\ge \alpha_{0}$ such that $q\in \left[g_{\alpha_{q}},g_{\alpha_{0}}\right]_{H}\subset \left[g_{\alpha_{q}},h\right]$, so that $$0_{E}\le h-q\le h-g_{\alpha_{q}}=h-h\wedge\left(e-e_{\alpha_{q}}\right)=\left(h-e+e_{\alpha_{q}}\right)^{+}\le e_{\alpha_{q}}.$$ This implies $Th-Tq=T\left(h-q\right)\le \tau_{T}\left(h-q\right)\le \tau_{T}\left(e_{\alpha_{q}}\right)$, so that $Th-\tau_{T}\left(e_{\alpha_{q}}\right)\le Tq$, hence $\left(Th-\tau_{T}\left(e_{\alpha_{q}}\right)\right)^{+}\le \left(Tq\right)^{+}\le \left|Tq\right|$, for every $q\in Q$. Thus, $\bigwedge\limits_{\alpha}\left(Th-\tau_{T}\left(e_{\alpha}\right)\right)^{+}\le \bigwedge\limits_{q\in Q}\left(Th-\tau_{T}\left(e_{\alpha_{q}}\right)\right)^{+}\le \bigwedge\limits_{q\in Q}\left|Tq\right|=0_{F}$, as required.
\end{proof}

Denote $L_{oc}^{E}\left(H,F\right):=L^{E}_{ob}\left(H,F\right)\cap L_{oc}\left(H,F\right)$, which is a full subspace of $L_{oc}\left(H,F\right)$. Note that if $H$ is majorizing, then $L_{oc}^{E}\left(H,F\right)=L_{oc}\left(H,F\right)$. We now prove a version of Riesz-Kantorovich theorem for $L_{oc}\left(E,F\right)$.

\begin{theorem}\label{veksler0}
The restriction to $H$ is an order isomorphism from $L_{oc}\left(E,F\right)$ onto $L_{oc}^{E}\left(H,F\right)$. Both of these spaces are Dedekind complete vector lattices where the lattice operations as well as supremuma of increasing nets are evaluated by Riesz-Kantoronich formulas \cite[Theorem 1.59]{at}.
\end{theorem}
\begin{proof}
Let us show that $L_{oc}\left(E,F\right)$ and $L_{oc}^{E}\left(H,F\right)$ satisfy the conditions of Lemma \ref{rk}. Note that these are full subspaces of $L\left(E,F\right)$ and $L\left(H,F\right)$ respectively by virtue of Proposition \ref{bsc0} and part (i) of Proposition \ref{sdsc}. First, the condition (iii) is satisfied by Lemma \ref{rkv}. Assume that $S\in L_{oc}\left(E,F\right)$ vanishes on $H$. Then since $H$ is order dense in $E$, by part (iv) of Proposition \ref{sdsc} it is dense with respect to order convergence. Therefore, $\ker S=E$, and so $S=0_{L\left(E,F\right)}$. Hence, the condition (ii) is satisfied.

Next, it is clear that if $S\in L_{oc}\left(E,F\right)$, then $\left.S\right|_{H}\in L^{E}_{ob}\left(H,F\right)$. By part (ii) of Proposition \ref{sdsc} the embedding $J$ of $H$ into $E$ is supremum-continuous. According to part (i) of Proposition \ref{4.10} it is then order continuous, and so $\left.S\right|_{H}=SJ\in L^{E}_{oc}\left(H,F\right)$. Thus, the condition (i) is satisfied.

We conclude that the restriction to $H$ is an order isomorphism from $L_{oc}\left(E,F\right)$ into $L_{oc}^{E}\left(H,F\right)$, which are Dedekind complete vector lattices where the lattice operations and suprema of increasing nets are given by Riesz-Kantoronich formulas.
\end{proof}

The following result may be viewed as a version of Veksler's (see \cite[Theorem 1.65]{ab}) theorem as well as \cite[Theorem 1.56]{ab} for pervasive Archimedean OVS's.

\begin{theorem}\label{veksler}
The restriction to $E$ is an order isomorphism from $L_{oc}\left(E^{\delta},F\right)$ onto $L_{oc}\left(E,F\right)$. In fact, for $T\in L_{ob}\left(E,F\right)$ the following conditions are equivalent:
\item[(i)] $T$ is order continuous;
\item[(ii)] If $e_{\alpha}\downarrow 0_{E}$, then $Te_{\alpha}\too 0_{F}$;
\item[(iii)] If $Q\downarrow 0_{E}$, then $\bigwedge\left|TQ\right|= 0_{F}$;
\item[(iv)] $T^{\pm}$ exist in $L\left(E,F\right)$ and are order continuous;
\item[(v)] $\left|T\right|$ exists in $L\left(E,F\right)$ and is order continuous;
\item[(vi)] $T$ is a restriction of a (unique) order continuous operator from $E^{\delta}$ to $F$.
\end{theorem}
\begin{proof}
The first claim together with (i)$\Leftrightarrow$(vi) follows from Theorem \ref{veksler0}.

(i)$\Rightarrow$(ii)$\Rightarrow$(iii) and (iv)$\Rightarrow$(v) are straightforward. (iii)$\Rightarrow$(iv) follows from Lemma \ref{rkv} and part (ii) of Proposition \ref{addit}. (v)$\Rightarrow$(i) follows from the fact (Proposition \ref{bsc0} combined with part (i) of Proposition \ref{sdsc}) that $L_{oc}\left(E,F\right)$ is full in $L\left(E,F\right)$.
\end{proof}

Note that $T\in L_{oc}\left(E,F\right)_{+}$ has a unique positive extension to $E^{\delta}$, since this extension has to be order continuous, by virtue of Proposition \ref{ocodm}. The following example shows that even a zero functional on $E$ can have a non-positive order bounded extension to $E^{\delta}$, and so the extension is not unique in the class of order bounded functionals.

\begin{example}\label{example}
Let $E:=c$, so that $E^{\delta}=\ell_{\8}$. Note that the latter is a Banach lattice, and $E$ is a norm closed sublattice. Let $G:=E^{\delta}\slash E$ (a Banach space quotient), and let $Q:E^{\delta}\to G$ be the quotient map. Take $0_{G^{*}}\ne \nu\in G^{*}$. Then $Q^{*}\nu\in \ell_{\8}^{*}$ is a non-zero order bounded functional, whose restriction to $E$ is zero. The only order continuous extension of $\left.Q^{*}\nu\right|_{E}$ to $E^{\delta}$ is zero, which means that $Q^{*}\nu$ is not order continuous.
\qed\end{example}

Consider a slightly different version of Veksler's theorem, where a version of $E$-boundedness is a part of the conclusion rather than assumption.

\begin{proposition}\label{veksler2}
Assume that $E$ is a vector lattice. Let $T:H\to F$ be an order continuous operator. Then, $G:=\left\{e\in E,~ T\left[0_{E},\left|e\right|\right]_{H}\mbox{ is order bounded}\right\}$ is an ideal of $E$ which contains $H$. Moreover, $T$ uniquely extends to an order continuous operator from $G$ into $F$. This is a maximal order bounded extension of $T$ within $E$.
\end{proposition}
\begin{proof}
It is clear that $G$ is a solid set which contains $H$ and is stable under scalar multiplication. Let $e,g\in G_{+}$ so that there are $f,h\in F_{+}$ with $-f\le T\left[0_{E},e\right]_{H}\le f$ and $-h\le T\left[0_{E},g\right]_{H}\le h$. Then, $T\left[0_{E},e\right]_{H}+T\left[0_{E},g\right]_{H}\subset f+h-F_{+}$, and so using Proposition \ref{ardp} and arguing as in the proof of part (ii) of Proposition \ref{rdp}, one can show that $T\left[0_{E},e+g\right]_{H}\le f+h$. Analogously, $T\left[0_{E},e+g\right]_{H}\ge - f-h$ thus proving $e+g\in G$. We conclude that $G$ is an ideal, and $T$ is $G$-bounded.

Existence and uniqueness of an extension of $T$ to $G$ follows from Theorem \ref{veksler0}. It is clear that if $e$ is included in an order bounded extension of $T$, then $e\in G$. This shows maximality of the obtained extension.
\end{proof}

An example of a construction described in the preceding result is the so-called $T$-universal completion considered in \cite{klw}. We now present a number of results which do not require Dedekind completeness of the target space.

\begin{corollary}\label{veksler1}
Let $G$ be a directed Archimedean ordered vector space. For $T\in L_{ob}\left(E,G\right)$ the following conditions are equivalent:
\item[(i)] $T$ is order continuous;
\item[(ii)] If $e_{\alpha}\downarrow 0_{E}$, then $Te_{\alpha}\too 0_{G}$;
\item[(iii)] $T$ is a restriction of an order continuous operator from $E^{\delta}$ to $G^{\delta}$.\medskip

Moreover, the extension in (iii) is unique, and it is positive if and only if $T$ is. If $G$ is a lattice, then the conditions above are equivalent to the fact that if $Q\downarrow 0_{E}$, then $\bigwedge\left|TQ\right|= 0_{G}$.
\end{corollary}
\begin{proof}
(i)$\Rightarrow$(ii) is straightforward. (iii)$\Rightarrow$(i) follows from the fact that order convergence on $E$ and $G$ are the restrictions of order convergences on $E^{\delta}$ and $G^{\delta}$, respectively.

For (ii)$\Rightarrow$(iii) first note that the condition (ii) remains valid if we view $T$ as an operator from $T$ into $G^{\delta}$. It then follows from Theorem \ref{veksler} that $T$ is a restriction of a unique order continuous operator $T^{\delta}:E^{\delta}\to G^{\delta}$. Moreover, $T^{\delta}$ is positive if and only if $T$ is.

It is also clear that (ii) implies the last condition. That condition implies (iii) by the same argument as in (ii)$\Rightarrow$(iii).
\end{proof}

\begin{corollary}
Assume that $G$ is a monotone complete directed ordered vector space. Any monotone continuous $T\in L\left(E,G\right)$ has a unique monotone continuous extension $T^{\delta}:E^{\delta}\to G$.
\end{corollary}
\begin{proof}
First, $G$ is Archimedean (see \cite[Corollary 2.21(i)]{brtgww}). Therefore, according to Corollary \ref{veksler1} and part (i) of Proposition \ref{sdsc} there is a unique monotone continuous extension $T^{\delta}:E^{\delta}\to G^{\delta}$ of $T$, and so we only need to show that $T^{\delta}E^{\delta}\subset G$. Let $h\in E^{\delta}_{+}$. According to the comment after Proposition \ref{4.4} there is $P\subset \left[0_{E^{\delta}},h\right]_{E}$ such that $P\uparrow h$. As $T^{\delta}$ is monotone continuous, we get $T^{\delta}h=\bigvee T^{\delta}P=\bigvee TP$. Since by Corollary \ref{monoton} $G$ is monotone closed in $G^{\delta}$ we conclude that $T^{\delta}h\in G$.
\end{proof}

We remind the overarching assumptions of this subsection that $E$ is a pervasive Archimedean ordered vector space and $H\subset E$ is an order dense linear subspace.

\begin{proposition}\label{rucom}
Assume that $G$ is a relatively uniformly complete directed ordered vector space and $H$ is majorizing in $E$ and such that $\overline{H}^{\mathrm{ru}}=E$. Then, every $T\in L_{oc}\left(H,G\right)$ has a unique order bounded extension to $E$. This extension is order continuous.
\end{proposition}
\begin{proof}
According to Corollary \ref{veksler1} there is an order continuous extension $T^{\delta}:H^{\delta}\to G^{\delta}$. Since $H$ is order dense and majorizing in $E$, it is supremum-dense, and so $E\subset H^{\delta}$. It follows that $\widehat{T}:=\left.T^{\delta}\right|_{E}$ is an order continuous extension of $T$. Note that order continuous operators are order bounded, and so relatively uniformly continuous. Since ru convergence on $G$ is the restriction of ru convergence on $G^{\delta}$, and $G$ is relatively uniformly complete, it follows that it is relatively uniformly closed in $G^{\delta}$. Hence, $\widehat{T}^{-1}G$ is a relatively uniformly closed subset of $E$ which contains $H$. Thus, $\widehat{T}^{-1}G=E$, and so $\widehat{T}$ is an order continuous operator from $E$ into $G$. If $S:E\to G$ is an order bounded extension of $T$, then arguing as in the proof of Theorem \ref{veksler0} one can show that $S=\widehat{T}$.
\end{proof}

\begin{remark}
Note that a relatively uniformly dense subspace of a pervasive Archimedean space is not necessarily order dense: take polynomials as a subspace of $\Co\left[0,1\right]$.
\qed\end{remark}

\begin{corollary}\label{rucompl}
Assume that $H$ is majorizing in $E$ such that $\overline{H}^{\mathrm{ru}}=E$. If $F$ is a directed Archimedean ordered vector space, then $T\in L_{ob}\left(E,F\right)$ is order continuous if and only if $\left.T\right|_{H}$ is.
\end{corollary}
\begin{proof}
As in other results, note that there is a unique order bounded operator from $E$ to $F^{\delta}$ which extends $\left.T\right|_{H}$. Therefore, this extension must agree with $T$, and by Proposition \ref{rucom} it is order continuous.
\end{proof}

\subsection{Uo-continuous operators}\label{s6}

In this subsection we will extend the concept of the universal completion to directed Archimedean OVS's. More details on the topic see in Subsection \ref{sbuc}. In order to prove some universal properties we first need to obtain few more extension results.

Let $G$ be a vector lattice. We say that a net $\left(g_{\alpha}\right)_{\alpha\in A}\subset G$ \emph{unbounded order (uo)} converges to $g\in G$ (denoted $g_{\alpha}\touo g$) if $h\wedge\left|g_{\alpha}-g\right|\too 0_{G}$, for every $h\in G_{+}$. It is easy to see that $g_{\alpha}\too g$ if and only if $g_{\alpha}\touo g$ and $\left(g_{\alpha}\right)_{\alpha\in A}$ is eventually order bounded. If $H\subset G$ is an order dense sublattice, then uo convergence on $H$ is the restriction to $H$ of uo convergence on $G$ (this follows from e.g. \cite[Theorem 5.3]{erz} and the fact that every order dense sublattice is \emph{regular}, i.e. embedded in a supremum-continuous way).

Recall that $G$ is called \emph{universally complete} if it is Dedekind complete, and any collection of disjoint elements in $G$ is order bounded. Every Archimedean vector lattice $G$ has a \emph{universal completion}, i.e. an essentially unique universally complete vector lattice $G^{u}$ which contains $G$ as an order dense sublattice. Note that if $H$ contains $G$ as an order dense sublattice, then $H^{u}=G^{u}$, in particular $G^{\delta u}=G^{u}$.

An Archimedean vector lattice is universally complete if and only if it is uo-complete (see \cite[Theorem 17]{azouzi}). We can now improve this result slightly (a finite dimensional directed Archimedean non-lattice shows that pervasiveness in the following characterization cannot be dropped).

\begin{proposition}
An ordered vector space is a universally complete vector lattice if and only if it is pervasive Archimedean and uo-complete.
\end{proposition}
\begin{proof}
Necessity is clear. For sufficiency, note that by Corollary \ref{ruur} uo-completeness implies order completeness, and so according to \cite[Proposition 3.11]{brtgww} $E$ is a vector lattice. Now universal completeness follows from the comment above.
\end{proof}

Throughout the subsection $E,F$ are directed Archimedean ordered vector spaces. Define the \emph{universal completion} $E^{u}$ of $E$ to be the universal completion of $E^{\delta}$, and define the \emph{unbounded order (uo) convergence} on $E$ to be the restriction of uo convergence on $E^{\delta}$. Note that if $G$ is a vector lattice such that $E\subset G\subset E^{\delta}$, then $G^{u}=E^{\delta u}$, and the uo convergence there is the restriction of that on $E^{\delta}$. Hence, in the definition of both the universal completion and uo convergence $E^{\delta}$ may be replaced with any such $G$. In particular, the definition of universal completion is consistent with the one given in \cite{stennder}. Also, note that the embedding of $E$ into $E^{u}$ is supremum-continuous as a composition of supremum continuous embeddings of $E$ into $E^{\delta}$ and of $E^{\delta}$ into $E^{u}$.

\begin{question}[Cf. Theorem 3.7.17 in \cite{kv}]
Find an intrinsic description of uo convergence on $E$, i.e. given just in terms of $E$.
\end{question}

We denote the space of uo-continuous linear operators from $E$ into $F$ by $L_{uo}\left(E,F\right)$; also let $L_{uob}\left(E,F\right):=L_{uo}\left(E,F\right)\cap L_{ob}\left(E,F\right)$.

\begin{proposition}\label{uocont}
$L_{uob}\left(E,F\right)\subset L_{oc}\left(E,F\right)$. Every supremum-continuous operator is uo-continuous.
\end{proposition}
\begin{proof}
Let $T\in L_{uob}\left(E,F\right)$ and assume that $e_{\alpha}\too e$. Then $e_{\alpha}\touo e$ and $\left(e_{\alpha}\right)_{\alpha\in A}$ is eventually order bounded. Since $T\in L_{uo}\left(E,F\right)$ it follows that $Te_{\alpha}\touo Te$. As $T\in L_{ob}\left(E,F\right)$, it follows that $\left(Te_{\alpha}\right)_{\alpha\in A}$ is eventually order bounded. We conclude that $Te_{\alpha}\too Te$. Thus, $T\in L_{oc}\left(E,F\right)$.

Assume that $T\in L\left(E,F\right)$ is supremum-continuous. By Theorem \ref{supcont1} there is a supremum-continuous operator $T^{\delta}:E^{\delta}\to F^{\delta}$ which extends $T$. Such operator is uo-continuous (see \cite[Theorem 5.2]{erz}). Hence its restriction $T$ is also uo-continuous.
\end{proof}

\begin{theorem}\label{uocont2}
Assume that $E$ is pervasive and $F$ is a universally complete vector lattice. Then every $T\in L_{uob}\left(E,F\right)$ has a unique order bounded uo-continuous extension to $E^{u}$.
\end{theorem}
\begin{proof}
First, we show that $T$ admits a unique order continuous extension, and then show its uo-continuity. The proof of the first claim is based on the fact (see \cite[Theorem 6.4]{bctv}) that a set $A\subset E^{\delta}$ is order bounded in $E^{u}$ is and only if it is uo-bounded. By Proposition \ref{uocont} $T$ is order continuous, and so according Theorem \ref{veksler0} it is enough to show that $L_{uob}\left(E,F\right)\subset L^{E^{u}}\left(E,F\right)$. For every $g\in E^{u}_{+}$ we have that $\left[0_{E^{u}},g\right]_{E}$ is  uo-bounded in $E^{\delta}$, hence in $E$. As $T$ is uo-continuous, it preserves uo-bounded sets, and so $T\left[0_{E^{u}},g\right]_{E}$ is uo-bounded, thus order bounded in $F^{u}=F$.\medskip

To prove uo-continuity of the extension $T^{u}$ we use the filter approach to convergences (see \cite{bctv}). Assume that $\Fo$ is uo-null in $E^{u}$, and without loss of generality we may assume that $\Fo$ has a solid base. Let $\left.\Fo\right|_{E}:=\left\{A\cap E,~ A\in\Fo\right\}$. Since uo convergence on $E$ is the restriction of uo convergence on $E^{u}$, it follows that $\left.\Fo\right|_{E}\touo 0_{E}$, hence $T\left.\Fo\right|_{E}\touo 0_{F}$, and so $\left|T\left.\Fo\right|_{E}\right|\wedge f\too 0_{F}$, for every $f\in F_{+}$. Note that $\left|T\left.\Fo\right|_{E}\right|\wedge f$ is the filter on $F$ generated by the sets of the form $\rho_{f}\left(T\left(A\cap E\right)\right)$, where $A\in\Fo$, and $\rho_{f}:F\to F_{+}$ is a map defined by $\rho_{f}\left(g\right):=\left|g\right|\wedge f$. Hence, $B\in\left|T\left.\Fo\right|_{E}\right|\wedge f$ iff there is $A\in \Fo$ such that $\rho_{f}\left(T\left(A\cap E\right)\right)\subset B$.  Moreover, since $\Fo$ has a solid base, in the last assertion $A$ may be assumed solid.

For $f\in F_{+}$ the fact that $\left|T\left.\Fo\right|_{E}\right|\wedge f\too 0_{F}$ means that there is $P_{f}\subset F$ with $\bigwedge P_{f}=0_{F}$ such that $\left[-p,p\right]\in \left|T\left.\Fo\right|_{E}\right|\wedge f$, for every $p\in P_{f}$. In other words, for every $p\in P_{f}$ there is a solid $A_{f,p}$ such that $\rho_{f}\left(T\left(A_{f,p}\cap E\right)\right)\subset \left[-p,p\right]$, or equivalently  $T\left(A_{f,p}\cap E\right)\subset B_{f,p}$, where $B_{f,p}:=\left\{g\in F,~ \left|g\right|\wedge f\le p\right\}$. Note that $B_{f,p}$ is order closed, for every $f,p\in F_{+}$.

Since $E$ is order dense in $E^{u}$, according to part (iv) of Proposition \ref{sdsc} for every $e\in E_{+}^{u}$ we have that $\left[-e,e\right]_{E}$ is dense in $\left[-e,e\right]$ with respect to order convergence on $E^{u}$. It follows that for every solid $A\subset E^{u}$ we have that $A\cap E$ is dense in $A$ with respect to order convergence on $E^{u}$. For every $f\in F_{+}$ and $p\in P_{f}$ the set $A_{f,p}$ is solid, and $T^{u}\left(A_{f,p}\cap E\right)=T\left(A_{f,p}\cap E\right)\subset B_{f,p}$, and since $T^{u}$ is  is order continuous, and $B_{f,p}$ is order closed, it follows that $T^{u}A_{f,p}\subset\overline{T^{u}\left(A_{f,p}\cap E\right)}^{1}_{o}\subset B_{f,p}$. Arguing as before, we see that $\left|T^{u}\Fo\right|\wedge f\too 0_{F}$, for every $f\in F_{+}$, hence $T^{u}\Fo\touo 0_{F}$.
\end{proof}

If $F$ is merely directed Archimedean applying Theorem \ref{uocont2} to $J\circ T$, where $J:F\to F^{u}$ is the inclusion, yields the following.

\begin{corollary}\label{uocont3}
Assume that $E$ is pervasive. Then every $T\in L_{uob}\left(E,F\right)$ has a (unique) extension $T^{u}\in L_{uob}\left(E^{u},F^{u}\right)$.
\end{corollary}

\begin{remark}
It follows that if $E$ is pervasive and $F$ is universally complete, then $L_{uob}\left(E,F\right)$ is the set of restrictions of members of $L_{uob}\left(E^{u},F\right)$ to $E$. Note however, that this class may be strictly smaller than the restrictions of members of $L_{oc}\left(E^{u},F\right)$ to $E$ (see Example \ref{onuo}).\medskip

Consider the case when $F=\R$, so that order and uo convergences agree on $F$. Then every uo-continuous operator is order continuous, hence order bounded, and so $L_{uo}\left(E,\R\right)=L_{uob}\left(E,\R\right)$. Moreover, $L_{uo}\left(E,\R\right)$ is the linear span of the coordinate functionals (see \cite[Proposition 2.2]{glx}), which are supremum-continuous. Thus, $L_{uo}\left(E,\R\right)$ is the span of supremum-continuous operators from $E$ into $\R$.\medskip

On the other hand, it is possible for both inclusions $L_{uob}\left(E,F\right)\subset L_{uo}\left(E,F\right)$ and $L_{uob}\left(E,F\right)\subset L_{oc}\left(E,F\right)$ to be strict, even when both $E$ and $F$ are Dedekind complete Banach lattices.

Consider the Krengel's operator $T$ on $\ell_{2}$, defined in \cite[Example 5.6]{ab}, where it is proven that this operator has no modulus in $L\left(\ell_{2},\ell_{2}\right)$. As $\ell_{2}$ is a Dedekind complete vector lattice, it follows from e.g. part (ii) of Proposition \ref{rdp} that every operator in $L_{ob}\left(\ell_{2},\ell_{2}\right)$ has a modulus. We conclude that $T$ is not order bounded. To see that $T$ is uo-continuous, first observe that by \cite[Proposition 4.10]{bctv} uo convergence agrees with coordinatewise convergence on $\ell_{2}$. It can be easily seen from the construction of $T$ that it is coordinatewise continuous. This operator is in fact defined on $\R^{\N}=\ell_{2}^{u}$ and is order bounded there. This shows that the converse to Corollary \ref{uocont3} is false.

Since $\ell_{2}$ is an order continuous Banach lattice, it follows that $L_{oc}\left(\ell_{2},\R\right)=\ell_{2}^{*}$ can be identified with $\ell_{2}$ in the usual way. However, any $\nu\in L_{uob}\left(\ell_{2},\R\right)$ has an order continuous extension to $\ell_{2}^{u}=\R^{\N}$. It is easy to see that functionals with this property can be identified with $c_{00}\varsubsetneq \ell_{2}$.
\qed\end{remark}

\begin{example}\label{onuo}
Let $X:=\N_{0}^{\N}$ and let $Y:=\bigcup\limits_{n\in\N}\N^{n}_{0}$ -- the set of finite sequences of natural numbers. For $n\in\N$ and $y\in \N^{n}_{0}$ let $X_{y}$ be the set of all elements of $x$ whose first $n$ coordinates agree with $y$. For $z\in X$ let $X^{z}:=\left\{x\in X,~ x\ge z\right\}$. We endow $X$ with the so-called Hechler topology, i.e. the topology with a basis of clopen sets of the form $X^{z}_{y}:=X^{z}\cap X_{y}$, where $z\in X$ and $y\in Y$. Note that this collection of sets indeed forms a clopen basis for a Tychonoff Baire topology on $X$ (see \cite[Lemma 1.2]{lr}).

Fix $z\in X$. For $n\in\N$ let $z^{n}:=\1_{\left[n,+\8\right)}z$ -- $z$ with nullified first $n-1$ coordinated. Let $U_{z}$ be the set of all $x\in X$ which are eventually larger than $z$, i.e. $U_{z}:=\bigcup\limits_{n\in\N}X^{z^{n}}$. Clearly, $U_{z}$ is open. Note also that $U_{z}$ is dense in $X$. Indeed, for every $n\in\N$, $y\in \N_{0}^{n}$ and $w\in X$ such that $X_{y}^{w}\ne\varnothing$ we have $u\vee z^{n+1}\in X_{y}^{w}\cap X^{z^{n+1}}\subset X_{y}^{w}\cap U_{z}$, for any $u\in X_{y}^{w}$.

Let $S\left(X\right)$ be the space of densely defined continuous functions on $X$ (see \cite[Section 2]{wickstead}). For $n\in\N$ let $p_{n}:X\to\N_{0}$ be the projection on the $n$-th coordinate. Let us show that for every $\left(r_{n}\right)_{n\in\N}\subset\R_{+}$ the element $\sum\limits_{n\in\N}r_{n}2^{-p_{n}}:=\bigvee\limits_{n\in\N}\sum\limits_{k=1}^{n}r_{k}2^{-p_{k}}$ is defined in $S\left(X\right)$. Let $z\in X$ be such that $z\left(n\right)\ge \log_{2}\left(r_{n}\right)+n$, for every $n\in\N$. Then, for every $n\in\N$, $m\ge n$ and $x\in X^{z^{n}}$, we have $p_{m}\left(x\right)\ge p_{m}\left(z\right)\ge \log_{2}\left(r_{m}\right)+m$, hence $r_{m}2^{-p_{m}\left(x\right)}\le 2^{-m}$. It follows that for every $n\in\N$ the series $\sum\limits_{m\ge n}r_{m}2^{-p_{m}}$ converges absolutely on $X^{z^{n}}$, hence its sum is a continuous function, and so $\sum\limits_{m\in\N}r_{m}2^{-p_{m}}$ is also a continuous function on $X^{z^{n}}$. We conclude that the series converges to a continuous function on a dense open set $U_{z}$. It is easy to see that this function must be $\bigvee\limits_{n\in\N}\sum\limits_{k=1}^{n}r_{k}2^{-p_{k}}$.

Let $T:\R^{\N}\to S\left(X\right)$ be defined by $T\left(r_{n}\right)_{n\in\N}:=\sum\limits_{n\in\N}r_{n}2^{-p_{n}}$, for $\left(r_{n}\right)_{n\in\N}\subset\R_{+}$, and extended by the Kantorovich extension lemma to the entire $\R^{\N}$. Since the pointwise topology is complete metrizable locally solid and order continuous, it follows from part (iii) of Proposition \ref{locfull} that every positive operator into any Archimedean vector lattice is order continuous: if $P\downarrow\0$ in $\R^{\N}$ it contains a relatively uniformly null sequence, which is mapped by $T$ to a sequence with infimum zero, hence $\bigwedge TP=0_{S\left(X\right)}$. Let us show however, that $T$ is not uo-continuous.

For a finite $A\subset \N$ let $e_{A}:=\left|A\right|\1_{\N\backslash A}\in\R^{\N}$. Since uo convergence on $\R^{\N}$ is the coordinatewise convergence, it follows that $e_{A}\touo\0$. We now prove that $\1\wedge Te_{A}\not\too \0$, since for every finite $A\subset\N$ we have that $\bigvee\limits_{B\supset A}\left(\1\wedge Te_{B}\right)=\1$. Indeed, for $x\in X$ find $n>A$ and $m\ge 2^{x\left(n\right)}$; then $\left[Te_{A\cup \left[n+1,n+m\right]}\right]\left(x\right)\ge m2^{-x\left(n\right)}\ge 1$. We conclude that $Te_{A}\not\touo \0$ in $S\left(X\right)$, and so $Te_{A}\not\touo \0$ in $S\left(X\right)^{u}$. Thus, $T$ is not uo-continuous.
\qed\end{example}

\section{Completions of a directed Archimedean ordered vector space}\label{secot}

In this section we add two more completions of directed Archimedean OVS, on top of the Dedekind and Riesz completions. We then provide a categorical overview of the topic.

\subsection{The universal completion}\label{sbuc}

In this section we develop the concept of the universal completion of a directed Archimedean ordered vector space introduced in Subsection \ref{s6} a bit further. Throughout the section $E,F$ are directed Archimedean ordered vector spaces. We start with an auxiliary result. Part (i) of it is essentially \cite[Lemma 5.2(ii)]{imhoff}.

\begin{lemma}\label{lem}
Assume that $E,F$ are vector lattices, let $H\subset E$ be an order dense subspace and let $T\in L\left(E,F\right)$. Then:
\item[(i)] If $T$ is a homomorphism and $\left.T\right|_{H}$ is injective, then $T$ is also injective.
\item[(ii)] If $T$ is order continuous, and $T\left|h\right|=\left|Th\right|$, for every $h\in H$, then $T$ is a homomorphism.
\item[(iii)] If $T$ is order continuous, and $Tg\wedge Th=0_{F}$, for every $g,h\in H$ such that $g\wedge h=0_{E}$, then $T$ is a homomorphism.
\end{lemma}
\begin{proof}
(i): Assume that $e\in E$ is such that $Te=0_{F}$. Then, $T\left|e\right|=\left|Te\right|=0_{F}$, and so for every $h\in H\cap \left[0_{E},\left|e\right|\right]$ we also have $Th=0_{F}$. As $\left.T\right|_{H}$ is injective, it follows that $H\cap \left[0_{E},\left|e\right|\right]=\left\{0_{E}\right\}$, and since $H$ is order dense we conclude that $e=0_{E}$.\medskip

(ii): Let $\Phi:E\to F$ be defined by $\Phi\left(e\right):=T\left|e\right|-\left|Te\right|$, for $e\in E$. Since $T$ is order continuous, the same is true for $\Phi$, and so $\ker\Phi$ is order closed. By assumption $H\subset \ker\Phi$, and since by part (iv) of Proposition \ref{sdsc} the latter is order dense, hence dense with respect to order convergence, we conclude that $\ker\Phi=E$. It follows that $T\left|e\right|=\left|Te\right|$, for every $e\in E$, and so $T$ is a homomorphism.\medskip

(iii): Assume that $e,f\in E$ are such that $e\wedge f=0_{E}$. By Proposition \ref{4.4} there are $P\subset \left[0_{E}, e\right]_{H}$ and $Q\subset \left[0_{E}, f\right]_{H}$ such that $P\uparrow e$ and $Q\uparrow f$. Then, $p\wedge_{E} q=0_{E}$, for every $p\in P$ and $q\in Q$, hence $Tp\wedge_{F}Tq=0_{F}$. By order continuity of $T$ we get $Te=\bigvee TP$ and $Tf=\bigvee TQ$. We conclude that $Te\wedge Tf =0_{F}$. Since $e,f$ were arbitrary, $T$ is a homomorphism (see \cite[Theorem 2.14]{ab}).
\end{proof}

Part (ii) of the following result extends \cite[Theorems 7.17 and 7.19]{ab0}).

\begin{theorem}\label{supcont4}
$T\in L\left(E,F\right)$ is supremum-continuous if and only if there is an order continuous homomorphism $T^{u}:E^{u}\to F^{u}$ which extends $T$. This extension is unique, and it is injective if and only if $T$ is an order embedding.
\end{theorem}
\begin{proof}
Sufficiency follows from Theorem \ref{supcont1} and the fact that the embedding of $E$ into $E^{u}$ is supremum-continuous. Assume that $T\in L\left(E,F\right)$ is supremum-continuous. By Theorem \ref{supcont1} there a supremum-continuous operator $T^{\delta}:E^{\delta}\to F^{\delta}$ which extends $T$. By Proposition \ref{uocont} we have $T^{\delta}\in L_{uob}\left(E^{\delta},F^{\delta}\right)$, and so by Corollary \ref{uocont3} there is an order continuous extension $T^{u}:E^{u}\to F^{u}$ of $T^{\delta}$. The fact that $T^{u}$ is a homomorphism follows from part (ii) of Lemma \ref{lem}.\medskip

Let us prove the uniqueness of extension. Assume that an order continuous homomorphism $S:E^{u}\to F^{u}$ extends $T$. Then, $\left.S\right|_{E^{\delta}}$ is an order continuous homomorphism into $F^{u}$ which extend $T$. Viewing $T$ as a supremum-continuous operator from $E$ into $F^{u}$, by Theorem \ref{supcont1} there is a unique order continuous homomorphism $T^{\delta}:E^{\delta}\to F^{u}$, which extends $T$. It follows that $T^{\delta}=\left.S\right|_{E^{\delta}}$. Using uniqueness clause in Theorem \ref{uocont2} yields $S=T^{u}$.\medskip

By part (i) of Lemma \ref{lem}, $T^{u}$ is injective iff $T^{\delta}$ is, which in turn is equivalent to the fact $T$ is an order embedding, by virtue of Theorem \ref{supcont1}.
\end{proof}

\begin{proposition}\label{uo}
Assume that $H\subset E$ is a subspace such that the embedding from $H$ into $E$ is supremum continuous. Then, uo convergence on $H$ is the restriction of uo convergence on $E$.
\end{proposition}
\begin{proof}
According to Theorem \ref{supcont1} the embedding of $H$ into $E$ extends to an injective order continuous homomorphism from $H^{\delta}$ into $E^{\delta}$. Hence, we may view $H^{\delta}$ as a regular sublattice of $E^{\delta}$. By e.g. \cite[Theorem 5.3]{erz} uo convergence on $H^{\delta}$ is the restriction of uo convergence on $E^{\delta}$. Thus, the result follows from the definition of uo convergence on $H$ and $E$.
\end{proof}

\begin{remark}
The converse to Proposition \ref{uo} is false even under various additional assumptions. For example, let $F=\Co\left[0,1\right]$ which is a vector lattice, and let $E$ be the space of affine functions on $\left[0,1\right]$. Then, $E$ is order isomorphic to $\R^{2}$, and so it is a Dedekind complete vector lattice. Both uo convergences and order convergence on $E$ and $F^{\delta}$ are the usual topology on $\R^{2}$. However, the embedding of $E$ into $F$ does not preserve even absolute value of elements, and so is definitely not supremum continuous.
\qed\end{remark}

Let us now present an intrinsic characterization of the universal completion.

\begin{proposition}\label{univintr}
Let $F$ be a universally complete vector lattice which contains $E$. Then, there is an isomorphism from $F$ onto $E^{u}$ which preserves $E$ if and only if $E$ is supremum-dense in $F_{+}$. Moreover, if $E$ is pervasive, then under these conditions $E$ is order dense in $F$.
\end{proposition}
\begin{proof}
Necessity: We may assume that $F=E^{u}=E^{\delta u}$. Let $f\in F_{+}$. Since $E^{\delta}$ is order dense in $E^{u}$, according to e.g. \cite[Proposition 3.1]{brtgww} we have that $f=\bigvee \left[0_{F}, f\right]_{E^{\delta}}$. For every $g\in E^{\delta}$ we have $g=\bigvee g^{\Downarrow_{E}}$ (supremum in $E^{\delta}$, hence in $E^{u}$ since the embedding is supremum-continuous). We have that $$f=\bigvee \left[0_{F}, f\right]_{E^{\delta}}=\bigvee\limits_{g\in \left[0_{F}, f\right]_{E^{\delta}}}\bigvee g^{\Downarrow_{E}}=\bigvee\left(\bigcup\limits_{g\in \left[0_{F}, f\right]_{E^{\delta}}} g^{\Downarrow_{E}}\right).$$ It is left to show that $\bigcup\limits_{g\in \left[0_{F}, f\right]_{E^{\delta}}} g^{\Downarrow_{E}}=f^{\Downarrow_{E}}$. Let $e\in E$. Clearly, if there is $g\in \left[0_{F}, f\right]_{E^{\delta}}$ such that $e\le g$, then $e\le g\le f$, hence $e\in f^{\Downarrow_{E}}$. Assume that $e\in f^{\Downarrow_{E}}$. If $e\le 0_{F}$, then $e\in 0_{F}^{\Downarrow_{E}}$. Otherwise, $g:=e^{+}\in \left[0_{F}, f\right]_{E^{\delta}}$ fulfills $e\in g^{\Downarrow_{E}}$.\medskip

Sufficiency: First, by part (ii) of Proposition \ref{sdsc} the inclusion $J:E\to F$ is supremum-continuous. It is easy to show that $E$ is supremum-dense in the ideal $H$ generated by $E$ in $F$. Then $H$ is a Dedekind complete vector lattice, and so it is isomorphic to $E^{\delta}$. Let $R:H\to E^{\delta}$ be the corresponding order isomorphism.

It is clear that $H$ is supremum-dense in $F_{+}$, which implies that it is order dense by part (iii) of Proposition \ref{sdsc}. Then $F$ is a universally complete vector lattice which contains $H$ as an order dense sublattice, and so $F=H^{u}$. Thus, by Theorem \ref{supcont4}, $R$ may be extended to an injective order continuous homomorphism $R^{u}:F\to E^{\delta u}=E^{u}$. At the same time $R^{-1}$ can also be extended, to $R^{-1u}:E^{u}\to F$. It is easy to show by supremum-continuity that these extensions are inverses of each other, which means that $R^{u}$ is an order isomorphism.
\end{proof}

We now refine Theorem \ref{supcont1} to the case of pervasive spaces. Recall that order continuous homomorphisms between vector lattices are precisely supremum-continuous linear maps (see e.g. \cite[Theorem 5.2]{erz}). We now extend this result to operators from a pervasive Archimedean OVS.

\begin{theorem}\label{supcont2}
Let $E$ be a pervasive Archimedean OVS. For an operator $T:E\to F$ the following conditions are equivalent:
\item[(i)] $T$ is supremum-continuous;
\item[(ii)] $T$ is a uo continuous pre-homomorphism;
\item[(iii)] $T$ is an order continuous pre-homomorphism;
\item[(iv)] $T$ is a $\uparrow$-continuous pre-homomorphism;
\item[(v)] $T$ is bounded supremum-continuous.\medskip

Moreover, the extensions $T^{\delta}$ from Theorem \ref{supcont1} and $T^{u}$ from Theorem \ref{supcont4} are injective if and only if $T$ is.
\end{theorem}
\begin{proof}
First, (i)$\Rightarrow$(v) is trivial. (i)$\Rightarrow$(ii)$\Rightarrow$(iii) follows from Proposition \ref{uocont}. (iii)$\Rightarrow$(iv) follows from parts (ii) and (iii) of Proposition \ref{4.10}. The last claim follows from part (i) of Lemma \ref{lem}.


(iv)$\Rightarrow$(i): First, by part (i) of Proposition \ref{sdsc} $T$ is order continuous. According to Proposition \ref{riesz} there is a vector lattice homomorphism $T^{\rho}:E^{\rho}\to F^{\rho}$ which extends $T$. Since $E\subset E^{\rho} \subset E^{\delta}$, we have that $E$ is order dense and majorizing in $E^{\rho}$. As $\left.\widetilde{T}\right|_{E}=T$ is order continuous, $T^{\rho}$ is also order continuous, according to Proposition \ref{ocodm}. By the comment before the proposition it follows that $T^{\rho}$ is supremum continuous. It now follows from Theorem \ref{supcont1} that $T$ is supremum-continuous.\medskip

(v)$\Rightarrow$(i): First note that $T$ is monotone continuous, and so it is positive and order continuous by part (i) of Proposition \ref{sdsc}. Hence, by Theorem \ref{veksler}, there is an order continuous $T^{\delta}:E^{\delta}\to F^{\delta}$ which extends $T$. The restriction of $T^{\delta}$ to $E$ clearly satisfies the conditions of part (iii) of Lemma \ref{lem}, hence $T^{\delta}$ is a homomorphism. By the comment before the theorem $T^{\delta}$ is supremum-continuous, and so the same is true for $T$, according to Theorem \ref{supcont1}.
\end{proof}

\begin{remark}\label{bs}
It is clear that the implications (i)$\Rightarrow$(v) and (i)$\Rightarrow$(ii)$\Rightarrow$(iii)$\Rightarrow$(iv) in Theorem \ref{supcont2} do not require the assumption of pervasiveness. However, (ii)+(v)$\Rightarrow$(i) fails in general, even when $\dim E<\8$ and $F=\R$, as we will see in \cite{bb}. In the same article we construct an uo-continuous pre-homomorphism which fails to be a strict pre-homomorphism. An example of a directed order continuous strict pre-homomorphism which fails to be order continuous is also presented there. In Example \ref{bscnr} we construct a bounded supremum continuous functional which fails to be a pre-homomorphism (and hence fails to be supremum continuous). While bounded supremum continuity implies directed order continuity, it does not imply order continuity, as Example \ref{bscnoc} demonstrates.

Let $E$ and $\nu$ be as in Example \ref{onuor}, in particular $\nu$ is order continuous pre-homomorphism. Let us show that it fails to be uo-continuous. Indeed, for $m\in\N$ let $h_{m}$ be defined by $h_{m}\left(n\right):=0$, for $n<m$, $h_{m}\left(m\right):=2^{m}-1$, $h_{m}\left(n\right):=1$, for $n>m$. Then, $E\ni h_{m}\touo\0$, but $\nu\left(h_{m}\right)=1\not\to 0$.

The last claim in Theorem \ref{supcont2} also fails for non-pervasive spaces. Indeed, in \cite[Example 5.4]{imhoff} the author has constructed a bijective supremum-continuous operator, which is not an order isomorphism, hence not an order embedding, and so by Theorem \ref{supcont1} $T^{\delta}$ is not injective.
\qed\end{remark}

\begin{remark}
It follows from Theorem \ref{supcont2} that the converse to Proposition \ref{uo} is true if $H$ is pervasive, and the embedding into $E$ satisfies $A^{\Uparrow_{H}\Downarrow_{H}}\subset A^{\Uparrow\Downarrow}$, for every finite $A\subset E$ (such subspaces are studied in \cite{stennder}).
\qed\end{remark}

\subsection{The relative uniform completion}\label{sbru}

In this section we adapt some of the results presented in \cite{bt} to Archimedean ordered vector spaces. The proofs which are analogous to the vector lattice case are omitted. Throughout the section $E$ is a directed Archimedean ordered vector space, embedded into another directed Archimedean vector space $G$.

We define $E^{\mathrm{ru}}$ to be the relative uniform closure of $E$ in $E^{\delta}$, equivalently, the intersection of all relatively uniform subspaces $E\subset G\subset E^{\delta}$.

\begin{corollary}[Cf. Corollary 4.7, \cite{bt}]
If $G\subset E^{\mathrm{ru}}$, then $G^{\mathrm{ru}}=E^{\mathrm{ru}}$.
\end{corollary}

The following result is a direct consequence of Propositions \ref{ruex} and \ref{ruexp}.

\begin{corollary}[Cf. Theorem 5.3, \cite{bt}]
Let $F$ be a directed Archimedean relatively uniformly complete ordered vector space. Then, any positive / bi-positive / regular / supremum-continuous operator $T:E\to F$ has a unique extension $T^{\mathrm{ru}}:E^{\mathrm{ru}}\to F$ of the same type.
\end{corollary}

\begin{corollary}[Cf. Theorem 6.1 and Corollary 6.2, \cite{bt}]\label{intersection}
Assume that $G$ is relatively uniformly complete. Then, $E^{\mathrm{ru}}$ is isomorphic to the intersection of all relatively uniformly complete subspaces $E\subset H\subset G$. If additionally $E$ is majorizing in $G$, then $E^{\mathrm{ru}}$ is isomorphic to $\overline{E}^{\mathrm{ru}}$.
\end{corollary}

\begin{corollary}\label{emberu}
Assume that $E$ is majorizing in $G$ and such that $G=\overline{E}^{\mathrm{ru}}$. Then, $G$ embeds into $E^{\mathrm{ru}}$.
\end{corollary}

We now discuss a few more results about extensions of operators. The following is sketched in \cite[Remark 3.6]{bx}, but for the reader's convenience, some additional details are provided.

\begin{proposition}\label{ruriesz}
\item[(i)] $E^{\rho\mathrm{ru}}=\overline{\spa E^{+}}^{\mathrm{ru}}=\overline{\spa \left|E\right|}^{\mathrm{ru}}$, where the positive part, absolute value, linear span and relative uniform closure are taken in $E^{\delta}$.
\item[(ii)] If $F$ is a relatively uniformly complete vector lattice, and $T:E\to F$ is a pre-homomorphism, there is a unique vector lattice homomorphism $T^{\rho\mathrm{ru}}:E^{\rho\mathrm{ru}}\to F$ which extends $T$.
\end{proposition}
\begin{proof}
(i): Recall that $E^{\rho}$ is the vector lattice generated by $E$ in $E^{\delta}$. In fact, since $E=E_{+}-E_{+}=\spa E_{+}$, it follows that $E^{\rho}=\vlt E_{+}$. Moreover, since $E^{\rho\mathrm{ru}}$ is the relative uniform closure of $E^{\rho}$ in $E^{\delta}$, it follows that $E^{\rho\mathrm{ru}}=\overline{\vlt E_{+}}^{\mathrm{ru}}$. Now, applying \cite[Corollary 3.2]{bx} to $F=E^{\delta}$, $A=E_{+}$, and either positive part or absolute value for $\Phi$ yields $\overline{\vlt E_{+}}^{\mathrm{ru}}=\overline{\spa E^{+}}^{\mathrm{ru}}=\overline{\spa \left|E\right|}^{\mathrm{ru}}$.\medskip

(ii): By Proposition \ref{riesz} there is a vector lattice homomorphism $T^{\rho}:E^{\rho}\to F$ which extends $T$, while according to \cite[Theorem 5.3]{bt} there is a vector lattice homomorphism $T^{\rho\mathrm{ru}}:E^{\rho\mathrm{ru}}\to F$ which extends $T^{\rho}$, thus extends $T$. If $S:E^{\rho\mathrm{ru}}\to F$ is a vector lattice homomorphism which extends $T$, then $\ker\left(T^{\rho\mathrm{ru}}-S\right)$ is a relatively uniformly closed sublattice of $E^{\rho\mathrm{ru}}$ which contains $E$, hence $\ker\left(T^{\rho\mathrm{ru}}-S\right)=E^{\rho\mathrm{ru}}$, and so $S=T^{\rho\mathrm{ru}}$.
\end{proof}

We will refer to $E^{\rho\mathrm{ru}}$ as the \emph{ru-Riesz} completion of $E$.

\begin{example}
If $E$ is an order unit space with the designated unit $e$, and $K_{E}$ is the state space of $\left(E,e\right)$, endowed with the weak* topology, then $E^{\rho\mathrm{ru}}\simeq\Co\left(\overline{\ext K_{E}}\right)$ (see \cite[Section 2.5]{kv}).
\qed\end{example}

\begin{corollary}[Cf. Theorem 5.2, \cite{bt}]\label{ruexpre}
Assume that $G$ is directed, $E$ is majorizing in $G$ and such that $\overline{E}^{\mathrm{ru}}=G$. Let $F$ be a directed Archimedean ordered vector spaces, and let $T:G\to F$ be an order bounded operator. Then, $T$ is a pre-homomorphism if and only if $\left.T\right|_{E}$ is.
\end{corollary}
\begin{proof}
First, it follows from Corollary \ref{emberu} that $G\subset E^{\mathrm{ru}}\subset E^{\rho\mathrm{ru}}$. According to part (ii) of Proposition \ref{ruriesz} $\left.T\right|_{E}$ extends to a vector lattice homomorphism from $E^{\rho\mathrm{ru}}$ into $F^{\delta}$, whose restriction $S:G\to F^{\delta}$ is an order bounded operator, which extends $\left.T\right|_{E}$. By Proposition \ref{ruex} such an extension is unique, and so $T=S$. Hence, $T$ extends to a vector lattice homomorphism from $E^{\rho\mathrm{ru}}$ into $F^{\delta}$, which according to Proposition \ref{riesz} implies that $T$ is a pre-homomorphism.
\end{proof}

For $h\in E_{+}$ we denote $E_{h}:=\bigcup\limits_{r>0}\left[-rh,rh\right]$ endowed with the norm $\|\cdot\|_{h}$ which is the Minkowski functional of $\left[-h,h\right]$. Note that then the latter set is precisely the unit ball of $E_{h}$. Also, if $E\subset G$, then $E_{h}\subset G_{h}$ isometrically. We also denote the closure of $E_{h}$ in $G_{h}$ with respect to $\|\cdot\|_{h}$ by $\overline{E_{h}}$. It follows from Corollary \ref{ruh} that $E$ is relatively uniformly complete if and only if $E_{h}$ is norm complete, for every $h\in E_{+}$.

Also, let $G_{E}:=\bigcup\limits_{e\in E_{+}}\left[-e,e\right]$, which is the full hull of $E$, i.e. the smallest full subset (and subspace) of $G$ which contains $E$.

\begin{proposition}[Cf. Proposition 3.12, \cite{bt}]\label{loc-adh}
  For $g\in G$ the following conditions are equivalent:
  \item[(i)] $g$ is a relative uniform limit in $G$ of a sequence in $E$ with a regulator in $E$;
  \item[(ii)] There exists $h\in E_{+}$ such that $g\in \overline{E_{h}}$;
  \item[(iii)] $g$ is in the relative uniform adherence of $E$ in $G_{E}$.
\end{proposition}

It is easy to see that if $f\le h$ in $E_{+}$ then $\overline{E_{f}}\subset \overline{E_{h}}$.

\begin{corollary}[Cf. Corollary 3.13, \cite{bt}]\label{IuX-cofin}
  Let $A\subset E_{+}$ be majorizing, i.e. for every $e\in E$ there is $a\in A$ with $a\ge e$. Then, $\bigcup\limits_{h\in A}\overline{E_{h}}$ is the relative uniform adherence of $E$ in $G_{E}$. In particular, if $E_{h}$ is norm-closed in $G_{h}$, for every $h\in A$, then $E$ is relatively uniformly closed in $G_{E}$.
\end{corollary}

Note that if $rA\subset A$, for every $r>0$, then $\bigcup\limits_{h\in A}\overline{E_{h}}=\bigcup\limits_{h\in A}\left(\overline{E_{h}}\cap\left[-h,h\right]\right)$.

\begin{proposition}[Cf. Proposition 7.1, \cite{bt}]
  If $E$ is a subspace of a relatively uniformly complete space $G$ then $E^{\mathrm{ru}}$ is isomorphic to the relative uniform closure of $E$ in $G_{E}$.
\end{proposition}

Let $E^{\mathrm{ru}}_{1}$ denote the relative uniform adherence of $E$ in $E^{\mathrm{ru}}$.

\begin{proposition}[Cf. Proposition 7.2, \cite{bt}]\label{ru-IX}
If $E$ is a subspace of a relatively uniformly complete space $G$ then $E^{\mathrm{ru}}_{1}$ is isomorphic to the relative uniform adherence of $E$ in $G_{E}$.
\end{proposition}

Applying Corollary \ref{IuX-cofin} to $E$ in $E^{\mathrm{ru}}$ we get the following result.

\begin{corollary}\label{IuX-cofin1}
  Let $A\subset E_{+}$ be majorizing, i.e. for every $e\in E$ there is $a\in A$ with $a\ge e$. If $E_{h}$ is norm complete, for every $h\in A$, then $E$ is relatively uniformly complete.
\end{corollary}

The criterion above can be modified in presence of the RDP.

\begin{proposition}
Assume that $E$ has the RDP. Then:
\item[(i)] If $E$ has an order unit $h$, then it is norm complete if and only if whenever $\left(e_{n}\right)_{n\in\N}\subset E_{+}$ is such that $\sum\limits_{n\in\N}\|e_{n}\|_{h}<\8$, it then follows that there is $\bigvee\limits_{n\in\N}\sum\limits_{k=1}^{n}e_{k}$.
\item[(ii)] Let $A\subset E_{+}$ is such that there is not proper full directed subspace of $E$ containing $A$. If $E_{h}$ is norm complete, for every $h\in A$, then $E$ is relatively uniformly complete.
\end{proposition}
\begin{proof}
(i) is proven similarly to \cite[Theorem 1]{bb0}, but instead of using positive and negative parts, it relies on the observation that if $e\in E_{+}$ and $f\in \left[-e,e\right]$, there are $g,h\in \left[0_{E},e\right]$ with $f=g-h$. Indeed, just take $g:=\frac{e+f}{2}$ and $h:=\frac{e-f}{2}$.\medskip

(ii): Let $C$ be the set of all $h\in E_{+}$ such that $E_{h}$ is norm complete. Since relative uniform completeness is inherited by full subspaces, it follows that if $0_{E}\le g\le \in C$, then $g\in C$. Clearly $C$ is closed under scalar multiplication. Let us show that it is also closed under addition. Assume that $g,h\in C$ and let $\left(e_{n}\right)_{n\in\N}\subset E_{+}$ is such that $\sum\limits_{n\in\N}r_{n}<\8$, where $r_{n}:=\|e_{n}\|_{g+h}$, for $n\in\N$. It follows that for every $n\in\N$ we have $e_{n}\in r_{n}\left[0_{E},g+h\right]$, and so the RDP yields $f_{n}\in r_{n}\left[0_{E},g\right]$ such that $e_{n}-f_{n}\in r_{n}\left[0_{E},h\right]$. Then, $\sum\limits_{n\in\N}\|f_{n}\|_{g}<\8$ and $\sum\limits_{n\in\N}\|e_{n}-f_{n}\|_{h}<\8$, and so there are $\bigvee\limits_{n\in\N}\sum\limits_{k=1}^{n}f_{k}$ and $\bigvee\limits_{n\in\N}\sum\limits_{k=1}^{n}\left(e_{k}-f_{k}\right)$. It is easy to see that the sum of these two suprema is $\bigvee\limits_{n\in\N}\sum\limits_{k=1}^{n}e_{k}$. It now follows from (i) that $E_{g+h}$ is norm complete, and so $g+h\in C$. According to \cite[Corollary 4.3]{brtgww} that $C-C$ is a full directed subspace of $E$, hence our assumption gives $C-C$, and so $C=E_{+}$ (here we are applying \cite[Corollary 4.3]{brtgww} again). We conclude that $E$ is relatively uniformly complete.
\end{proof}

\begin{proposition}[Cf. Proposition 7.5, \cite{bt}]
Let $E$ be a subspace of $G$. The following conditions are equivalent:
\item[(i)] $E$ is majorizing, relatively uniformly dense in $G$ (in the sense of adherence) and every relatively uniformly Cauchy sequence in $E$ has a limit in $G$;
\item[(ii)] Every $\mathrm{ru}_{E}$-Cauchy sequence in $E$ has a $\mathrm{ru}_{G}$-limit in $G$, and every $g\in G$ is a $\mathrm{ru}_{G}$-limit of a $\mathrm{ru}_{E}$-Cauchy sequence in $E$;
\item[(iii)] $G$ is isomorphic to $Cauchy\left(E\right)\slash Null\left(E\right)$, where $Cauchy\left(E\right),Null\left(E\right)\subset E^{\N}$ consist of relatively uniformly Cauchy and relatively uniformly null sequences in $E$, respectively;
\item[(iv)] $G$ is isomorphic to $E^{\mathrm{ru}}_{1}$.
\end{proposition}
\begin{proof}
(iv)$\Rightarrow$(i)$\Rightarrow$(ii) is straightforward. (iii)$\Rightarrow$(iv) is proven similarly to \cite[Proposition 7.5]{bt}, except for bi-positivity of the operator $T$. Assume that $\left(e_{n}\right)_{n\in\N}\subset E$, $h\in E_{+}$ and $g\in G_{+}$ are such that $-\frac{1}{n}h\le e_{n}-g\le \frac{1}{n}h$. Let $f_{n}:=e_{n}+\frac{1}{n}h\ge g \ge 0_{G}$, for $n\in\N$. This sequence clearly defines a positive element in $Cauchy\left(E\right)\slash Null\left(E\right)$, equivalent to $\left(e_{n}\right)_{n\in\N}$, and so it is the preimage of $g$ with respect to $T$.\medskip

(ii)$\Rightarrow$(iii): Define $S:G\to Cauchy\left(E\right)\slash Null\left(E\right)$ by $g\mapsto \left(e_{n}\right)+Null\left(E\right)$, where $\left(e_{n}\right)$ is a $\mathrm{ru}_{E}$-Cauchy sequence in $E$ whose $\mathrm{ru}_{G}$-limit is $g$. It is easy to see that this map is a well-defined injective homomorphism. Its surjectivity follows immediately from the first assumption in (ii).
\end{proof}

\begin{proposition}
If $E$ has the RDP, then $E^{\mathrm{ru}}_{1},E^{\mathrm{ru}}$ have the RDP. 
\end{proposition}
\begin{proof}
It follows from Corollary \ref{IuX-cofin} that $\bigcup\limits_{h\in E_{+}}\overline{E_{h}}$, where each $\overline{E_{h}}$ is take in $E^{\delta}$. Proposition \ref{aip} guarantees that $\overline{E_{h}}$ has the RDP, for every $h\in E_{+}$. Since the RDP is equivalent to the interpolation property (\cite[Theorem 1.54]{at}), it is easy to see that a union of directed collection of subspaces of $E^{\delta}$ with the RDP has the RDP. This proves that $E^{\mathrm{ru}}_{1}$ has the RDP.

It is now easy to prove by transfinite induction that $E^{\mathrm{ru}}_{j}$ has the RDP, for every ordinal $j\le\omega_{1}$: if $j$ is a successor, then we have that $E^{\mathrm{ru}}_{j}=\left(E^{\mathrm{ru}}_{j-1}\right)^{\mathrm{ru}}_{1}$ has the RDP by the assumption of the induction and the first step; if $j$ is a limit ordinal the directed union argument does the job. Finally, we conclude that $E^{\mathrm{ru}}=E^{\mathrm{ru}}_{\omega_{1}}$ has the RDP.
\end{proof}

\begin{question}
Assume that $E$ has the ADP. Does $E^{\mathrm{ru}}$ have the RDP? In particular, if $E$ is relatively uniformly complete and has the ADP, does it then have the RDP?
\end{question}

The next characterization of $\left(E^{\rho}\right)^{\mathrm{ru}}_{1}$ is inspired by \cite[Theorem 2]{veksl}. Here $G$ is the universal completion of $E$ represented as $C^{\8}\left(K\right)$, for some compact Hausdorff extremally disconnected space $K$ (see \cite[Section 7.2]{ab0} for the details).

\begin{proposition}
Let $g\in C^{\8}\left(K\right)$. Then, $g\in \left(E^{\rho}\right)^{\mathrm{ru}}_{1}$ if and only if there is $h\in E_{+}$ such that $h\ge \pm g$ and such that whenever $x,y\in \overline{\supp} h$ are such that $\lim\limits_{z\to x}\left(-1\vee\frac{e\left(z\right)}{h\left(z\right)}\wedge 1\right)=\lim\limits_{z\to y}\left(-1\vee\frac{e\left(z\right)}{h\left(z\right)}\wedge 1\right)$, for every $e\in E$, it then follows that $\lim\limits_{z\to x}\frac{g\left(z\right)}{h\left(z\right)}=\lim\limits_{z\to y}\frac{g\left(z\right)}{h\left(z\right)}$.
\end{proposition}
\begin{proof}
First, recall that $F:=E^{\rho}$ is the vector lattice generated by $E$ in $E^{\delta}$. Since the latter is a sublattice of $E^{u}$, it follows that $F$ is the vector lattice generated by $E$ in $E^{u}=C^{\8}\left(K\right)$. Also, note that $E^{u}_{F}=E^{\delta}=F^{\delta}$, and so our goal is to describe the relative uniform adherence of $F$ in $E^{u}_{F}$. According to the comment after Corollary \ref{IuX-cofin} the latter is $\bigcup_{h\in E_{+}}\left(\overline{F_{h}}\cap\left[-h,h\right]\right)$. Hence, it is enough to show that $g\in\overline{F_{h}}$, for $h\in E$ with $h\ge\pm g$, if and only if whenever $x,y\in \overline{\supp} h$ are such that $\lim\limits_{z\to x}\left(-1\vee\frac{e\left(z\right)}{h\left(z\right)}\wedge 1\right)=\lim\limits_{z\to y}\left(-1\vee\frac{e\left(z\right)}{h\left(z\right)}\wedge 1\right)$, for every $e\in E$, it then follows that $\lim\limits_{z\to x}\frac{g\left(z\right)}{h\left(z\right)}=\lim\limits_{z\to y}\frac{g\left(z\right)}{h\left(z\right)}$.\medskip

For $f\in E^{u}$ define $\hat{h}\left(f\right):=-h\vee f\wedge h$. Let us show that $F_{h}$ is the vector lattice generated by $\hat{h}\left(E\right)$ in $E^{u}$. Clearly, $\hat{h}\left(E\right)\subset F$ and $\hat{h}\left(E\right)\subset E^{u}_{h}$, hence $\vlt \hat{h}\left(E\right)\subset F\cap E^{u}_{h}=F_{h}$. On the other hand, assume that $f\in F_{h}$, so that $-h\le \frac{1}{r}f\le h$, for some $r>0$. Since $\frac{1}{r}f\in F$, there are $\left(e_{ij}\right)_{i,j=1}^{m,n}\subset E$ such that $\frac{1}{r}f=\bigvee\limits_{i=1}^{m}\bigwedge\limits_{j=1}^{n}e_{ij}$. Then, $$\frac{1}{r}f=\hat{h}\left(\frac{1}{r}f\right)=-h\vee \bigvee\limits_{i=1}^{m}\bigwedge_{j=1}^{n}e_{ij}\wedge h=\bigvee\limits_{i=1}^{m}\bigwedge_{j=1}^{n}\left(-h\vee e_{ij}\wedge h\right)=\bigvee\limits_{i=1}^{m}\bigwedge_{j=1}^{n}\hat{h}\left(e_{ij}\right)\in\vlt \hat{h}\left(E\right),$$
and so $f=r \frac{1}{r}f\in \vlt \hat{h}\left(E\right)$.\medskip

Let $L:=\overline{\supp} h$. Consider the map $J:\Co\left(L\right)\to \Co^{\8}\left(K\right)$ as described in \cite[Examples 2.2 and 8.8]{bt}. This map is an isometric isomorphism onto the ideal generated by $h$ in $\Co^{\8}\left(K\right)$. It follows that $J^{-1}\left(\overline{F_{h}}\cap\left[-h,h\right]\right)=\overline{J^{-1}F_{h}}\cap\left[-\1,\1\right]=\overline{\vlt}J^{-1}\hat{h}\left(E\right)\cap\left[-\1,\1\right]$.
Note that $\1=J^{-1}$, and so $\overline{\vlt}J^{-1}\hat{h}\left(E\right)$ is a closed vector lattice generated by a collection of functions containing $\1$. Recall that if $\1\in A\subset \Co\left(L\right)$, then $f\in\overline{\vlt} A$ if and only if  $f\left(x\right)=f\left(y\right)$ whenever $x,y\in L$ are such that $a\left(x\right)=a\left(y\right)$, for every $a\in A$ (this follows from e.g. \cite[Corollary 2.2]{bt0}). Hence, $g\in \overline{F_{h}}\cap\left[-h,h\right]$ $\Leftrightarrow$ $J^{-1}g\in \overline{\vlt}J^{-1}\hat{h}\left(E\right)\cap\left[-\1,\1\right]$ $\Leftrightarrow$ $\left[J^{-1}g\right]\left(x\right)=\left[J^{-1}g\right]\left(y\right)$, whenever $x,y\in L$ are such that $J^{-1}\hat{h}\left(x\right)=J^{-1}\hat{h}\left(y\right)$, for every $e\in E$. Due to the nature of $J$, the latter condition translates to $\lim\limits_{z\to x}\frac{g\left(z\right)}{h\left(z\right)}=\lim\limits_{z\to y}\frac{g\left(z\right)}{h\left(z\right)}$, whenever $x,y\in L$ are such that $\lim\limits_{z\to x}\left(-1\vee\frac{e\left(z\right)}{h\left(z\right)}\wedge 1\right)=\lim\limits_{z\to y}\left(-1\vee\frac{e\left(z\right)}{h\left(z\right)}\wedge 1\right)$, for every $e\in E$.
\end{proof}

\begin{question}
Find a description of $E^{\rho\mathrm{ru}}$ in terms of $C^{\8}\left(K\right)$.
\end{question}

\subsection{Categorical interpretations}\label{cat}

Recall that a category $\mathcal{D}$ is a \emph{subcategory} of a category $\mathcal{C}$ if $\left|\mathcal{D}\right|\subset \left|\mathcal{C}\right|$ (we use this notation to denote the class of objects of a category), and for every $D_{1},D_{2}\in \left|\mathcal{D}\right|$ we have $\mathcal{D}\left(D_{1},D_{2}\right)\subset\mathcal{C}\left(D_{1},D_{2}\right)$ (this notation signifies the set of morphisms from $D_{1}$ into $D_{2}$). If the latter inclusion is an equality for all $D_{1},D_{2}$, we say that $\mathcal{D}$ is a \emph{full subcategory} of $\mathcal{C}$.

A full subcategory $\mathcal{D}$ of a category $\mathcal{C}$ is called \emph{reflective} if for every object $C\in \left|\mathcal{C}\right|$ there is $C^{\mathcal{D}}\in \left|\mathcal{D}\right|$ and a morphism $\iota_{C}\in\mathcal{C}\left(C,C^{\mathcal{D}}\right)$ such that for every $D\in\left|\mathcal{D}\right|$ and $\varphi\in \mathcal{C}\left(C,D\right)$ there is a unique $\varphi^{\mathcal{D}}\in\mathcal{D}\left(C^{\mathcal{D}},D\right)$ such that
$\varphi=\varphi^{\mathcal{D}}\circ\iota_{C}$. The standard ``diagram chasing'' argument shows that the pair $\left(C^{\mathcal{D}},\iota_{C}\right)$ is essentially unique. This pair is called the \emph{reflector} of $C$ in $\mathcal{D}$. If $\iota_{C}$ is clear from the context, we will simply refer to $C^{\mathcal{D}}$ as the reflector.\medskip

Let $\mathbf{DAOVS}_{\bigvee}$, $\mathbf{DAOVS}_{\uparrow}$, $\mathbf{DAOVS}_{oc}$, $\mathbf{DAOVS}_{+}$, $\mathbf{DAOVS}_{reg}$, $\mathbf{DAOVS}_{ob}$, $\mathbf{DAOVS}_{uob}$ and $\mathbf{DAOVS}_{pre}$ be the categories where the objects are directed Archimedean ordered vector spaces and the morphisms are $\bigvee$-continuous, $\uparrow$-continuous, order continuous, positive, regular, order bounded and uo-continuous order bounded operators, and pre-homomorphisms, respectively. Further, let $\mathbf{PAOVS}_{ob}$, $\mathbf{RAOVS}_{ob}$, $\mathbf{AVL}_{ob}$, $\mathbf{RVL}_{ob}$, $\mathbf{DVL}_{ob}$, $\mathbf{UVL}_{ob}$ stand for the full subcategories of $\mathbf{DAOVS}_{ob}$ whose objects are pervasive, relatively uniformly complete, lattices, relatively uniformly complete lattices, Dedekind complete lattices and universally complete lattices, respectively. Finally, $\mathbf{O}_{m}:=\mathbf{DAOVS}_{m}\cap \mathbf{O}_{ob}$, i.e. the full subcategory of $\mathbf{DAOVS}_{m}$, whose objects are the same as in $\mathbf{O}_{ob}$, where $\mathbf{O}\in \left\{\mathbf{DAOVS, PAOVS, RAOVS, AVL, RVL, DVL, UVL}\right\}$ and $m\in \left\{\bigvee, \uparrow, oc, +, reg, ob, uob, pre\right\}$. Note that in $\mathbf{AVL}:=\mathbf{AVL}_{pre}$ and $\mathbf{RVL}:=\mathbf{RVL}_{pre}$ the morphisms are vector lattice homomorphisms.

We can now summarize many results of the article in the following form:

\begin{corollary}\label{reflections}
\item[(i)] $\mathbf{DVL}_{\bigvee}$ and $\mathbf{UVL}_{\bigvee}$ are reflective in $\mathbf{DAOVS}_{\bigvee}$ with the Dedekind and universal completions being the reflectors.
\item[(ii)] $\mathbf{DVL}_{\uparrow}$ and $\mathbf{DVL}_{oc}$ are reflective in $\mathbf{PAOVS}_{\uparrow}$ and $\mathbf{PAOVS}_{oc}$, respectively, with the Dedekind completion being the reflector in both cases. $\mathbf{UVL}_{uob}$ is reflective in $\mathbf{PAOVS}_{uob}$ with the universal completion being the reflector.
\item[(iii)] $\mathbf{AVL}$ and $\mathbf{RVL}$ are reflective in $\mathbf{DAOVS}_{pre}$ with the Riesz and ru-Riesz completions being the reflectors.
\item[(iv)] $\mathbf{RAOVS}_{m}$ is reflective in $\mathbf{DAOVS}_{m}$, where $m\in \left\{\bigvee, +, reg, pre\right\}$ with the relative uniform completion being the reflector.
\end{corollary}
\begin{proof}
(i) follows from Theorems \ref{supcont1} and \ref{supcont4}. (ii) follows from Theorem \ref{veksler}, part (i) of Proposition \ref{sdsc}, and Theorem \ref{uocont2}. (iii) follows from Proposition \ref{riesz} and part (ii) of Proposition \ref{ruriesz}. (iv) follows from Propositions \ref{ruex}, \ref{ruexp} and Corollary \ref{ruexpre}.
\end{proof}

Throughout the section $E$ is a directed Archimedean OVS. Let $\iota_{E}^{\rho}:E\to E^{\rho}$, $\iota_{E}^{\mathrm{ru}}:E\to E^{\mathrm{ru}}$, $\iota_{E}^{\rho\mathrm{ru}}:E\to E^{\rho\mathrm{ru}}$, $\iota_{E}^{\delta}:E\to E^{\delta}$ and $\iota_{E}^{u}:E\to E^{u}$ be the usual embeddings.

\begin{proposition}
Let $F$ be a universally complete vector lattice, and let $\iota:E\to F$ be linear. The following conditions are equivalent:
\item[(i)] There is an isomorphism $J:F\to E^{u}$ such that $J\iota=\iota_{E}^{u}$;
\item[(ii)] $\iota$ is an order embedding and $\iota\left(E\right)$ is supremum-dense in $F_{+}$;
\item[(iii)] $\left(F,\iota\right)$ is the reflector of $E$ in $\mathbf{UVL}_{\bigvee}$.
\end{proposition}
\begin{proof}
(i)$\Rightarrow$(iii) follows immediately from Theorem \ref{supcont4}. The converse follows from uniqueness of the reflector. (i)$\Leftrightarrow$(ii) follows from Proposition \ref{univintr}. The last claim follows from part (iii) of Proposition \ref{sdsc}.
\end{proof}

Analogously we can characterize $E^{\delta}$ as follows.

\begin{proposition}
Let $F$ be a Dedekind complete vector lattice, and let $\iota:E\to F$ be linear. The following conditions are equivalent:
\item[(i)] There is an isomorphism $J:F\to E^{\delta}$ such that $J\iota=\iota_{E}^{\delta}$;
\item[(ii)] $\iota$ is an order embedding and $\iota\left(E\right)$ is supremum-dense in $F$;
\item[(iii)] $\left(F,\iota\right)$ is the reflector of $E$ in $\mathbf{DVL}_{\bigvee}$.
\end{proposition}

In the case of pervasive spaces, we get the following characterizations of the two completions.

\begin{proposition}
Assume that $E$ is pervasive. Let $F$ be a Dedekind complete vector lattice, and let $\iota:E\to F$ be linear. The following conditions are equivalent:
\item[(i)] There is an isomorphism $J:F\to E^{\delta}$ such that $J\iota=\iota_{E}^{\delta}$;
\item[(ii)] $\iota$ is an order embedding and $\iota\left(E\right)$ is order dense and majorizing in $F$;
\item[(iii)] $\left(F,\iota\right)$ is the reflector of $E$ in $\mathbf{DVL}_{oc}$;
\item[(iv)] $\left(F,\iota\right)$ is the reflector of $E$ in $\mathbf{DVL}_{\uparrow}$.
\end{proposition}

\begin{proposition}
Assume that $E$ is pervasive. Let $F$ be a universally complete vector lattice, and let $\iota:E\to F$ be linear. The following conditions are equivalent:
\item[(i)] There is an isomorphism $J:F\to E^{u}$ such that $J\iota=\iota_{E}^{u}$;
\item[(ii)] $\iota$ is an order embedding and $\iota\left(E\right)$ is order dense in $F$;
\item[(iii)] $\left(F,\iota\right)$ is the reflector of $E$ in $\mathbf{UVL}_{uob}$.
\end{proposition}

\begin{proposition}
Let $F$ be an Archimedean vector lattice, and let $\iota:E\to F$ be linear. The following conditions are equivalent:
\item[(i)] There is an isomorphism $J:F\to E^{\rho}$ such that $J\iota=\iota_{E}^{\rho}$;
\item[(ii)] $\iota$ is an order embedding and $\iota\left(E\right)$ is supremum-dense in $F$ and $\vlt\left(\iota\left(E\right)\right)=F$;
\item[(iii)] $\iota$ is a pre-homomorphism, order embedding and $\vlt\left(\iota\left(E\right)\right)=F$;
\item[(iv)] $\left(F,\iota\right)$ is the reflector of $E$ in $\mathbf{AVL}$.
\end{proposition}
\begin{proof}
(i)$\Rightarrow$(ii)$\Rightarrow$(iii) are straightforward. (i)$\Leftrightarrow$(iv) has been already proven in Corollary \ref{reflections}.

(iii)$\Rightarrow$(i): According to Proposition \ref{riesz} there is an injective vector lattice homomorphism $\iota^{\rho}:E^{\rho}\to F$ which extends $\iota$. Then, $\iota^{\rho}E^{\rho}=\vlt\left(\iota\left(E\right)\right)=F$, and so $\iota^{\rho}$ is surjective, hence an isomorphism, and so $J:=\iota^{\rho-1}$ satisfy the requirements.
\end{proof}

\begin{proposition}[cf. Proposition 6.11, \cite{bt}]\label{char}
Let $F$ be a relatively complete ordered vector space, and let $\iota:E\to F$ be linear. The following conditions are equivalent:
\item[(i)] There is an isomorphism $J:F\to E^{\mathrm{ru}}$ such that $J\iota=\iota_{E}^{\mathrm{ru}}$;
\item[(ii)] $\iota$ is an order embedding and $\iota\left(E\right)$ is majorizing in $F$ and $\overline{\iota\left(E\right)}^{\mathrm{ru}}=F$;
\item[(iii)] $\iota$ is an order embedding and there are no proper uniformly complete subspaces of $F$ which contain $\iota\left(E\right)$;
\item[(iv)] $\left(F,\iota\right)$ is the reflector of $E$ in $\mathbf{RAOVS}_{m}$, where $m\in \left\{\bigvee, +, reg, pre\right\}$.
\end{proposition}
\begin{proof}
(i)$\Rightarrow$(ii)$\Rightarrow$(iii) is straightforward. (iii)$\Rightarrow$(i) follows from Corollary \ref{intersection}. (i)$\Leftrightarrow$(iv) has been already proven in Corollary \ref{reflections}.
\end{proof}

\begin{proposition}
Let $E$ be an Archimedean ordered vector space, let $F$ be a relatively complete vector lattice, and let $\iota:E\to F$ be linear. The following conditions are equivalent:
\item[(i)] There is an isomorphism $J:F\to E^{\rho\mathrm{ru}}$ such that $J\iota=\iota_{E}^{\rho\mathrm{ru}}$;
\item[(ii)] $\iota$ is an order embedding and $\iota\left(E\right)$ is majorizing in $F$ and $\overline{\vlt\iota\left(E\right)}^{\mathrm{ru}}=F$;
\item[(iii)] $\iota$ is an order embedding and there are no proper uniformly complete sublattices of $F$ which contain $\iota\left(E\right)$;
\item[(iv)] $\left(F,\iota\right)$ is the reflector of $E$ in $\mathbf{RVL}$.
\end{proposition}

Let us present a partial characterization of $E^{\mathrm{ru}}$ which implements the idea of the completion as a complete object which has the original object as a dense subobject. In this case ``subobject'' is meant not only algebraically but also topologically.

\begin{proposition}[cf. Proposition 6.11, \cite{bt}]\label{char0}
Assume that $F\supset E$ is relatively uniformly complete. Then in the list below we have (i)$\Rightarrow$(ii)$\Rightarrow$(iii) and none of the implications are reversible:
\item[(i)] $F=\overline{E}_{1}^{\mathrm{ru}}$ and $\mathrm{ru}_{E}=\left.\mathrm{ru}_{F}\right|_{E}$;
\item[(ii)] $F\simeq E^{\mathrm{ru}}$, with the isomorphism preserving $E$;
\item[(iii)] $F=\overline{E}^{\mathrm{ru}}$ and $\mathrm{ru}_{E}=\left.\mathrm{ru}_{F}\right|_{E}$.
\end{proposition}
\begin{proof}
Note that according to Corollary \ref{intersection} the condition (ii) is equivalent to the fact that $E$ is majorizing in $G$ and $F=\overline{E}^{\mathrm{ru}}$. In the light of this, (ii)$\Rightarrow$(iii) is straightforward. For the failure of the converse see Example \ref{norucomp}. For (ii)$\not\Rightarrow$(i) see \cite[Section 9]{bt}.\medskip

(i)$\Rightarrow$(ii): It follows from \cite[Theorem 7.9]{bctv} (see also \cite[Remark 3.10]{bt}) that every set in $E$ which is order bounded in $F$ is order bounded in $E$. Let us show that $E$ is majorizing in $F=\overline{E}^{\mathrm{ru}}_{1}$. Let $f\in F$ so that there is a sequence $\left(e_{n}\right)_{n\in\N}\subset E$ which converges to $f$. Then this sequence is order bounded in $F$, hence in $E$, and so there is $e\in E$ such that $e_{n}\le e$, for every $n\in\N$. Passing to the limit yields $f\le e$.
\end{proof}

\begin{example}\label{norucomp}
Let $G:=\left\{f\in\Co\left[0,1\right],~ \forall n\in\N:~ f\left(\frac{1}{n}\right)=0\right\}$, which is an ideal in $\Co\left[0,1\right]$. Let $\Po\left[0,1\right]$ be the space of polynomials on $\left[0,1\right]$, which is a norm dense subspace of $\Co\left[0,1\right]$. Since a nonzero polynomial can have only finite number of roots, we have $\Po\left[0,1\right]\cap G=\left\{\0\right\}$.

Inside of the vector lattice $\Co\left[0,1\right]^{\N}$ consider the ideals $G^{\N}$ and $H:=\bigcup\limits_{n\in\N}\Co\left[0,1\right]^{n}$. Then, $F:=G^{\N}+H$ is also an ideal in $\Co\left[0,1\right]^{\N}$, and so it is a relatively complete Archimedean vector lattice (here we are using a simple observation that relative uniform completeness passes down to ideals). Let $E:=\bigcup\limits_{n\in\N}\Po\left[0,1\right]^{n}$, considered as a subspace of $H$, hence of $F$. It is not hard to show that $E$ is directed. It is also clear that $E$ is not majorizing in $F$, which implies that $F$ fails to satisfy the condition (ii) in Proposition \ref{char}.\medskip

We claim that for every $f\in F_{+}$ the set $\left[0_{F},f\right]_{E}$ is order bounded in $E$. Represent $f=\left(f_{n}\right)_{n\in\N}$, where $f_{n}\in\Co\left[0,1\right]$, and there is $m\in\N$ such that $f_{n}\in G$, for every $n\ge m$. Let $e\in \left[0_{F},f\right]_{E}$, so that $e=\left(e_{n}\right)_{n\in\N}$, where $e_{n}\in\Po\left[0,1\right]$ and $\0\le e_{n}\le f_{n}$, for every $n\in\N$. If $n\ge m$, it then follows that $e_{n}\in G$, hence $e_{n}=\0$. For $n<m$ let $u_{n}:=\|f_{n}\|_{\8}\1\in\Po\left[0,1\right]$, and let $u:=\left(u_{1},...,u_{m-1},\0,\0,...\right)\in E$. For $n=1,...,m-1$ we have $e_{n}\le f_{n}\le u_{n}$, while for $n\ge m$ we have $e_{n}=\0= u_{n}$. Hence, $e\le u$, and since $e$ was arbitrary, we conclude that $\left[0_{F},f\right]_{E}$ is order bounded in $E$.

Next, we show that $\overline{E}^{\mathrm{ru}}_{2}=F$. Since $\Po\left[0,1\right]$ is norm dense in $\Co\left[0,1\right]$, it is easy to see that $H\subset \overline{E}^{\mathrm{ru}}_{1}$. It is enough to show that $F=\overline{H}^{\mathrm{ru}}_{1}$. Let $f=\left(f_{n}\right)_{n\in\N}\in F$, and let $m\in\N$ be such that $f_{n}\in G$, for every $n\ge m$. For $n\in\N$ let $h_{n}:=\left(f_{1},...,f_{n},\0,\0,...\right)\in H$ and $v_{n}:=n\left|f_{n}\right|$. Note that $v_{n}\in G$, for every $n\ge m$, hence $v:=\left(v_{n}\right)_{n\in\N}\in F$. Since $\left|f-h_{n}\right|\le \frac{1}{n+1}v$, for every $n\in\N$, we conclude that $h_{n}\toru f$, as required.

It is now clear that $\overline{E}^{\mathrm{ru}}=F$, but also $\mathrm{ru}_{E}=\left.\mathrm{ru}_{F}\right|_{E}$ (see \cite[Proposition 3.9 or Remark 3.10]{bt}), and so $E$ and $F$ satisfy the condition (iii) in Proposition \ref{char0}.
\qed\end{example}

\subsection{The relative uniform completion as a direct limit}\label{sublim}

A \emph{directed system} in a category $\mathcal{C}$ is an analogue of a monotone net. Namely, it is a pair $\left(\left(C_{\alpha}\right)_{\alpha\in A},\left(\varphi_{\alpha,\beta}\right)_{\alpha\le\beta}\right)$ where $C_{\alpha}\in\left|\mathcal{C}\right|$ and $\varphi_{\alpha,\beta}\in \mathcal{C}\left(C_{\alpha},C_{\beta}\right)$, for every $\alpha\le\beta$, such that $\varphi_{\beta,\gamma}\circ \varphi_{\alpha,\beta}=\varphi_{\alpha,\gamma}$, for every $\alpha\le\beta\le\gamma$.

The \emph{direct limit} of $\left(\left(C_{\alpha}\right)_{\alpha\in A},\left(\varphi_{\alpha,\beta}\right)_{\alpha\le\beta}\right)$ is a pair $\left(C,\left(\iota_{\alpha}\right)_{\alpha\in A}\right)_{\alpha\in A}$, where $C\in\left|\mathcal{C}\right|$ and $\iota_{\alpha}\in \mathcal{C}\left(C_{\alpha},C\right)$, for every $\alpha$, such that $\iota_{\alpha}=\iota_{\beta}\circ \varphi_{\alpha,\beta}$, for every $\alpha\le\beta$, and whenever $\left(C',\left(\iota'_{\alpha}\right)_{\alpha\in A}\right)_{\alpha\in A}$ is such that $\iota'_{\alpha}=\iota'_{\beta}\circ\varphi_{\alpha,\beta}$, for every $\alpha\le\beta$, there is a unique $\iota:C\to C'$ such that $\iota'_{\alpha}=\iota\circ \iota_{\alpha}$, for every $\alpha\in A$.

Direct limit may be regarded as the analogue of the extremum of a monotone net. We note that often (but not always) if objects in $\mathcal{C}$ are structured sets, morphisms are structured maps, and there is a set $X$ which contains all $C_{\alpha}$'s and $\varphi_{\alpha,\beta}$'s are inclusion maps, then $C=\bigcup\limits_{\alpha\in A}C_{\alpha}$ (as a set).

Assume that $\left(\left(C_{\alpha}\right)_{\alpha\in A},\left(\varphi_{\alpha,\beta}\right)_{\alpha\le\beta}\right)$ and $\left(\left(D_{\alpha}\right)_{\alpha\in A},\left(\psi_{\alpha,\beta}\right)_{\alpha\le\beta}\right)$ are directed systems in $\mathcal{C}$ and for every $\alpha\in A$ there is $\theta_{\alpha}:C_{\alpha}\to D_{\alpha}$ such that $\theta_{\beta}\varphi_{\alpha,\beta}=\psi_{\alpha,\beta}\theta_{\alpha}$, for every $\alpha\le\beta$. Assume further that $\left(C,\left(\iota_{\alpha}\right)_{\alpha\in A}\right)_{\alpha\in A}$ and $\left(D,\left(\kappa_{\alpha}\right)_{\alpha\in A}\right)_{\alpha\in A}$ are the direct limits of the two systems. It is not hard to show that there is $\theta:C\to D$ such that $\theta\iota_{\alpha}=\kappa_{\alpha}\theta_{\alpha}$, for every $\alpha\in A$.\medskip

Now let $E$ be a directed Archimedean OVS considered as a subspace of $E^{\delta}$. For every $e\le f$, let $\varphi_{ef}:E_{e}\to E_{f}$ be the inclusion map. Let $\overline{\varphi}_{ef}:\overline{E_{e}}\to \overline{E_{f}}$ also be the inclusions. It is clear that $\overline{\varphi}_{fg}\circ\overline{\varphi}_{ef}=\overline{\varphi}_{eg}$, for every $e\le f\le g$. Note that $\overline{E_{e}}\simeq E_{e}^{\mathrm{ru}}$ and $\overline{\varphi}_{ef}=\varphi_{ef}^{\mathrm{ru}}$, for every $0_{E}\le e\le f$. Since $E_{+}$ is a directed set, it follows that $\left(\left(\overline{E_{e}}\right)_{e\in E_{+}},\left(\overline{\varphi}_{ef}\right)_{e\le f}\right)$ is a directed system in $\mathbf{DAOVS}_{+}$. We leave it to the reader to verify that in this category the direct limit of a directed system where the transition morphisms are inclusions is the union. This yields the following result.

\begin{proposition}\label{limit}
\item[(i)] $E^{\mathrm{ru}}_{1}$ is the direct limit in $\mathbf{DAOVS}_{+}$ of the system described above.
\item[(ii)] If $E$ is pervasive then the limit can also be taken in $\mathbf{PAOVS}_{+}$, $\mathbf{PAOVS}_{pre}$, $\mathbf{PAOVS}_{\uparrow}$ or  $\mathbf{PAOVS}_{\bigvee}$ (and the same for lattices).
\item[(iii)] $E^{\mathrm{ru}}$ is the direct limit of $\left(E^{\mathrm{ru}}_{j}\right)_{j<\omega_{1}}$ in $\mathbf{DAOVS}_{+}$ $\mathbf{DAOVS}_{pre}$ or $\mathbf{DAOVS}_{\bigvee}$. If $E$ is pervasive, the the limit can also be taken in $\mathbf{PAOVS}_{\uparrow}$.
\end{proposition}
\begin{proof}
(i) is justified by the discussion above and Corollary \ref{IuX-cofin}.

(ii): Let us first consider the case of $\mathbf{PAOVS}_{pre}$. Let $F$ be a pervasive Archimedean ordered vector space and assume that a pre-homomorphism $T_{e}:\overline{E_{e}}\to F$ is given, for every $e\in E_{+}$, such that whenever $e\le f$ we have that $Te=Tf\circ\overline{\varphi}_{ef}$. Note that the latter condition simply means that $Te=\left.Tf\right|_{\overline{E_{e}}}$. It is then easy to see that there is a map $T:E^{\mathrm{ru}}_{1}\to F$ which agrees with $T_{e}$ on $\overline{E_{e}}$, for every $e\in E_{+}$. It is also not hard to show that $T$ is a positive linear operator, and so it is left to show that it is a pre-homomorphism. In fact, by virtue of Corollary \ref{ruexpre}, it is enough to show that $\left.T\right|_{E}$ is a pre-homomorphism.

Fix $e\in E$ and let $J_{e}:E_{e}\to \overline{E_{e}}$ be the inclusion. Note that since $\overline{E_{e}}\simeq E_{e}^{\mathrm{ru}}$, it follows that $E$ is supremum-dense in $\overline{E_{e}}$, and so $J_{e}$ is supremum-continuous. Therefore, $$\left.\left(\left.T\right|_{E}\right)\right|_{E_{e}}=\left.T\right|_{E_{e}}=\left.\left(\left.T\right|_{\overline{E_{e}}}\right)\right|_{E_{e}}=\left.T_{e}\right|_{E_{e}}=T_{e}J_{e}$$ is a pre-homomorphism, as a composition of two pre-homomorphisms. It follows now from \cite[Theorem 5.6 and Proposition 4.2(a)]{stennder} that $\left.T\right|_{E}$ is a pre-homomorphism, as required.\medskip

We now consider the case of $\mathbf{PAOVS}_{\uparrow}$. We start with the same configuration with $T_{e}$ being $\uparrow$-continuous, for every $e\in E_{+}$, and we construct $T:E^{\mathrm{ru}}_{1}\to F$ as above. Again, it is enough to show that $\left.T\right|_{E}$ is $\uparrow$-continuous, this time by virtue of Proposition \ref{ocodm}. Let $Q\subset E_{+}$ be such that $Q\downarrow 0_{E}$ and pick $q_{0}\in Q$. Since the latter set is directed, we have that $Q\cap\left[0_{E},q_{0}\right]\downarrow 0_{E}$, and in particular $Q\cap\left[0_{E},q_{0}\right]\downarrow 0_{E}$ in $E_{q_{0}}$. Since $J_{q_{0}}$ is supremum-continuous, and $T_{q_{0}}$ is $\uparrow$-continuous, it follows that $T\left(Q\cap\left[0_{E},q_{0}\right]\right)=T_{q_{0}}J_{q_{0}}\left(Q\cap\left[0_{E},q_{0}\right]\right)\downarrow 0_{F}$, which implies $\bigwedge TQ=0_{F}$.\medskip

Finally, the case of $\mathbf{PAOVS}_{\bigvee}$ follows similar footsteps: by the previous two steps $T$ is an $\uparrow$-continuous pre-homomorphism, hence it is supremum-continuous, by virtue of Theorem \ref{supcont2}.\medskip

(iii) is proven analogously: we take advantage of the fact that the restriction of the limiting operator to $E$ has to have the required continuity and that $E$ is relatively uniformly dense in $E^{\mathrm{ru}}$.
\end{proof}

\begin{remark}
Note that $\varphi_{ef}$ may fail to be a pre-homomorphism, or even preserve finite suprema (see \cite[Examples 3.14, 3.16 and 4.3(1)]{stennder}). It may also fail to be $\uparrow$-continuous: take $E$ as in Example \ref{odensen} and note that $E=E_{f+\1}$, while $F=E_{\1}$. Due to these pathologies we can only consider $\mathbf{DAOVS}_{+}$ in part (i) of Proposition \ref{limit}.
\qed\end{remark}

\begin{remark}
Assume that $E$ is a lattice. If $\overline{\varphi}_{ef}$ were interval preserving for every $e\le f$, then the direct limit in (i) would be already relatively uniformly complete, and hence equal to $E^{\mathrm{ru}}$ (see \cite[Proposition 2.14]{filter}; see also \cite[Proposition 3.14]{bt} for the argument which does not use directed limits). However, \cite[Example 6.13]{bt} shows that this assumption does not hold even when $E^{\mathrm{ru}}=E^{\mathrm{ru}}_{1}$, and furthermore \cite[Section 9]{bt} shows that $E^{\mathrm{ru}}\ne E^{\mathrm{ru}}_{1}$, in general.
\qed\end{remark}

\begin{remark}
Proposition \ref{limit} shows that we can construct $E^{\mathrm{ru}}$ ``completely categorically''. It is not hard to show that $\mathbf{DAOVS}_{+}$ has direct limits (see also \cite{vv} for various categories of vector lattices). Start with $E_{0}:=E$, and define $E_{1}$ as the direct limit of $\left(\left(\overline{E_{e}}\right)_{e\in E_{+}},\left(\overline{\varphi}_{ef}\right)_{e\le f}\right)$. Note that $\overline{\varphi}_{ef}:\overline{E_{e}}\to \overline{E_{f}}$ can be defined by extending the inclusion from $E_{e}$ into $E_{f}$ by continuity, and does not refer to the embedding of $\overline{E_{e}}$ and $\overline{E_{f}}$ into $E^{\delta}$. Also, since $E$ is the limit of $\left(\left(E_{e}\right)_{e\in E_{+}},\left(\varphi_{ef}\right)_{e\le f}\right)$, there is an embedding from $\iota_{0}:E_{0}\to E_{1}$, as described before Proposition \ref{limit}. This construction can be iterated. In particular, if $j$ is a limit ordinal, we construct $E_{i}$ as the direct limit of $\left(E_{i}\right)_{i<j}$. In the end, $E^{\mathrm{ru}}=E_{\omega_{1}}$.
\qed\end{remark}

\section{Some counterexamples}\label{sbce}

First, we show that the additional hypothesis is not superfluous in Propositions \ref{rdp}.

\begin{example}\label{adp}
In this example we denote $\ell_{\8}:=\ell_{\8}\left(\Z\right)$. Let $$G:=\left\{f\in\ell_{\8},~\forall i=1,2,3,4:~\exists \lambda_{i}\left(f\right):=\lim\limits_{n\to\8}f\left(4n+i\right) ~\&~ \exists \lambda_{0}\left(f\right):=\lim\limits_{n\to\8}f\left(-n\right)\right\},$$ which is a sublattice of $\ell_{\8}$. Note that $G$ is stable with respect to alteration of finite number of coordinates. Also, denote $\lambda:=\lambda_{1}+\lambda_{2}-\lambda_{3}-\lambda_{4}$.


Let $j:\Z\to\R$ be $\I_{\Z}^{+}$, i.e. $j\left(n\right)=n$, for $n\in\N$ and $j\left(n\right)=0$, for $n\in\Z\backslash \N$, and let $F:=G+\R j$. Define $\mu:F\to\R$ by $\mu\left(f\right):=\lim\limits_{n\to\8}\frac{f\left(n\right)}{n}$. It is easy to see that $\mu$ is well defined and $\ker\mu=G$. Note that if $f=g+rj$, where $g\in G$, then $r=\mu\left(f\right)$, and $g=Pf:=f-\mu\left(f\right)j$. Let $$E:=\left\{f\in F,~ \lambda\left(Pf\right)=0,~ \lambda_{0}\left(f\right)=\mu\left(f\right)+\sum\limits_{i=1}^{4}\lambda_{i}\left(Pf\right) \right\}.$$ Note that since finitely supported functions belong to $E$, the latter is an order dense subspace of $\R^{\N}$, hence it is pervasive. It is not hard to show that $j\vee \1$ is a strong unit in $E$ with $\left\|f\right\|_{j\vee \1}=\bigvee\limits_{n\in\Z}\frac{\left|f\left(n\right)\right|}{n\vee 1}$.\medskip

Let us show that for every $e\in E_{+}$ there are $e',e''\in \left[\0,e\right]$ such that $e=e'+e''$, $\lambda_{0}\left(e'\right)= \lambda_{0}\left(e\right)$ and $\mu\left(e''\right)=\mu\left(e\right)$ (and then $\mu\left(e'\right)=0=\lambda_{0}\left(e''\right)$). If $\mu\left(e\right)=0$, put $e':=e$ and $e'':=\0$. Otherwise define $e':\Z\to\R$ by $e'\left(n\right):=e\left(n\right)$, for $n\le 0$ and $e'\left(n\right):=e\left(n\right)\wedge\frac{1}{4}\lambda_{0}\left(e\right)$, for $n\in\N$. Then, since $e\left(n\right)\to+\8$, we have $e'\left(n\right)=\frac{1}{4}\lambda_{0}\left(e\right)$, eventually, and so $\lambda_{i}\left(e'\right)=\frac{1}{4}\lambda_{0}\left(e\right)=\frac{1}{4}\lambda_{0}\left(e'\right)$, for every $i=1,2,3,4$, justifying $e'\in E$. It is clear that $e'\in \left[\0,e\right]\cap \ker\mu$.\medskip

We now prove that $E$ satisfies the ADP. Let $u,v\in E_{+}$ and let $w\in\left[\0,u+v\right]$. Then, $\lambda_{0}\left(w\right)\le \lambda_{0}\left(u\right)+\lambda_{0}\left(v\right)$ and $\mu\left(w\right)\le \mu\left(u\right)+\mu\left(v\right)$, and so it is easy to find $a,b,c,d\in\left[0,1\right]$ such that $\lambda_{0}\left(w\right)= a\lambda_{0}\left(u\right)+b\lambda_{0}\left(v\right)$ and $\mu\left(w\right)=c\mu\left(u\right)+d\mu\left(v\right)$. Then, we in fact have $\lambda_{0}\left(w\right)=\lambda_{0}\left(au'+cu''+bv'+dv''\right)$ and $\mu\left(w\right)=\mu_{0}\left(au'+cu''+bv'+dv''\right)$. Let $e:=au'+cu''\in \left[\0,u\right]$ and $f:=bv'+dv''\in \left[\0,v\right]$.

Fix $n\in\N$. Let $e_{n},f_{n}:\Z\to\R$ be defined by $e_{n}\left(m\right):=w\left(m\right)\wedge u\left(m\right)$ and $f_{n}\left(m\right):=\left(w\left(m\right)-u\left(m\right)\right)^{+}$, for $\left|m\right|\le n$ and $e_{n}\left(m\right):=e\left(m\right)$ and $f_{n}\left(m\right):=f\left(m\right)$ otherwise. Since altering finitely many coordinates does not affect containment in $E$, we have $e_{n},f_{n}\in E$, with $\lambda_{0}\left(e_{n}\right)=\lambda_{0}\left(e\right)$, $\lambda_{0}\left(f_{n}\right)=\lambda_{0}\left(f\right)$, $\mu\left(e_{n}\right)=\mu\left(e\right)$, $\mu\left(f_{n}\right)=\mu\left(f\right)$. Moreover coordinatewise inequalities yield $e_{n}\in \left[\0,u\right]$ and $f_{n}\in \left[\0,v\right]$. Since $e_{n}+f_{n}$ agrees with $w$ on $\left[-n,n\right]$ and with $e+f$ outside, it follows that  $$\|w-e_{n}-f_{n}\|_{j\vee\1}=\bigvee\limits_{m<-n}\left|w\left(m\right)-e\left(m\right)-f\left(m\right)\right|\vee\bigvee\limits_{m>n}\frac{1}{m}\left|w\left(m\right)-e\left(m\right)-f\left(m\right)\right|.$$
Note that $\lim\limits_{m\to\8} w\left(-m\right)=\lambda_{0}\left(w\right)=\lambda_{0}\left(e+f\right)=\lim\limits_{m\to\8} \left( e\left(-m\right)+f\left(-m\right)\right)$ and $\lim\limits_{m\to\8} \frac{1}{m} w\left(m\right)=\mu\left(w\right)=\mu\left(e+f\right)=\lim\limits_{m\to\8} \frac{1}{m}\left( e\left(m\right)+f\left(m\right)\right)$. Hence, $e_{n}+f_{n}\to w$ with respect to $\|\cdot\|_{j\vee\1}$, as required.\medskip

Let $$H:=E\cap G=\ker\mu\cap E=\left\{f\in \ker\lambda, \lambda_{0}\left(f\right)=\sum\limits_{i=1}^{4}\lambda_{i}\left(f\right)\right\},$$ which is a full relatively uniformly closed subspace of $E$, since $\mu$ is a positive functional. Moreover, it is easy to see that $\1+3\1_{\Z\backslash\N}$ is a strong unit in $H$, and so the latter is directed.\medskip

Let us prove that for every $i=1,2,3,4$ the functional $\lambda_{i}$ is $E$-bounded on $H$. It is enough to show that this functional is bounded on $\left[\0,j\vee\1\right]_{H}$. For $f$ in that set note that $f\ge\0$, hence $\lambda_{i}\left(f\right)\ge 0$. On the other hand, $\left.f\right|_{\Z\backslash\N}\le \1$, and so $1\ge \lambda_{0}\left(f\right)=\sum\limits_{i=1}^{4}\lambda_{i}\left(f\right)$, implying $\lambda_{i}\left(f\right)\le 1$.

We claim that the space of $E$-bounded functionals on $H$ does not satisfy the RDP, and so cannot be a lattice. Let  $\nu_{1}:=\lambda_{1}+\lambda_{3}$ and $\nu_{2}:=\lambda_{1}+\lambda_{2}=\lambda_{3}+\lambda_{4}$. It is clear that $\lambda_{1},\lambda_{3}\le \nu_{1},\nu_{2}$.

For $i=1,2,3,4$ let $X_{i}:=\left\{4n+i\right\}_{n\in\N}\subset\Z$, and for $i=1,2$ and $j=3,4$ put $h_{ij}:=\1_{X_{i}\cup X_{j}}+2\1_{\Z\backslash\N}\in H$. Then, we have the following values of $\lambda_{1},\lambda_{3},\nu_{1},\nu_{2}$ at $h_{13},h_{14},h_{23},h_{24}$:
$$
\begin{array}{c|cccc}
 & \lambda_{1} & \lambda_{3} & \nu_{1} & \nu_{2} \\ \hline
h_{13} & 1 & 1 & 2 & 1 \\
h_{14} & 1 & 0 & 1 & 1\\
h_{23} & 0 & 1 & 1 & 1\\
h_{24} & 0 & 0 & 0 & 1
\end{array}
$$

Assume that a linear functional $\theta$ on $H$ satisfies $\lambda_{1},\lambda_{3}\le\theta\le \nu_{1},\nu_{2}$. These inequalities force $\theta\left(h_{13}\right)=\theta\left(h_{14}\right)=\theta\left(h_{23}\right)=1$ and $\theta\left(h_{24}\right)=0$, but then $1=\theta\left(h_{13}\right)+\theta\left(h_{24}\right)=\theta\left(\1+3\1_{\Z\backslash\N}\right)=\theta\left(h_{14}\right)+\theta\left(h_{23}\right)=2$, contradiction.
\qed\end{example}\medskip

Let us present an order bounded functional on an order dense, majorizing and relatively uniformly dense (in the sense of adherence) subspace of an Archimedean vector lattice, that cannot be extended to an order bounded functional on the whole space. This example shows that the additional hypothesis is not superfluous in Propositions \ref{ruex}.

\begin{example}\label{noext}
Let $I$ the set of cardinality continuum. Note that the collection of countable subsets of $I$ also has cardinality continuum, as well as the set of all countably supported functions on $I\times \N$. For $f:I\to\R$ and $g:\N\to\R$ we use the notation $f\otimes g$ for a function defined by $\left[f\otimes g\right]\left(i,n\right):=f\left(i\right)g\left(n\right)$.

Let $g:=\I_{\N}$, the identity map on $\N$, and for $i\in I$ let $g_{i}:=\delta_{i}\otimes g$. For $p>0$ we have $g_{i}^{p}=\delta_{i}\otimes g^{p}$. Let $G:=\spa\left\{g^{\frac{1}{n}}\right\}_{n\in\N}$ and for $i\in I$ let $G_{i}:=\delta_{i}\otimes G=\spa\left\{g_{i}^{\frac{1}{n}}\right\}_{n\in\N}$. It is not hard to show that if $\0\ne h\in G$, then either $\lim\limits_{n\to\8}h\left(n\right)=+\8$ or $\lim\limits_{n\to\8}h\left(n\right)=-\8$

We consider $c_{0}\left(I,\ell_{\8}\right)$ as a space of functions on $\R^{I\times \N}$, where it is an ideal. Other ideals that we consider are $c_{00}\left(I,\ell_{\8}\right)$ and $c_{0}\left(I\times\N\right)$. Let $E:=c_{0}\left(I,\ell_{\8}\right)+\R\1+\bigplus\limits_{i\in I}G_{i}$. 

We claim that the latter is a sublattice of $\R^{I\times \N}$. Every $e\in E$ is of the form $e=f+r\1+\sum\limits_{k=1}^{l}\delta_{i_{k}}\otimes h_{k}+\sum\limits_{k=l+1}^{m}\delta_{i_{k}}\otimes h_{k}$, where $f\in c_{0}\left(I,\ell_{\8}\right)$, $r\in \R$, $i_{1},...,i_{n}\in I$ and $h_{k}\in G$ are such that $h_{k}\left(n\right)\to +\8$ and if $k\le l$ and $h_{k}\left(n\right)\to-\8$, otherwise. Then, $e^{+}=r^{+}\1+\sum\limits_{k=1}^{l}\delta_{i_{k}}\otimes h_{k}+f^{+}+f_{0}$, where $f_{0}\in c_{0}\left(I,\ell_{\8}\right)$. Indeed, for $i\notin\left\{i_{1},...,i_{n}\right\}$ we have $f_{0}\left(i,n\right)=\left(f\left(i,n\right)+r\right)^{+}-f\left(i,n\right)^{+}-r^{+}$, hence $0\ge f_{0}\left(i,n\right)\ge -\left|f_{0}\left(i,n\right)\right|$; if $i=i_{k}$, for $k\le l$, we have $$f_{0}\left(i,n\right)=\left(r+f\left(i,n\right)+h_{k}\left(i,n\right)\right)^{+}-f\left(i,n\right)^{+}-r^{+}-h_{k}\left(n\right)=r+f\left(i,n\right)-f\left(i,n\right)^{+}-r^{+},$$ eventually, while if $i=i_{k}$, for $k> l$, we have $$f_{0}\left(i,n\right)=\left(r+f\left(i,n\right)+h_{k}\left(n\right)\right)^{+}-f\left(i,n\right)^{+}-r^{+}=-f\left(i,n\right)^{+}-r^{+},$$ eventually, in both cases $f_{0}\left(i,\cdot\right)\in\ell_{\8}$; combining three cases yields $f_{0}\in c_{0}\left(I,\ell_{\8}\right)$.

Note that elements of $c_{0}\left(I,\ell_{\8}\right)$ are countably supported functions, and there are continuum many such functions on $I\times\N$. Hence, it is possible to label $c_{0}\left(I,\ell_{\8}\right)=\left\{e_{i},~ i\in I\right\}$. For $i\in I$ and $k\in\N$ let $h_{i,k}:=e_{i}+2^{-k}g_{i}^{\frac{1}{k+1}}$, and let $$H:=\spa \left( \left\{\1\right\}\cup\left\{g_{i}\right\}_{i\in I}\cup\left\{h_{i,k}\right\}_{i\in I,~k\in\N}\right)+c_{00}\left(I\times\N\right).$$ Clearly, $H$ is majorizing in $E$, and since $c_{00}\left(I\times\N\right)\subset H$, the latter is order dense in $E$.

In order to show that $E$ is the relative uniform adherence of $H$ in $E$, it is enough to find a sequence in $H$ which relatively uniformly converges to every element of $c_{0}\left(I,\ell_{\8}\right)$, as well as to $g_{i}^{\frac{1}{k+1}}$, for every $i\in I$ and $k\in \N$. For every $e\in c_{0}\left(I,\ell_{\8}\right)$ there is $i\in I$ such that $e=e_{i}$, and then $H\ni h_{i,k}\toru e_{i}$. For every $i\in I$ and $k,m\in\N$ we have that $H\ni 2^{k}\left(h_{i,k}-h_{i,m}\right)=g_{i}^{\frac{1}{k+1}}-2^{k-m}g_{i}^{\frac{1}{m+1}}\xrightarrow[m\to\8]{\mathrm{ru}}g_{i}^{\frac{1}{k+1}}$.\medskip

Let $F:= \R\1+c_{0}\left(I\times\N\right)+ c_{00}\left(I,\ell_{\8}\right)+ \bigplus\limits_{i\in I}G_{i}$. Choose distinct $\left\{i_{k}\right\}_{k\in\N}$ and let $v:=\sum\limits_{k\in\N}\frac{1}{k}\delta_{i_{k}}\otimes \1$ (pointwise). Clearly, $v\in c_{0}\left(I,\ell_{\8}\right)\subset E$. Let us show that $v\notin F$. Assume the contrary. The elements of $F$ are constants outside of a countable subset of $I\times\N$, and for $v$ such constant is $0$. Hence, $v\in c_{0}\left(I\times\N\right)+ c_{00}\left(I\right)\otimes \left(G+\ell_{\8}\right)$, and so $v=f+\delta_{j_{1}}\otimes h_{1}+...+\delta_{j_{n}}\otimes h_{n}$, where $f\in c_{0}\left(I\times\N\right)$ and $h_{n}\in G+\ell_{\8}$, for some $j_{1},...,j_{n}\in I$. There is $k\in\N$ such that $i_{k}\ne j_{1},...,j_{n}$, and so for every $m\in\N$ we have $f\left(i_{k},m\right)=v\left(i_{k},m\right)=\frac{1}{k}$, contradicting $f\in c_{0}\left(I\times\N\right)$.

There is a linear (no boundedness assumptions) functional $\nu$ on $E$ which vanishes on $F$ and such that $\nu\left(v\right)=1$. Note that for every $i\in I$ and $k\in \N$ we have that $g_{i}^{\frac{1}{k+1}}\in G_{i}\subset F$, hence  $\nu\left(h_{i,k}\right)=\nu\left(e_{i}\right)$. Assume that $\left.\nu\right|_{H}$ has an order bounded (hence, relatively uniformly continuous) extension $\mu$ to $E$. Then, for every $i\in I$ we have $\mu\left(e_{i}\right)=\lim\limits_{k\to\8}\mu\left(h_{i,k}\right)=\lim\limits_{k\to\8}\nu\left(h_{i,k}\right)=\nu\left(e_{i}\right)$, hence $\mu$ and $\nu$ agree on $c_{0}\left(I,\ell_{\8}\right)$. For $n\in\N$ let $v_{n}:=\sum\limits_{k\ge n}\frac{1}{k}\delta_{i_{k}}\otimes \1$; note that $v_{n}\in c_{0}\left(I,\ell_{\8}\right)$. For every $n\in\N$ we have that $v-v_{n}\in c_{00}\left(I,\ell_{\8}\right)\subset F$, therefore $\nu\left(v-v_{n}\right)=0$, hence $nv_{n}\in\left[\0,\1\right]$ but $\mu\left(nv_{n}\right)=\nu\left(nv_{n}\right)=n\nu\left(v\right)=n$ shows that $\mu$ is not order bounded, contradiction.\medskip

Let us show however that $\left.\nu\right|_{H}$ is order bounded. First, an auxiliary fact. Let $V\subset \ell_{\8}\left(I\times\N\right)\cap E$ with $\dim V<\8$ and let $W:=V+c_{0}\left(I\times\N\right)$. Since $c_{0}\left(I\times\N\right)$ is a closed subspace of $\ell_{\8}\left(I\times\N\right)$, it is a closed subspace of $W$ of finite codimension. We claim that $\nu$ it is norm bounded on $W$. Indeed, since $c_{0}\left(I\times\N\right)\subset F\subset\ker\nu$, it follows that $\ker\nu\cap W$ is a closed subspace of $W$, and so $\left.\nu\right|_{W}\in W^{*}$.

In order to prove that $\left.\nu\right|_{H}$ is order bounded, it is enough to show that $\nu$ is bounded on $\left[-h,h\right]\cap H$, for every $h$ of the form $\1+g_{j_{1}}+...+g_{j_{m}}$, where $j_{1},...,j_{m}\in I$. Let $J:=\left\{j_{1},...,j_{m}\right\}$.

Assume that $e\in H$ is such that $\left|e\right|\le h$. We have that $e=r\1+\sum r_{i}g_{i}+ \sum r_{i,k}h_{i,k} + f$, where $f\in c_{00}\left(I\times\N\right)$. Fix $i\notin J$. For every $n\in\N$ we have that $$1=h\left(i,n\right)\ge \left|e\left(i,n\right)\right|=\left|f\left(i,n\right)+r+r_{i}n+\sum\limits_{k} r_{i,k}2^{-k}n^{\frac{1}{k+1}}+\sum\limits_{j,k} r_{j,k} e_{j}\left(i,n\right)\right|.$$
Note that $f\left(i,\cdot\right)+r\1+\sum\limits_{j,k} r_{j,k} e_{j}\left(i,\cdot\right)\in\ell_{\8}$. Let $s$ be the norm of that function there. We have $1+s\ge \left|r_{i}n+\sum\limits_{k} r_{i,k}2^{-k}n^{\frac{1}{k+1}}\right|$, for every $n\in\N$. Clearly, this is possible only if $r_{i}=r_{i,k}=0$, for every $k\in\N$.

It follows that $\left[-h,h\right]\cap H\subset c_{00}\left(I\times\N\right)+\spa\left\{\1,e_{j_{1}},...,e_{j_{m}}\right\}+ \bigplus\limits_{i\in J}G_{i}$.

Let $P,Q:\R^{I\times\N}\to \R^{I\times\N}$ be defined by $Pf:=\1_{J\times \N}f$, and $Q:=\I_{\R^{I\times\N}}-P$. We claim that $PE\subset F$ and $QE\subset E$. Since $Q$ is a band projection, it preserves ideals of $\R^{I\times\N}$ such as $c_{0}\left(I,\ell_{\8}\right)$, $c_{00}\left(I,\ell_{\8}\right)$, $c_{0}\left(I\times\N\right)$, $c_{00}\left(I\times\N\right)$. Also, $P\ell_{\8}\left(I\times\N\right)\subset c_{00}\left(I,\ell_{\8}\right)\subset F$. Hence,  $Qc_{0}\left(I,\ell_{\8}\right)\subset E$, $Pc_{0}\left(I,\ell_{\8}\right)\subset  F$,  $P\1\in   F\subset E$, $Q\1=\1-P\1\in E$. Finally, for every $i\in I$ and $h\in G_{i}$ we have $P$ either preserves or nullifies $h$; in both cases $Ph,Qh\in F$.

Let $W:=V+c_{0}\left(I\times\N\right)$, where $V=Q\spa\left\{\1,e_{j_{1}},...,e_{j_{m}}\right\}$. According to the auxiliary fact above, we have $\left.\nu\right|_{W}\in W^{*}$. Let $t$ be the norm of $\nu$ on $W$. Note that $$Q\left(\left[-h,h\right]\cap H\right)\subset Q\left(c_{00}\left(I\times\N\right)+\spa\left\{\1,e_{j_{1}},...,e_{j_{m}}\right\}+ \bigplus\limits_{i\in J}G_{i}\right)\subset W,$$ and so $Q\left(\left[-h,h\right]\cap H\right)\subset Q\left[-h,h\right]\cap W\subset \left[-\1,\1\right]\cap W=\Ba_{W}$. Since $PE\subset\ker\nu$, it follows that $\nu\left(e\right)=\nu\left(Qe\right)$, for every $e\in E$. Thus, for $f\in \left[-h,h\right]\cap H$ we have $Qf\in \Ba_{W}$, and so $\left|\nu\left(f\right)\right|=\left|\nu\left(Qf\right)\right|\le t$.
\qed\end{example}\medskip
Next, a bounded supremum continuous operator which fails to be a pre-homomorphism.

\begin{example}\label{bscnr}
Let $\Po\left(\R\right)$ be the space of polynomials on $\R$ with the pointwise order, and let $E:=\left\{f\in\Po\left(\R\right),~ f\left(-1\right)+f\left(1\right)=2f\left(0\right)\right\}$. As a subspace of $\Co\left(\R\right)$ it is Archimedean. To show that $E$ is directed, let $e\in E$ and let $n>\deg e$ be an even natural number. Let $h:\R\to\R$ be defined by $h\left(t\right):=t^{n}-t^{2}+1$. Then, $h\left(-1\right)=h\left(0\right)=h\left(1\right)=1$ implies that $h\in E$. If $\left|t\right|\ge 1$ we have $h\left(t\right)\ge 1$. If $\left|t\right|< 1$, then $t^{2}<1$, and so $h\left(t\right)> 0$. Hence, $h\in E_{+}$. Moreover, since $h$ does not vanish and $\deg h>\deg e$, it follows that $\frac{e}{h}\in\Co_{0}\left(\R\right)$, hence it is bounded by some $r\ge 0$. We conclude that $-rh\le e\le rh$, and so $E$ is directed.

Next, note that for any $e,f\in \Po\left(\R\right)$ with $-f\le e\le f$, it follows that $\deg e\le \deg f$. Indeed, otherwise we have $\lim\limits_{t\to\8}\frac{f\left(t\right)}{e\left(t\right)}=0$, but $e\ne \0$ and $\left|e\right|\le \left|f\right|$ leads to $\limsup\limits_{t\to\8}\frac{\left|f\left(t\right)\right|}{\left|e\left(t\right)\right|}\ge 1$, contradiction.

Hence, for every $f\in E_{+}$ we have that $\dim E_{f}\le \deg f+1$. Therefore, $\left[\0,f\right]_{E}$ is a closed bounded set in a finite dimensional normed space $\left(E,\|\cdot\|_{f}\right)$, hence it is compact. As the embedding of $E_{f}$ into $\Co\left(\R\right)$ is continuous, it follows that $\left[\0,f\right]_{E}$ is compact in $\Co\left(\R\right)$, hence equicontinuous.\medskip

Let us show that $\delta_{0}$ is bounded supremum continuous. Let $f\in E_{+}$ and let $P\subset \left[\0,f\right]_{E}$ be such that $\bigwedge P=\0$. Assume that there is $\varepsilon\in\left(0,1\right)$ such that $p\left(0\right)\ge\varepsilon$, for every $p\in P$. Since $\left[\0,f\right]_{E}$ is equicontinuous, there is $\eta>0$ such that $p\left(t\right)\ge\frac{\varepsilon}{2}$, for every $p\in P$ and $\left|t\right|\le \eta$. Let $h:\R\to\R$ be defined by $h\left(t\right):=t^{2}\left(t^{2}-1\right)^{2}\left(\eta^{2}-t^{2}\right)$. Since $h\left(\pm 1\right)=0=h\left(0\right)$, it follows that $h\in E$, but $h\not\le\0$, as $h\left(\frac{\eta}{2}\right)>0$. Also, $h$ is bounded on $\left[-\eta,\eta\right]$ by some $s>0$. For every $p\in P$ and $\left|t\right|\ge\eta$ we have $p\left(t\right)\ge 0\ge \frac{\varepsilon}{2s}h\left(t\right)$; if $\left|t\right|\le\eta$, then $p\left(t\right)\ge \frac{\varepsilon}{2}\ge \frac{\varepsilon}{2s}h\left(t\right)$. We conclude that $\frac{\varepsilon}{2s}h\not\le \0$ is a lower bound for $P$, contradiction.\medskip

Let $\nu:E\to\R$ be defined by $\nu:=2\delta_{1}+\delta_{-1}$. Note that $0_{E'}\le\nu\le 4\delta_{0}$, and so it follows from Proposition \ref{bsc0} that $\nu$ is bounded supremum continuous. It is left to exhibit a failure of $\nu$ to be a pre-homomorphism. Let $e,f:\R\to \R$ be defined by $e\left(t\right):=\left(\frac{3}{5}-t\right)\left(5t^{2}+3t+2\right)$, $f\left(t\right):=\left(t-\frac{2}{5}\right)\left(-\frac{5}{4}t^{2}-\frac{1}{2}t+\frac{11}{4}\right)$ and $g\left(t\right):=t^{2}\left(t^{2}-1\right)^{2}$. Note that $g\ge\0$. We have the following evaluations of $\delta_{-1},\delta_{0},\delta_{1},\nu$ at $e,f,g$:
$$
\begin{array}{c|cccc}
 & \delta_{-1} & \delta_0 & \delta_1 & \nu \\ \hline
e & \dfrac{32}{5} & \dfrac65 & -4 & -\dfrac85 \\[3mm]
f & -\dfrac{14}{5} & -\dfrac{11}{10} & \dfrac35 & -\dfrac85\\[3mm]
g & 0 & 0 & 0 & 0
\end{array}
$$
Hence, $e,f,g\in E$. The zeros of $g$ are precisely $\left\{-1,0,1\right\}$, and so $g\in E_{+}$. Since $\deg g=6> 3=\deg e=\deg f$, arguing as above, we see that there is $t_{0}\ge 1$ such that $\left|e\left(t\right)\right|\le g\left(t\right)$, whenever $\left|t\right|\ge t_{0}$. Since $e\left(-1\right),e\left(0\right),f\left(1\right)>0$, there is $\varepsilon>0$ such that $e$ is positive on the $\varepsilon$-neighborhoods of $-1$ and $0$, while $f$ is positive on the $\varepsilon$-neighborhood of $1$. Let $U$ be the union of these three open sets. On the compact sets $\left[-t_{0},t_{0}\right]\backslash U$ the continuous functions $e,f,g$ are bounded below by some $u,v,w\in\R$, where $w>0$. Let $r\ge 1$ be such that $u+rw,v+rw\ge 0$. Then, $rg+e\vee f\in\Co\left(\R\right)_{+}$. Indeed, if $t\in U$, then $rg\left(t\right)+e\left(t\right)\vee f\left(t\right)\ge e\left(t\right)\vee f\left(t\right)\ge 0$; if $t\in \left[-t_{0},t_{0}\right]\backslash U$, then $rg\left(t\right)+e\left(t\right)\vee f\left(t\right)\ge w+u\vee v\ge 0$; if $\left|t\right|\ge t_{0}$, then $rg\left(t\right)+e\left(t\right)\vee f\left(t\right)\ge g\left(t\right)-g\left(t\right)=0$. Hence, $\left\{rg+e,rg+f\right\}^{\Uparrow}\subset E_{+}$, and so $\0\in \left\{rg+e,rg+f\right\}^{\Uparrow\Downarrow}$. At the same time, $\nu\left(rg+e\right)=\nu\left(rg+f\right)=-\frac{8}{5}$, and so the definition of a pre-homomorphism is not satisfied for $\nu$.
\qed\end{example}\medskip

A bounded supremum continuous functional on an order unit space which fails to be order continuous.

\begin{example}\label{bscnoc}
Let $E:=\R^{3}_{\8}\oplus_{\8} c_{0}\oplus_{\8} c_{0}$. For $n\in\N$ define $\nu_{n},\mu_{n},\lambda_{n}:E\to\R$ by $\nu_{n}\left(x,y,z,f,g\right):=x+f\left(n\right)$, $\mu_{n}\left(x,y,z,f,g\right):=y+g\left(n\right)$ and $\lambda_{n}\left(x,y,z,f,g\right):=\frac{x+y}{2}+\frac{z-f\left(n\right)-g\left(n\right)}{n+1}$, and $\rho:E\to \R$ by $\rho\left(x,y,z,f,g\right):=\frac{x+y}{2}+z$. Let $C:=\left\{\nu_{n},\mu_{n},\lambda_{n}\right\}_{n\in\N}\cup\left\{\rho\right\}\subset 3\Ba_{E^{*}}$. Note that $C_{\top}=\left\{0_{E}\right\}$. Indeed, $\nu_{n}\left(x,y,z,f,g\right)=0$ $\Rightarrow$ $-x=f\left(n\right)\to 0$, hence $x=0$, and so $f\left(n\right)=0$, for every $n\in\N$; similarly $y=0$ and $g=\0$; finally $\rho\left(x,y,z,f,g\right)=0$ implies that $z=0$.

Let $E_{+}:=\left\{e\in E,~ \forall c\in C:~ c\left(e\right)\ge 0\right\}$. Clearly, $E_{+}$ is a closed convex wedge in $E$, and the fact that $C_{\top}=\left\{0_{E}\right\}$ implies that this wedge is a cone. Hence, $E$ ordered via $E_{+}$ is an Archimedean ordered vector space and $C\subset E^{*}_{+}$. We claim that $u:=\left(1,1,0,\0,\0\right)$ is a strong unit. Indeed, we have $c\left(u\right)=1$, for every $c\in C$; also $\Ba_{E}\subset\left[-3u,3u\right]$, for if $e\in \Ba_{E}$, then $\left|c\left(e\right)\right|\le 3$, hence $c\left(3u\pm e\right)\ge 0$, for every $c\in C$.

Let us show that on the other hand, $\left[-u,u\right]\subset 2\Ba_{E}$. Suppose $e=\left(x,y,z,f,g\right)\in E$ is such that $-u\le e\le u$, or, equivalently, that $\left|c\left(e\right)\right|\le 1$, for every $c\in C$. Then $1\ge \left|\nu_{n}\left(e\right)\right|\ge \left|x\right|-\left|f\left(n\right)\right|\to \left|x\right|$ implies $\left|x\right|\le 1$, an then $1\ge\left|f\left(n\right)\right|-\left|x\right|$ implies $\left|f\left(n\right)\right|\le 2$, for every $n\in\N$. Similarly, $\left|y\right|\le 1$, and $\left|g\left(n\right)\right|\le 2$, for every $n\in\N$. Finally, $\left|z\right|\le \left|\rho\left(e\right)\right|+\frac{\left|x\right|+\left|y\right|}{2}\le 2$. We conclude that $\|\cdot\|_{u}\simeq \|\cdot\|$, and so $E$ is a complete order unit space.\medskip

Fix $n\in\N$. Let us show that $\lambda_{n}$ is supremum continuous. Suppose that $P\subset E$ is such that $\bigwedge P=0_{E}$, but $\lambda_{n}\left(P\right)\ge\varepsilon$, for some $\varepsilon>0$. Let $f\in c_{0}$ be the $n$-th element of the standard basis, and let $e:=\left(1,1,0,\left(n+1\right)f,\left(n+1\right)f\right)\in E$. Then, for every $\lambda_{n}\ne c\in C$ we have $c\left(e\right)\ge 0$, but $\lambda_{n}\left(e\right)=-1<0$, and so $e\notin E_{+}$. We claim that $-\varepsilon e\le P$, which would contradict $\bigwedge P=0_{E}$. Indeed, for every $\lambda_{n}\ne c\in C$ and $p\in P$ we have $c\left(-\varepsilon e\right)\le 0\le  c\left(p\right)$, while $\lambda_{n}\left(-\varepsilon e\right)=\varepsilon\le c\left(p\right)$.\medskip

Next, we show that the projections $\varphi,\psi$ onto the first and the second coordinates are bounded supremum continuous. First, if $e=\left(x,y,z,f,g\right)\in E_{+}$, then for every $n\in\N$ we have $0\le x+f\left(n\right)\to x$, and so $\varphi\in E^{*}_{+}$. Similarly, $\psi\in E^{*}_{+}$. Moreover, $2\lambda_{n}\to \varphi+\psi$ (in norm), and since all involved elements are positive, there is a subsequence $\left(2\lambda_{n_{m}}\right)_{m\in\N}$ which converges to $\varphi+\psi$ relatively uniformly (the proof of Proposition \ref{ru1} goes through, since the elements of the net are positive). By part (ii) of Proposition \ref{preserve} we conclude that $\varphi+\psi$ is bounded supremum continuous, and since $0_{E^{*}}\le \varphi\le \varphi+\psi$, it follows from part (i) of Proposition \ref{preserve} that $\varphi$ is bounded supremum continuous.\medskip

Let us prove however, that $\varphi$ is not order continuous. For $n\in\N$ let $f_{n}\in c_{0}$ be the $n$-th element of the standard basis, and let $g_{n}:=\sum\limits_{k=1}^{n}f_{k}$, $w:=\left(2,2,-2,\0,\0\right)$, $p_{n}:=2\left(1,1,0,-f_{n},\0\right)$, $q_{n}:=2\left(1,1,0,\0,-f_{n}\right)$, $r_{n}:=\left(2,2,0,\left(n+1\right)f_{n},\left(n+1\right)f_{n}\right)$, $v_{n}:=\left(-1,1,0,g_{n},-g_{n}\right)$.

Let $Q:=\left\{p_{n},q_{n},r_{n}\right\}_{n\in\N}\cup\left\{w\right\}$. Note that $\bigwedge c\left(Q\right)=0$, for every $c\in C$, which is easily seen to imply $\bigwedge Q=0_{E}$. Let us show that $v_{n}\xrightarrow[]{Q}0_{E}$. We have the following values of our functionals at our vectors:
$$
\begin{array}{c|ccccc}
 & p_m & q_m & r_m & w & v_m \\ \hline
\nu_n
 & 2-2\delta_{nm}
 & 2
 & 2+(m+1)\delta_{nm}
 & 2
 & -1+\1_{\{n\le m\}} \\[2pt]
\mu_n
 & 2
 & 2-2\delta_{nm}
 & 2+(m+1)\delta_{nm}
 & 2
 & 1-\1_{\{n\le m\}} \\[2pt]
\lambda_n
 & 2+\dfrac{2\delta_{nm}}{n+1}
 & 2+\dfrac{2\delta_{nm}}{n+1}
 & 2-2\delta_{nm}
 & 2-\dfrac{2}{n+1}
 & 0 \\[6pt]
\rho
 & 2 & 2 & 2 & 0 & 0
\end{array}
$$
It now follows that $-r_{m}\le v_{k}\le r_{m}$ and $-w\le v_{k}\le w$, for every $m,k\in\N$, and $-p_{m}\le v_{k}\le p_{m}$ and $-q_{m}\le v_{k}\le q_{m}$, for every $k\ge m$. However, $\varphi\left(v_{k}\right)=-1\not\to 0$ shows that $\varphi$ is not order continuous.
\qed\end{example}\medskip

A bounded supremum continuous operator from an order unit space with the RDP into an Archimedean vector lattice which fails to be order continuous.

\begin{example}\label{tong}
In this example we denote $\ell_{\8}:=\ell_{\8}\left(\Z\backslash\left\{0\right\}\right)$. Enumerate $X:=\Q\cap\left[0,1\right]:=\left\{x_{n}\right\}_{n\in\N}$. Let $T:\ell_{\8}\to \ell_{\8}\left(X\right)$ be defined by $\left[Tf\right]\left(x_{n}\right):=\frac{n}{n+1}f_{n}+\frac{1}{n+1}f_{-n}$, $n\in\N$, $f\in \ell_{\8}$. Clearly, $T$ is a positive operator.

Let $F$ be the set of elements of $\ell_{\8}\left(X\right)$ which extend to continuous functions on $\left[0,1\right]$ (in other words, uniformly continuous functions on $X$). Clearly, $E:=T^{-1}F$ is a linear subspace of $\ell_{\8}$ which contains $\1$. We also define $S:\Co\left[0,1\right]\to E$ by $\left(Sf\right)_{\pm n}:=f\left(x_{n}\right)$, $n\in\N$. Note that $TSf$ is the restriction of $f$ to $X$.

Let us show that $E$ has the RDP. Let $e,f,g,h\in E$ be such that $e,f\le g,h$. Without loss of generality we may assume that $-\1\le e,f,g,h\le \1$. Since $T$ is positive, we have $u:=T\left(g\wedge h\right)\ge T\left(e\vee f\right)=:v$.

There are $e',f',g',h'\in\Co\left[0,1\right]$, whose restriction to $X$ equals $Te,Tf,Tg,Th$, respectively. Note that since $e',f'\le g',h'$ on $X$, we also have $e',f'\le g',h'$ on $\left[0,1\right]$, and so $e'\vee f'\le g'\wedge h'$.

Fix $n\in\N$. We have  $u\left(x_{n}\right)=\frac{n}{n+1}\left(g_{n}\wedge h_{n}\right)+\frac{1}{n+1}\left(g_{-n}\wedge h_{-n}\right)$, whereas $\left[g'\wedge h'\right]\left(x_{n}\right)=\left(\frac{n}{n+1}g_{n}+\frac{1}{n+1}g_{-n}\right)\wedge \left(\frac{n}{n+1}h_{n}+\frac{1}{n+1}h_{-n}\right)$, and so $$a:=\left[g'\wedge h'\right]\left(x_{n}\right)-u\left(x_{n}\right)=\frac{1}{n+1}\left(n\left(g_{n}-h_{n}\right)^{+}+\left(g_{-n}-h_{-n}\right)^{+}\right)\wedge \left(n\left(h_{n}-g_{n}\right)^{+}+\left(h_{-n}-g_{-n}\right)^{+}\right),$$
Hence, applying the symmetric form of \cite[Theorem 1.7(4)]{ab0} yields
$$0\le a\le \frac{1}{n+1}\left(\left(h_{-n}-g_{-n}\right)^{+}+\left(g_{-n}-h_{-n}\right)^{+}\right)=\frac{1}{n+1}\left|g_{-n}-h_{-n}\right|\le \frac{2}{n+1}.$$
Consider the function $w:=g'\wedge h'-u$ on $X$, and extend it by $0$ on $\left[0,1\right]\backslash X$. This extended function is positive on $\left[0,1\right]$, and for every $r>0$, we have $w^{-1}\left[r,+\8\right)\subset \left\{q_{n}\right\}_{\frac{2}{n+1}\ge r}$ is finite, hence closed. It follows that $w$ is upper semicontinuous, and so $u':=g'\wedge h'-w$ is lower semicontinuous, agreeing with $u$ on $X$ and with $g'\wedge h'$ on $\left[0,1\right]\backslash X$. Similarly, we can define an upper semicontinuous $v'$, which agrees with $v$ on $X$ and with $e'\vee f'$ on $\left[0,1\right]\backslash X$. In particular, $v'\le u'$: if $x\in X$, then $v'\left(x\right)=v\left(x\right)\le u\left(x\right)=u'\left(x\right)$, and the case when $x\notin X$ is similar.

By Katetov-Tong theorem (see \cite[Exercise 1.17.15(b)]{engelking}) there is $j\in \Co\left[0,1\right]$ such that $v'\le j\le u'$. For every $n\in\N$ we have $$\frac{n}{n+1}\left(e_{n}\vee f_{n}\right)+\frac{1}{n+1}\left(e_{-n}\vee f_{-n}\right)\le j\left(x_{n}\right)\le \frac{n}{n+1}\left(g_{n}\wedge h_{n}\right)+\frac{1}{n+1}\left(g_{-n}\wedge h_{-n}\right),$$ and it is easy to find $i\in\ell_{\8}$ such that
$j\left(x_{n}\right)=\frac{n}{n+1}i_{n}+\frac{1}{n+1}i_{-n}$, for every $n\in\N$ (implying $i\in E$), and $e\vee f\le i\le g\wedge h$, as required.\medskip

Next, we show that the embedding of $E$ into $\ell_{\8}$ is supremum continuous. Let $Q\subset E$ be such that $\bigwedge Q=\0$. Assume that $e\in\ell_{\8}$ is a lower bound of $Q$, and then so is $e^{+}$. If there is $n\in\N$ such that $r:=e_{n}>0$, then $u:=r\1_{\left\{n\right\}}-rn\1_{\left\{-n\right\}}\in\ker T\subset E$ satisfies $\0\not\ge u\le e^{+}\le Q$, contradiction. If there is $n\in\N$ such that $t:=e_{-n}>0$, then a contradiction is reached by considering $v:=t\1_{\left\{-n\right\}}-\frac{t}{n}\1_{\left\{n\right\}}$. Hence, $e\le\0$.\medskip

We now prove that $\left.T\right|_{E}$ is bounded supremum continuous. Let $Q\subset E$ be order bounded and such that $\bigwedge Q=\0$. Without loss of generality we may assume that $Q\subset \Ba_{E}$, and the last step implies that $\bigwedge\limits_{q\in Q}q_{\pm n}=0$, for every $n\in\N$. Assume that $f\le TQ$ and fix $m\in\N$. For every $n\ge m$ and $q\in Q$ we have $f\left(x_{n}\right)\le \frac{n}{n+1}q_{n}+\frac{1}{n+1}q_{-n}\le q_{n}+\frac{1}{m}$. Hence $f\left(x_{n}\right)\le\frac{1}{m}$, for every $n\ge m$. Since $X\backslash\left\{x_{1},...,x_{m-1}\right\}$ is dense in $\left[0,1\right]$, it follows that $f\le \frac{1}{m}\1$, and since $m$ was arbitrary we get $f\le \0$. We conclude that $T$ is bounded supremum continuous.\medskip

Finally, we show that $\left.T\right|_{E}$ fails to be order continuous. Let $\Lambda$ be the collection of finite subsets of $\N$ and let $\Gamma:=\Lambda\times \N$. For $\gamma=\left(A,m\right)$, let $f_{\gamma}:=\1\wedge md\left(\cdot,A'\right)\in \Co\left[0,1\right]$, where $A'=\left\{x_{n}\right\}_{n\in A}$. Note that for every $n\in\N$ and $A\ni n$ and $m\in\N$, we have that $f_{\gamma}$ vanishes on $x_{n}$ and $0\le f_{\gamma}\left(x_{m}\right)\le 1$. We claim that  $f_{\gamma}\not\too \0$ in $F$, because for every $\gamma_{0}=\left(A,m\right)$ we have $\bigvee\limits_{\gamma\ge\gamma_{0}} f_{\gamma}=\1$. To prove the last assertion assume that $g\in F$ is such that $g\ge f_{\gamma}$, for every $\gamma\ge\gamma_{0}$. Then, for every $n\ge m$ we have $\gamma:=\left(A,n\right)\ge \gamma_{0}$, and so $g\left(x\right)\ge \1\wedge nd\left(x,A'\right)$. It follows that $g\left(x\right)\ge 1$, for every $x$ in the dense set $X\backslash A'$, hence $g\ge\1$.

For $n\in\N$ let $p_{\pm n}\in T^{-1}\1\subset E$ be defined by $p_{n}=\1_{\Z\backslash\left\{0,n\right\}} + n\1_{\left\{-n\right\}}$, $p_{-n}=\1_{\Z\backslash\left\{0,-n\right\}} + \frac{1}{n}\1_{\left\{n\right\}}$. It follows that $\bigwedge\limits_{n\in \pm\N} p_{n}=\0$. Also, for every $n,m\in\N$ and $A\ni n$ it follows $\0\le Sf_{\gamma}\le p_{\pm n}$, where $\gamma=\left(A,m\right)$. Therefore, $Sf_{\gamma}\too \0$, witnessed by $Q:=\left\{p_{n}\right\}_{n\in \pm\N}$. However, by the preceding step we have $f_{\gamma}=TSf_{\gamma}\not\too \0$ in $F$, hence $\left.T\right|_{E}$ is not order continuous.
\qed\end{example}

\subsection*{Acknowledgements}

The author would like to thank Koos Grobler, Anke Kalauch, Vladimir Troitsky and Onno van Gaans for helpful conversations on the topic of the paper. The author is also grateful to Marten Wortel for proofreading and Janko Stennder for contributing Example \ref{janko} and informing the author about \cite[Example 26]{kvgs}. The idea of using direct limits in order to construct $E^{\mathrm{ru}}$ in Proposition \ref{limit} has originated with Marcel de Jeu. Additionally David Mu\~noz-Lahoz gets credit for setting up the interface that enabled interactions with models GPT-5.5 and GPT-5.6 Sol, which generated ideas for Lemma \ref{rdpa}, Propositions \ref{rdpb}, \ref{rdp}, \ref{aip} and \ref{adpou} and Examples \ref{mo}, \ref{mo2}, \ref{nadpou}, \ref{odensen}, \ref{onuor}, \ref{onuo} and \ref{norucomp}, as well as all examples contained in Section \ref{sbce}.

A part of the work on the article, as well as the aforementioned conversations, took place during the author's attendance of COSA 2025 conference in Kruger Park in South Africa. The author is grateful to organizers for making it happen. The author's trip to South Africa was supported financially by the DSTI-NRF Centre of Excellence in Mathematical and Statistical Sciences (CoE-MaSS), South Africa. He is grateful to Jan Harm van der Walt for securing the funds. He was also supported by Pacific Institute for the Mathematical Sciences and grant PCI2024-155094-2 funded by MICIU/AEI/10.13039/50110001103

\end{document}